\documentclass[11pt]{amsart}
\usepackage{amsmath}
\usepackage[dvips]{graphicx}
\usepackage{amsfonts}
\usepackage{amsopn}
\usepackage{amsthm} 
\usepackage{amssymb}
\usepackage{color}
\usepackage{hyperref}

\usepackage{latexsym}
\usepackage{amsgen}

\newtheorem{theorem}{Theorem}[section]
\newtheorem{definition}[theorem]{Definition}
\newtheorem{lemma}[theorem]{Lemma}
\newtheorem{corollary}[theorem]{Corollary}
\newtheorem{proposition}[theorem]{Proposition}

\newtheorem{remark}[theorem]{Remark}

\numberwithin{equation}{section}

\newcommand{\Tr}{\rm Tr}

\newcommand{\Reals}{\mathbb{R}}

\newcommand{\Sk}{\mathcal{ {\rm Sym[\Mb^k]}} }
\newcommand{\Skt}{\mathcal{ {\rm Sym_2[\Mb^k]}} }

\newcommand{\Skcomplex}{\mathcal{ {\rm Sym[\mathbb C^k]}} }
\newcommand{\Sktcomplex}{\mathcal{ {\rm Sym_2[\mathbb C^k]}} }

\newcommand{\cSk}{ \mathcal{S}ym[k](\Reals)}
\newcommand{\cSkt}{ \mathcal{S}ym_2[k](\Reals)}
\newcommand{\cSkcomplex}{ \mathcal{S}ym[k](\mathbb C)}

\newcommand{\Mb}{\mathbb M}

\newcommand{\calPtwod}{\mathcal P_2(\Mb)}

\newcommand{\id}{{\bf id}}

\usepackage{amssymb}

\title[A partial Laplacian as an infinitesimal generator on the Wasserstein space]{A partial Laplacian as an infinitesimal generator on the Wasserstein space}
\author[Y.T. Chow]{Yat Tin~Chow}
\address{Department of Mathematics, UCLA, Los Angeles, CA 90095} \email{ytchow@math.ucla.edu}
\author[W. Gangbo]{Wilfrid~Gangbo}
\address{Department of Mathematics, UCLA, Los Angeles, CA 90095} \email{wgangbo@math.ucla.edu}

\begin{document} 
\begin{abstract}  We study stochastic processes on the Wasserstein space, together with their infinitesimal generators. 
One of these processes is modeled after Brownian motion and plays a central role in our work. Its infinitesimal generator defines a partial Laplacian on the space of Borel probability measures, and we use it to define heat flow on the Wasserstein space.   We verify a distinctive smoothing effect of this flow for a particular class of initial conditions. To this end, we will develop a theory of Fourier analysis and conic surfaces in metric spaces.  We note that the use of the infinitesimal generators has been instrumental in proving various theorems for Mean Field Games, and we anticipate they will play a key role in future studies of viscosity solutions of PDEs in the Wasserstein space.
\end{abstract}

\maketitle 
\tableofcontents

\section{Introduction} 
A fundamental result in stochastic analysis is that the Laplace operator is the infinitesimal generator of Brownian motion. That is,  for any twice continuously differentiable function $f:\mathbb{R}^d\rightarrow \mathbb{R}$ with bounded second derivatives,
$$
 \Tr {\rm (Hess \;  f)}(x)  =\lim_{t\rightarrow 0^+}\frac{\mathbb{E}f(x+\sqrt{2}W_t)-f(x)}{t}
$$
for each $x\in \mathbb{R}^d$.  Here, of course, $(W_t,t\ge 0)$ is a standard $d$--dimensional Brownian motion and $\Tr {\rm (Hess \;  f)} \equiv \Delta f$ 
is the Laplace operator.  This property is also closely related to the fact that 
$$
v(x,t)=\mathbb{E}f(x+\sqrt{2}W_t)
$$
is a solution in $(0,\infty) \times \mathbb{R}^d,$ of the heat equation 
\begin{equation}\label{ClassicalHeatFlow}
\partial_t v=\Delta v.
\end{equation}  

\par In this manuscript, we will lift these results from Euclidean space to the the space of Borel probability measures on $\mathbb{R}^d$ equipped with the Wasserstein metric. We will denote this space as $\mathcal P_2(\mathbb R^d)$.  Our motivation for pursing this line of research is to look for ways of regularizing partial differential equations in $\mathcal P_2(\mathbb R^d)$ such as Master Equations from the theory of mean field games. One example of such a Master Equation is 
\begin{align}\label{MasterEq}
&\partial_t{\mathcal V} +\int_{\mathbb R^d}\nabla_x{\mathcal V}\cdot \nabla_\mu {\mathcal V} \mu(dy) + \frac{1}{2}|\nabla_q {\mathcal V}|^2 +{\mathcal F}(q,\mu) 
-\int_{\mathbb R^d}\text{Tr}\left[\nabla^2_\mu {\mathcal V}\right]\mu(dy)\mu(dy')  \nonumber\\
&= 2\Delta_q{\mathcal V}  +2 \int_{\mathbb R^d}\Bigl( 
\nabla_y\cdot\left[\nabla_\mu {\mathcal V}(t,q,\mu)(y)\right] +\nabla_q\cdot\left[\nabla_\mu {\mathcal V}(t,q,\mu)\right](y) \Bigr) \mu(dy), &  
\end{align} 
which appeared in a recent work on Mean Fields Games \cite{CardaliaguetDLL}. 
This equation poses more challenges than we are currently prepared to address. So for now we address the smoothing effects (cf. Remark \ref{re:superposition}) of two of the three mechanisms at work in \eqref{MasterEq}, leaving out the mechanism induced by the underlying Lagrangian.

In what follows, we will define Brownian motion on the Wasserstein space as the stochastic process 
$$t \mapsto \mathbb B_t^{m}:=(\id+\sqrt 2 W_t)_\# m.$$
Here and below, $T_\# \mu$ is the push-forward measure defined for every Borel probability measure $\mu$ on $\mathbb R^d$ and Borel map $T:\mathbb R^d\rightarrow \mathbb R^d$ via the formula
$$
T_\# \mu(A):=\mu(T^{-1}(A)),
$$
After investigating properties of this new type of stochastic process, we will set   
\[
V(t, m)=\mathbb E\Bigl(U\bigl(\mathbb B_t^{m}\bigr) \Bigr)
\]
for $(t, m) \in (0,\infty) \times \calPtwod$ and assert    
\[
\partial_t V= \triangle_w V, \qquad V(0,\cdot)=U.
\]
We will take great care and precision in defining $\triangle_w V$ below in \eqref{eq:defn-Laplacian}. See Definition also \ref{defn:Laplacian}/Remark \ref{re:Laplacian2}.




In this manuscript, we start a Fourier analysis on the set of probability measures of finite second moments, the so--called the Wasserstein space. We later introduce  measures on the infinite dimensional metric space, which allow us to integrate by parts products of special functions defined on the Wasserstein space. Before introducing these measures, we devote some time to the study of the pairing of stochastic paths $(\mathbb B^m)_m$ on the Wasserstein $\mathcal P_2(\mathbb R^d)$, with a linear operator $\triangle_w$ defined on the set of twice differential functions on $\mathcal P_2(\mathbb R^d)$. Our pairing $(\mathbb B^m, \triangle_w)$ is reminiscent to the finite dimensional pairing $\bigl( (W+q)_{q \in \mathbb R^d}, \triangle\bigr)$, of the standard $d$--dimensional Brownian motions starting at $q \in \mathbb R^d$ with $\triangle,$ the Laplacian operator on $\mathbb R^d.$ The infinite dimensional linear operator $\triangle_w$ which is shown to be the infinitesimal generator of $\mathbb B$ is obtained as a partial  trace of the Hessian operator: If $U:\mathcal P_2(\mathbb R^d) \rightarrow \Reals$ is bounded and twice differentiable at $m$ and ${\rm Hess} \, U[m]$ denotes the Wasserstein Hessian of $U$ at $m$ then  
\begin{equation}\label{eq:defn-Laplacian}
\triangle_w U :=\sum_{i=1}^d {\rm Hess} \, U[m](\zeta_i, \zeta_i).  
\end{equation}
Here, $\zeta_i$ is the Wasserstein gradient of the $i$--th moment $m \mapsto \int_{\mathbb R^d} x_i m(dx).$ The stochastic path starting at a given point  $m \in \mathcal P_2(\mathbb R^d)$   is 
$$t \mapsto \mathbb B_t^{m}:=(\id+\sqrt 2 W_t)_\# m,$$ 
where, $(W_t)_{t\geq 0}$ is a standard $d$--dimensional Brownian motion starting at $0$. Setting  
\[
V(t, m)=\mathbb E\Bigl(U\bigl(\mathbb B_t^{m}\bigr) \Bigr)
\]
for $(t, m) \in (0,\infty) \times \calPtwod$ we assert    
\begin{equation}\label{eq:wasserstein-heat}
\partial_t V= \triangle_w V, \qquad V(0,\cdot)=U.
\end{equation}
The above facts, certainly contribute to the argument that $\triangle_w$ can be viewed as a partial Laplacian operator on the Wasserstein space and $t \mapsto \mathbb B_t^{m}:=(\id+\sqrt 2 W_t)_\# m$ can be compared to a Brownian motion starting at $m.$ We also study  perturbations $(\sigma^\epsilon, \triangle_{w, \epsilon})$ of $(\mathbb B, \triangle_w)$ for $\epsilon\geq 0$ by setting     
\begin{equation}\label{eq:general-heat1}
\sigma^{\epsilon, \beta}_t[m]:= (\id + \sqrt{2 \beta} W_t)_\# (G_t^\epsilon \ast m).
\end{equation}
Here, $G^\epsilon_t$ is the heat kernel  for the heat equation given by  
\begin{equation}\label{eq:aug.19.2017.ter}
G^\epsilon_t (z)={1\over \sqrt{4\pi \epsilon t}^d} \exp\Bigl(-{|z|^2 \over 4\epsilon t} \Bigr).
\end{equation}
The path $\sigma_t:=\mathbb B_t^{m,  \beta}$ solves the stochastic differential equation (cf. e.g.  \cite{CarmornaDe} \cite{CardaliaguetDLL})
\begin{equation}\label{eq:aug.19.2017.2bis}
d\sigma= {\rm div} \Bigl( \beta  \nabla \sigma  dt-\sqrt{2  \beta}\sigma dW \Bigr) \;\; \text{on} \quad (0,T) \times \mathcal P_2(\mathbb R^d), \qquad \sigma_0=m.
\end{equation}
Similarly, $\sigma^\epsilon_t:=\sigma^{\epsilon, \beta}_t[m]$ solves  the system of stochastic differential equations (cf. e.g. \cite{CarmornaDe} \cite{CardaliaguetDLL}) 
\begin{equation}\label{eq:general-heat3} 
d\sigma^\epsilon= {\rm div} \Bigl( (\epsilon+\beta)  \nabla \sigma^\epsilon  dt-\sqrt{2  \beta}\sigma^\epsilon dW \Bigr) \;\; \text{on} \quad (0,T) \times \mathcal P_2(\mathbb R^d), \quad \sigma_0^\epsilon=m.
\end{equation}
By definition equations \eqref{eq:aug.19.2017.2bis} or \eqref{eq:general-heat3} is satisfied if the linear identity \eqref{eq:aug.19.2017.2def} holds. Then, in Theorem \ref{th:measurable2}, we obtain a nonlinear version of \eqref{eq:aug.19.2017.2def},  It\^o's formula on  the Wasserstein space. 

Note  \eqref{eq:aug.19.2017.2bis} and \eqref{eq:general-heat3}  are the simplest versions of systems of equations recently considered by  \cite{CardaliaguetDLL} and \cite{CarmornaDe}. The purpose of these authors was the study of Mean Field Games involving a Hamiltonian, for large populations and therefore, their goals are quite different from ours.  Notwithstanding the fact that \eqref{eq:aug.19.2017.2bis} or \eqref{eq:general-heat3} are not the subject of our study, we  formulated them as  continuity  equations just  to display the velocities driving the Brownian motions.

Arguably,  $\mathbb B_t^{m}:=\mathbb B_t^{m,1}$   can be called a  Brownian  starting at $m$ and $\triangle_w$ can be called a partial Wasserstein Laplacian operator.  When restricted to finitely many symmetric products of $\mathbb R^d$, the partial Wasserstein Laplacian operator coincides with classical finite dimensional operators. For instance,  suppose $U$ is differentiable in a neighborhood of $\delta_a$ for some $a \in \mathbb R^d$ and $U$ is twice differentiable at $\delta_a$ as stated in Theorem \ref{th:hessian1}. Then $x\mapsto u(x):=U(\delta_x)$ is twice differentiable at $a$  (cf. Remark \ref{re:Laplacian2}) and 
\[
\triangle_w U[\delta_{a}]=\triangle u(a).
\]
However, if we define $v$ on $(\mathbb R^d)^k$ by 
\[
v(x) \equiv v(x_1, \cdots, x_k):=U\biggl({1\over k} \sum_{j=1}^k \delta_{x_j}\biggr)
\]
then unless $k=1$, 
\[
\triangle_w U\biggl[{1\over k} \sum_{j=1}^k \delta_{x_j} \biggr]=\sum_{j=1}^k \triangle_{x_j} v(x)+ \sum_{j \not= l}{\rm div}_{x_j}\bigl( \nabla_{x_l} v(x)\bigr)\not=  \sum_{j=1}^k \triangle_{x_j} v(x). 
\] 

Our study of the partial Laplacian operator will be mainly restricted to the set of $k$--polynomials and their graded sums. These are functions  on $\mathcal P_2(\mathbb R^d)$  of the form 
\[m \mapsto F_{\Phi}[m]:={1\over k}\int_{(\mathbb R^d)^k} \Phi(x) m(dx_1) \cdots m(dx_k), \]  
$\Phi \in C((\mathbb R^d)^k)$ being a symmetric function that grows at most quadratically at infinity. The set of such $\Phi$'s is denoted as $\mathcal{ {\rm Sym[\mathbb R^k]}}$ and the set of $F_{\Phi}$'s is denoted as $\cSk.$ Although at a first glance the sets of $k$--polynomials may appear to be too small, by the Stone Weierstrass Theorem, the subalgebra they generate is a dense subset of $C(\mathcal K)$ for the uniform convergence (cf. Remark \ref{re:stoneW}). Here, $\mathcal K$ is any locally compact subset of $\mathcal P_2(\mathbb R^d).$ As a consequence 
\begin{equation}\label{eq:oct08.2017.1} 
\bigoplus_{k} \Bigl[ \cSk \cap C_c\bigl((\mathbb R^{d})^k\bigr)\Bigr],
\end{equation}
the $\mathbb N$--graded sums of the set of $k$--polynomials, generates a subalgebra  which is a dense subset of $C(\mathcal K).$

The nonnegative real numbers  are contained in the spectrum of $-\triangle_w$.   For any $\beta\geq 0$, it is shown that the intersection of the kernel of $\triangle_w +\beta {\rm Id}$ with $\cSkcomplex$, is represented by a general conical surface, contained in the symmetric $k$--product of $\mathbb R^d.$ The surface in question is  the quotient space 
\begin{equation}\label{eq:sep03.2017.1} 
\Bigl\{(\xi_1, \cdots, \xi_k) \in (\mathbb R^d)^k\: \: 4\pi^2 \Big|\sum_{j=1}^k \xi_j \Big|^2= \beta\ \Bigr\}/P_k,
\end{equation} 
where $P_k$ is the set of permutations of $k$ letters. Note when $k>1$ and $\beta=0$, the surface degenerates into a linear space, and so, it has infinitely many elements, which means that the kernel of $\triangle_w$ has infinitely elements. More serious is the fact that the surface is unbounded,  which precludes the  Wasserstein Laplacian operator to have a smoothing property, unless restricted to an appropriate set of functions. Solving the simplest case of Poisson equation on the Wasserstein space amount to, given $$a \in \cSkcomplex \cap L^2((\mathbb R^d)^k) \cap  L^1((\mathbb R^d)^k),$$ solving  
\[
-4\pi^2 \big|\sum_{j=1}^k \xi_j \big|^2 b(\xi_1, \cdots, \xi_k)= a(\xi_1, \cdots, \xi_k). 
\]
This, obviously is not an elliptic equation as when $k>1$ and $\beta=0$,  the surface in \eqref{eq:sep03.2017.1} does not reduce to the null vector. 

In this manuscript, we also address the following natural and useful question:  suppose we know that a function $F:\mathcal P_2(\mathbb R^d) \rightarrow \Reals$ is of the form 
$F=F_\Phi$ for a symmetric function $\Phi: (\mathbb R^d)^k \rightarrow \Reals.$ Can we reconstruct $\Phi$?  We can convince ourselves that the problem reduces to expressing $\Phi$ as a sum, up to a multiplication constant, of the so--called  $k$--th defects  of $F$. When we do not require any differentiability property of $F,$ we reconstruct $\Phi$ by providing a polarization isomorphism based on the inclusion--exclusion principle, without any differentiation operations. Our arguments was inspired by works on vector spaces, which can be traced back to \cite{Greenberg} in a particular case, followed by generalization in  \cite{Ward}.

The Wasserstein Laplacian is the sum of two operators, one being nonpositive with a trivial kernel when restricted to $\mathcal{S}ym[k](\Reals)$. The latter alluded operator, which has a smoothing effect, associates to any smooth function $U: \mathcal P_2(\mathbb R^d) \rightarrow \Reals$, the function 
\[ m \rightarrow O[m]:= \int_{\mathbb R^d}{\rm div}_x \bigl(\nabla_w U[m](x)\bigr) m(dx).\] 
Here, $\nabla_w$ denotes the Wasserstein gradient operator. Given $\epsilon>0$ and a  twice continuously differentiable function $U$,  
\begin{equation}\label{eq:wasserstein-heat-late}
V(t, m):=\mathbb E\Bigl(U\bigl(\sigma^\epsilon_t[m]\bigr) \Bigr), \qquad \text{on} \quad (0,\infty) \times \mathcal P_2(\mathbb R^d).
\end{equation}
solves the initial value differential equation 
\[
\partial_t V=(\triangle_w+\epsilon O)V, \qquad V(0,\cdot)=U.
\]

For each $s\geq 0$, we define $H^s(\mathcal P_2(\mathbb R^d)),$ a space of functions on the Wasserstein space, in the spirit of the Sobolev functions. It has the virtue  that whenever $U \in H^s(\mathcal P_2(\mathbb R^d)),$ then $V(t, \cdot)$ given in \eqref{eq:wasserstein-heat-late} not only solves the previous differential equation but 
 $$V(t, \cdot) \in H^l\bigl(\mathcal P_2(\mathbb R^d)\bigr) \quad \forall\, l \; \geq 0.$$ 
 We identify a set in ${\mathcal H}^s(\mathcal P_2(\mathbb R^d)) \subset H^s(\mathcal P_2(\mathbb R^d))$ such that choosing $U_0$ in that set, even if $U_0$ is not three times differentiable, when $s, \epsilon>0$, then $\triangle_{w,\epsilon} V(t, \cdot)$ becomes twice differentiable (cf. Remark \ref{re:superposition}). This is an improved  smoothing effect in the $m$ variable. 
 
 The functions  $F$ in ${\mathcal H}^s(\mathcal P_2(\mathbb R^d))$  are pointwise sums of the infinite series $\sum_{k=1}^\infty 1/k! F_{\Phi_k}$, where $\Phi_k \in L^2((\mathbb R^d)^k)$ is a symmetric function which grows at most super linearly at $\infty$  and is such that its inverse Fourier transform $a_k$, satisfies \eqref{eq:aug23.2017.1}. Under  more stringent assumptions on $\Phi_k,$ so that $F \in {\mathcal H}^s(\mathcal P_2(\mathbb R^d))$, we prove that $F_{\Phi_k}$  is uniquely determined by a specific projection operator $\pi_k$ (cf. Remark \ref{re:unique-proj}) defined on a subset of the graded sum in \eqref{eq:oct08.2017.1} .

A bilinear form 
$$
\langle \cdot \;; \; \cdot \rangle_{H^0}: \mathcal H^0\bigl(\mathcal P_2(\mathbb R^d)\bigr) \times \mathcal H^0\bigl(\mathcal P_2(\mathbb R^d)\bigr) \rightarrow \Reals
$$  
which involving functions, their gradients and Laplacians, is provided in Proposition \ref{pr:integration-by-parts}.  Under appropriate conditions on $F, G \in \mathcal H^0(\mathcal P_2(\mathbb R^d))$, we assert  
\[
-\langle \triangle_w F; G \rangle_{H^0}= \left \langle \int_{\mathbb{R}^d} \nabla_w F[m](x)m(dx) ; \int_{\mathbb{R}^d} \nabla_w G[m](x)m(dx)  \right \rangle_{H^0} .
\]
In some cases, this turns into a integration by parts formula involving signed Radon measures $\mathbb P^{k ,R}$ on $\mathcal P_2(\mathbb R^d)  \times \mathcal P_2(\mathbb R^d) .$ When $F=\Phi$, $G=F_\Psi$ where $\Phi, \Psi$ are $k$--symmetric functions of class $C^3$ and supported by the ball of radius $R$, Theorem \ref{th:measureR}  shows the above  to be equivalent to 
\[
-\int_{\mathcal M^2} \triangle_w F_\Phi[m_1] F_\Psi[m_2]  d\mathbb P^{k ,R} =\int_{\mathcal M^2} D_2(\nabla_w F_\Phi, \nabla_w G_\Psi)  d\mathbb P^{k ,R}. 
\]
Here $D_2$ is the bilinear function 
\[
D_2(\nabla_w F, \nabla_w G)(m_1, m_2):=  \int_{\mathbb R^{2d}}  \langle \nabla_w F[m_1](q_1)  ;  \nabla_w G[m_2](q_2) \rangle m_1(dq_1) m_2(dq_2). 
\]

In the recent years, there have been many attempts to construct a ``full'' Laplacian on the Wasserstein space.  In \cite{VonRenesseS}, von Renesse and Sturm studied a canonical diffusion process on the Wasserstein space, when the underlying space is the one--dimensional torus. Then in \cite{Sturm}, Sturm constructed  entropic measures on Wasserstein spaces $\mathcal P(M)$, where the underlying set $M$ is a compact manifold of finite dimension. Unlike the case when $M$ is a one--dimensional set, the closability of the Dirichlet form associated to the entropic measures, remains to--date, an outstanding open question. We end this introduction by drawing the attention of the reader to a (far from being exhaustive) literature which studies infinite dimensional Laplacian operators  on flat spaces. The first one due to Levy  \cite{Levy}, relies on a concept of the mean of a function on a Hilbert space, to propose a Laplacian operator. No meaningful subset of the domain of definition of this operator was known until a later studied by Dorfman \cite{Dorfman}. This author proves, when the Hessian has the form ${\rm Hess}U(x)=r(x) I + T(x)$, where $r$ is uniformly continuous and $T$ satisfies a so--called $N$--property, then $U$ belongs to the domain of definition of Levy's Laplacian operator. Other definitions of a Laplacian operators  on a Hilbert space appeared in the literature. For instance, \cite{Umemura} consideres a Hilbert space $\mathbb D$ and a nuclear space $\mathbb L$ and defined Laplacians on subsets of $L^2(\mathbb L^*).$ The previously mentioned studies raised many new challenging questions, none of which we attempt to pursue here. We rather  take a different turn and search for infinitesimal generators of stochastic paths on the Wasserstein space, which are partial traces. Our interests include the study of  infinitesimal generators which could have smoothing effects on partial differential equations on the Wasserstein space. A first evidence to this fact steams out of the College of France lectures by P--L. Lions \cite{PLLions2007-2008}  \cite{PLLions2008-2009} and the pioneering work by Cardaliaguet et. al.  \cite{CardaliaguetDLL}, in their study of the so--called {\it master equation} in mean field game systems. These are first or second order Hamilton--Jacobi equation on the Wasserstein space, with a non--local term.  Their equation incorporated terms which turned out to be $\triangle_{w, \epsilon} U$. When $\epsilon>0$, they referred  to games with individual noise and common noise. The presence of $\triangle_{w, \epsilon} U$ was instrumental for the well--posedness of the {\it master equation}. For more discussions on the topic, we refer to \cite{CardaliaguetDLL}  \cite{CarmornaDeI} \cite{CarmornaDeII}.

\section{Notation and Preliminaries} 

\subsection{Notation}
 
In this manuscript, if $(\mathcal S, {\rm dist})$ is a metric space, the domain of $U: \mathcal S \rightarrow \Reals \cup\{\pm \infty\}$ is the set ${\rm dom}(U)$ of $m \in \mathcal S$ such that $U[m] \in \mathbb R.$

A function $\rho: [0, +\infty ) \to [0,+\infty )$ is called a modulus if $\rho $ continuous, nondecreasing, sub--additive, and $\rho (0) =0$. It is a modulus of continuity for $U: \mathcal S \rightarrow \Reals$ if $|U(s_2)-U(s_1)| \leq \rho({\rm dist}(s_1, s_2))$ for any $s_1, s_2 \in \mathcal S.$

We denote $\mathbb R^d$ as $\mathbb M$ because it is more convenient to write expressions such as $\Mb^2$ than $(\mathbb R^d)^2$. This notation is also meant to emphasize the fact that most of our results proven in this manuscript are valid on spaces more general than $\mathbb R^d$.  Throughout this manuscript, $\calPtwod$ denotes the set Borel probability measures on $\Mb,$  of finite second moments. This is a length space when endowed with $W_2$, the Wasserstein distance.

Given $m, \nu \in \calPtwod$ we denote as $\Gamma(m, \nu)$ the set of Borel measures $\gamma$ on $\Mb^2,$ which have $m$ as their first marginal and $\nu$ as their second marginal.  We denote as $\Gamma_0(m, \nu)$, the set of $\gamma \in \Gamma(m, \nu)$ such that 
\[ 
W_2^2(m, \nu)= \int_{\Mb^2}|x-y|^2 \gamma(dx, dy). 
\]

We denote the first  (resp. second)  projection of $\Mb^2$ onto $\Mb$ as $\pi^1$ (resp. $\pi^2$) 
\[
\pi^1(x, y)=x, \quad (\text{resp.} \;\; \pi^2(x,y)=y).
\]

Let $L^2(m)$ denote the set of Borel maps $\zeta: \Mb \rightarrow \Mb$ such that $\|\zeta \|_{m}^2:=\int_{\Mb} |\zeta(x)|^2 m(dx)<\infty.$ This is a Hilbert space with the inner product $\langle \cdot ; \cdot\rangle_m$ such that   
\[
\langle\zeta_1 ; \zeta_2\rangle_m= \int_\Mb \langle \zeta_1(x) , \zeta_2(x) \rangle m(dx). 
\]
Let $T_m \calPtwod$ denote the closure of $\nabla C_c^\infty(\Mb)$ in $L^2(m)$, and let us denote the orthogonal projection of $L^2(m)$ onto $T_m \calPtwod$ as $\pi_m$. The union of all the sets $\{m\} \times L^2(m)$ is denoted as $\mathcal T \calPtwod$ and by an abuse of language, is referred to as the tangent bundle of $\calPtwod.$ 

Let $P_k$ denote the set of permutations of $k$ letters. If $a \in \mathbb C$ we denote its complex conjugate as $a^*.$

Let $\Sk$ be the set of $\Phi \in C(\Mb^k)$ such that there exists $C>0$ such that for any $x=(x_1, \cdots, x_k) \in \Mb^k,$ 
\[ 
|\Phi(x)| \leq (1+|x|^2) \quad   \text{and} \quad \Phi(x)= \Phi(x^\sigma) \quad \text{for any} \quad \sigma \in P_k.
\]   
In other words, $\Phi$ is well--defined on the $k$--symmetric product of $\Mb.$ In this case, we call  $\Phi$  is symmetric.  Let $\Skt$ be the set of $\Phi \in \Sk \cap C^2(\Mb^k)$ such that $\nabla^2 \Phi$ is bounded and uniformly continuous on $\Mb^k$ and $\nabla^2 \Phi$ has a modulus of continuity which is a concave function. Similarly, we define $\Skcomplex$ and $\Sktcomplex$ when the functions are taking complex values.  

When $\Phi \in \Sk$, the function $F_\Phi: \calPtwod \rightarrow \Reals$  set to be  
\begin{equation}\label{eq:aug19.2017.4}
F_\Phi[m]:={1 \over k} \int_{\Mb^k} \Phi(x_1, \cdots, x_k)m(dx_1) \cdots m(dx_k).
\end{equation} 
is well--defined. We denote as $\cSk$ the set of $F_\Phi$ such that $\Phi \in \Sk$ and denote as $\cSkt$ the set of $F_\Phi$ such that $\Phi \in \Skt.$

The Fourier transform of $\Phi\in \Sk$ is the function $\widehat \Phi \in \cSk$ defined by 
\[
\widehat \Phi(\xi)= \int_{\Mb^k} \exp \Bigl(-2\pi \sum_{j=1}^k \langle \xi_{j}, x_j\rangle \Bigr) \Phi(x) dx.
\]
We denote the Fourier inverse of $A \in L^2(\Mb^k)$ as $\check{A}.$

\subsection{Preliminaries}  Let $U: \calPtwod \rightarrow [-\infty, \infty]$ denote a  function with values in the extended real line. The recent work \cite{GangboT2017},   shows  two notions of Wasserstein subgradient which appeared in the literature to be equivalent. For that reason, we recall once more these definitions and state in Remark \ref{rem:infdef} that they are equivalent. 
 
\begin{definition}\label{defn:infdef} Let $m\in {\rm dom}(U)$ and let  $\zeta\in L^2(m)$. 
\begin{enumerate}
\item[(i)] We call $\zeta$ a subgradient of $U$ at $m$ and write $\zeta\in \partial_\cdot U[m]$ if  for any $\nu \in \calPtwod$
\begin{equation}\label{eq:infdef1} 
U[\nu]-U[m]\geq  \inf_{\gamma\in\Gamma_0(m, \nu)} \int_{\Mb^2} \zeta(x) \cdot (y-x) \gamma(dx, dx)  +o\big(W_2(m,\nu)\big),
\end{equation}
\item[(ii)] We call $\zeta$ a supergradient of $U$ at $m$ and write $\zeta \in \partial^\cdot U[m]$  if $-\zeta \in \partial_\cdot (-U)[m]$.
\end{enumerate}  
\end{definition}

\begin{remark}\label{rem:infdef} Let $m$ and $\zeta$ be as in Definition \ref{defn:infdef}.
\begin{enumerate}
\item[(i)] It has recently been shown \cite{{GangboT2017}} that $\zeta \in \partial_\cdot U[m]$ if and only if  \eqref{eq:infdef1} holds when we replace the ``inf''   by ``sup''.  
\item[(ii)] It is well--known that in case  both $\partial_\cdot U[m]$ and $\partial^\cdot U[m]$ are not empty then they coincide.
\end{enumerate}
\end{remark}

\begin{definition}\label{defn:infdef2} Let $m\in {\rm dom}(U).$ 
\begin{enumerate} 
\item[(i)] We say that $U$ is differentiable at $m$ if both $\partial_\cdot U [m]$ and $\partial^\cdot U [m]$ are nonempty. In this case, we set $\partial U [m]=\partial^\cdot U[m]$.
\item[(ii)] If $\partial_\cdot U[m]$ is nonempty, then it is a closed convex set in the Hilbert space $L^2(m)$ and so, it has a unique element of minimal norm. As customary done in convex analysis, we denote this element as $\nabla_w U[m]$ and refer to it as the Wasserstein gradient of $U$.
\item[(iii)] Thanks to Remark \ref{rem:infdef} (ii), if $\partial^\cdot U[m]$ is not empty, there is no confusion referring to its unique element of minimal norm as the Wasserstein gradient of $U$ at $m$.  
\item[(iv)] If $\nabla_w U[m]$ exists, $dU[m]:L^2(m) \rightarrow \Reals$ denotes the linear form  $\zeta \mapsto  \zeta \cdot (U[m]) :=  \langle \nabla_w U[m];\zeta \rangle_m.$ 
\end{enumerate}  
\end{definition}

\begin{remark}\label{rem:grad-linear} Note  that if $\phi \in C_c^\infty(\Mb)$, then according to Definition \ref{defn:infdef2}, the Wasserstein gradient of $F_\phi$ is $\nabla \phi.$ 
\end{remark}

\begin{remark}\label{rem:modulus} Assume $\rho:[0,\infty) \rightarrow [0,\infty)$ is a concave modulus. 
\begin{enumerate} 
\item[(i)] Then $\rho(t)/t$ is monotone nonincreasing and so, for any $t\geq 0$ and $\epsilon>0$, we have $\rho(t) \leq \rho(\epsilon) + t/\epsilon \rho(\epsilon).$ 
\item[(ii)] If $m, \nu \in \calPtwod$ and $\gamma \in \Gamma_0(m, \nu)$ then 
\[
 \int_{\Mb^2} |x-y| \rho\bigl(|x-y| \bigr)\gamma(dx, dy) \leq \rho\bigl(W^{1\over 2}_2(m, \nu)\bigr) W_2(m, \nu) \Bigl( W^{1\over 2}_2(m, \nu) +1\Bigr)
\]
\end{enumerate}  
\end{remark}
\proof{} (i) Mollifying $\rho$ if necessary, it is not a loss of generality to assume that $\rho$ is of class $C^1.$ We have $t^2 (\rho(t)/t)'= \rho'(t)t-\rho(t).$ But $R(t):=-\rho(t)$ is convex and so, for $t>0$ we have $R(0)- R(t)\geq R'(t)(0-t)$. This is equivalent to $t^2 (\rho(t)/t)'\leq0$ which proves the first part of the remark. If $s \in [0,\epsilon]$ and $t\in [\epsilon, \infty)$ then 
\[
\rho(s)\leq \rho(\epsilon) \leq \rho(\epsilon) + {s \over \epsilon} \rho(\epsilon) \quad \text{and so,} \quad {\rho(t) \over t} \leq {\rho(\epsilon) \over \epsilon} \leq  { \rho(\epsilon)\over t}+ {\rho(\epsilon) \over \epsilon}.
\] 
This proves (i). 

(ii) Let $\gamma \in \Gamma_0(m, \nu)$. By (i) for any $\epsilon>0$, 
\[
\int_{\Mb^2} |x-y| \rho\bigl(|x-y| \bigr)\gamma(dx, dy) \leq \rho(\epsilon) \Bigl( {W_2^2(m, \nu)\over \epsilon} + \int_{\Mb^2} |x-y|   \gamma(dx, dy)\Bigr).
\]
We apply Cauchy--Schwarz inequality to obtain 
\[
\int_{\Mb^2} |x-y| \rho\bigl(|x-y| \bigr)\gamma(dx, dy) \leq \rho(\epsilon) W_2(m, \nu) \Bigl( {W_2(m, \nu) \over \epsilon}+1\Bigr).
\]
We conclude the proof by setting $\epsilon:= W^{1/2}_2(m, \nu)$. \endproof

Since $\Mb$ is not a compact set, the   space $\calPtwod$ is not a locally compact space (cf. e.g. \cite{AGS}). Suppose $\phi: \Mb \rightarrow [0,\infty]$ is a lower semicontinuous monotone nondecreasing function 
$$\phi: \Mb \rightarrow [0,\infty], \quad \lim_{|x| \rightarrow \infty} {\phi(x) \over |x|^2}=\infty.$$ 
Consider the  locally compact set 
\begin{equation}\label{eq:oct01.2017.4}
\mathcal P_\phi(\Mb):= \biggl\{m \in \calPtwod \; \Big | \; \int_{\Mb} \phi(x) m(dx)<\infty  \biggr\}.
\end{equation} 

\begin{remark}\label{re:stoneW}  The subalgebra generated by  $\bigl\{F_\Phi\; | \; \Phi \in C_c(\Mb)\bigr\},$
separates points in $\mathcal P_{\phi}(\Mb)$ and vanishes nowhere. Hence, by the Stone Weierstrass Theorem, it is a dense subset of $C\bigl( \mathcal P_{\phi}(\Mb)\bigr)$ for the uniform convergence. 
\end{remark}

\section{Traces of second order derivative: the partial Laplacian operators} 
\noindent Let $U: \calPtwod \rightarrow [-\infty, \infty]$ denote a  function with values in the extended real line.

\subsection{The Wasserstein Laplacian operator}  Consistent with Levi--Civita connection in \cite{Lott2008}, we have the following definition.

\begin{definition}\label{defn:infdef2hessian} Suppose $U$ is differentiable in a neighborhood of $m\in {\rm dom}(U)$ and for any $\zeta \in C_c^\infty(\Mb, \Mb)$,  $\nu \mapsto  \zeta \cdot (U[\nu])$ is differentiable at $m.$ 
\begin{enumerate} 
\item[(i)] We define ${\rm \bar Hess} \,U [m]: \nabla C_c^\infty(\Mb) \times \nabla C_c^\infty(\Mb) \rightarrow \Reals$ if the following exists:  
\[
{\rm \bar Hess} \,U [m](\zeta_1, \zeta_2)=\zeta_1 \cdot \Bigl( \zeta_2\cdot \bigl(U[m]\bigr)\Bigr)- \Bigl(\bar \nabla_{\zeta_1} \zeta_2\Bigr)\cdot (U[m])
\]
for any $\zeta_1, \zeta_2 \in \nabla C_c^\infty(\Mb).$ Here, $\bar \nabla_{\zeta_1} \zeta_2= \nabla \zeta_2 \zeta_1.$
\item[(ii)] If  there is a constant $C$ such that  $|{\rm Hess} \,U [m](\zeta_1, \zeta_2)| \leq C \, \|\zeta_1\|_m \, \|\zeta_2\|_m$ for all $\zeta_1, \zeta_2 \in \nabla C_c^\infty(\Mb)$ then  ${\rm \bar Hess}[m]$ has a unique extension onto $T_m \calPtwod \times T_m \calPtwod$ which we denote as ${\rm Hess} \,U [m].$ In that case, we say that $U$ has a Hessian at $m.$ 
\end{enumerate}  
\end{definition}

\begin{theorem}\label{th:hessian1} Suppose $U$ is differentiable in a neighborhood of $m\in {\rm dom}(U)$, $x \mapsto \nabla_w U[\nu](x)$ is continuous for $\nu$ in a neighborhood of $m$  and there exists a constant $C_m$ such that $|\nabla_w U[\nu](x)| \leq C_m(1+|x|)$ for any $x\in \Mb$ and any $\nu$ in the neighborhood of $m$. Suppose $\rho, \epsilon_2:[0,\infty) \rightarrow \Reals$ are nonnegative function (depending on $m$) such that $\lim_{t \rightarrow 0^+} \rho(t)=\lim_{t \rightarrow 0^+} \epsilon_2(t) =0$ and $\rho$ is  a concave modulus. Suppose there are Borel bounded matrix valued  functions $\tilde A[m] : \Mb \rightarrow \mathbb R^{d \times d}$ and  $A_{m m}: \Mb^2 \rightarrow \mathbb R^{d \times d}$ such that for any $\nu \in \calPtwod$ we have 
\[
\sup_{\gamma \in \Gamma_0(m, \nu)}\Bigl|
\nabla_w U[\nu](y) -\nabla_w U[m](x)- P_\gamma[m](x, y)\Bigr| \leq |x-y|\rho\bigl(|x-y|\bigr)+ W_2(m, \nu)\epsilon_2\bigl(W_2(m, \nu)\bigr) \,,
\]
where, for $\gamma \in \mathcal P(\Mb^2)$ and $x, y \in \Mb,$ we have set  
\begin{equation}\label{eq:aug19.2017.1}
P_\gamma[m](x, y):= \tilde A[m](x) (y-x)+ \int_{\Mb^2} A_{mm} (x,a) (b-a) \gamma(da, db)\,.
\end{equation}
Then, $U$ has a Hessian at $m$  and 
\[
{\rm Hess} \,U[m](\zeta_1, \zeta_2)=\int_{\Mb} \langle \tilde A[m](x) \zeta_1(x) , \zeta_2(x)\rangle m(dx) +  \int_{\Mb^2} \langle A_{mm}(x,a) \zeta_1(a) , \zeta_2(x)\rangle m(dx) m(da)
\]
for $\zeta_1, \zeta_2 \in T_m \calPtwod.$ 
\end{theorem}
\proof{} Fix $\zeta_1, \zeta_2 \in C_c^\infty(\Mb, \Mb).$ We are to show that the map $\nu \mapsto \Lambda(\nu):= dU[\nu](\zeta_2)$ is differentiable at $m$ and then show that $\zeta_1 \cdot \bigl(dU[\nu](\zeta_2)\bigr)$ exists.

Let $\nu \in \calPtwod,$ $\gamma \in \Gamma_0(m, \nu)$ and set $A_\nu:=\nabla_w U[\nu].$  Since $\zeta_2$ is of compact support, its first and second derivatives are bounded and so, we may choose bounded vector fields $v, w \in C(\Mb^2, \Mb)$ such that for any $x, y \in \Mb$ we have 
\begin{equation}\label{eq:aug18.2017.1}
\zeta_2(y)= \zeta_2(x)+\nabla \zeta_2(x)(y-x)+|y-x|^2 v(x,y), \quad \zeta_2(y)- \zeta_2(x)= |y-x| w(x,y).
\end{equation}
By assumption, for each $\gamma \in \mathcal P(\Mb^2)$, there exists  $l(\gamma, \nu,y) \in \Mb$ such that $|l(\gamma, \nu,y)|\leq 1$ and 
\begin{equation}\label{eq:aug18.2017.1b}
A_\nu(y)-A_m(x)- P_\gamma[m](x, y)=l(\gamma, \nu,y ) \bigl(|x-y|\rho\bigl(|x-y|\bigr)+ W_2(m, \nu)\epsilon_2\bigl(W_2(m, \nu)\bigr) \bigr)
\end{equation}
We have then by definition of the map $\Lambda$ that
\begin{align}
\Lambda(\nu)-\Lambda(m)
&= \int_{\Mb^2} \bigl(\langle A_\nu(y), \zeta_2(y) \rangle-\langle A_m(x), \zeta_2(x) \rangle \bigr)\gamma(dx, dy) \nonumber\\
&=  \int_{\Mb^2} \bigl(\langle A_\nu(y)-A_m(x) , \zeta_2(y) \rangle +   \langle A_m(x), \zeta_2(y)- \zeta_2(x) \rangle \bigr)\gamma(dx, dy). \nonumber
\end{align}
This, combined with \eqref{eq:aug18.2017.1} yields 
\begin{align}
\Lambda(\nu)-\Lambda(m) 
&=  \int_{\Mb^2} \langle A_\nu(y)-A_m(x) , \zeta_2(x) + |y-x| w(x,y) \rangle  \gamma(dx, dy) \nonumber\\
&+  \int_{\Mb^2}  \langle A_m(x), \nabla \zeta_2(x)(y-x)+|y-x|^2 v(x,y)  \rangle \bigr)\gamma(dx, dy)\nonumber\\
&=I +II. \label{eq:aug18.2017.2}
\end{align}
Note 
\begin{equation}\label{eq:aug18.2017.3}
\Bigl|II-  \int_{\Mb^2}  \langle \nabla \zeta_2^T (x)A_m(x), (y-x)  \rangle \gamma(dx, dy) \Bigr| = \Bigl| \int_{\Mb^2} |y-x|^2 v(x,y) \gamma(dx, dy) \Bigr| \leq W_2^2(m, \nu).
\end{equation}
By \eqref{eq:aug18.2017.1b} 
\begin{align}
I &=  \int_{\Mb^2} \Bigl\langle \tilde A[m](x) (y-x)+ \int_{\Mb^2} A_{mm} (x,a) (b-a) \gamma(da, db) , \zeta_2(x) \Bigr\rangle  \gamma(dx, dy) \nonumber\\ 
&+ \int_{\Mb^2} \Bigl\langle \tilde A[m](x) (y-x) +\int_{\Mb^2} A_{mm} (x,a) (b-a) \gamma(da, db), |y-x| w(x,y) \Bigr \rangle  \gamma(dx, dy) \nonumber\\
&+  \int_{\Mb^2} \Bigl \langle l(\gamma, \nu,y) \bigl(|x-y|\rho(|x-y|) +W_2(m, \nu)\epsilon_2 ( W_2(m, \nu) \bigr) , \zeta_2(y) \Bigr\rangle  \gamma(dx, dy) 
\nonumber\\ 
& =III +IV+V.       \label{eq:aug18.2017.4} 
\end{align}
We have by Cauchy--Schwarz inequality that 
\[
|V| \leq \|\zeta_2\|_m W_2(m, \nu) \epsilon_2(W_2(m, \nu))+ 
\|\zeta_2\|_{L^\infty} \int_{\Mb^2}  \bigl(|x-y|\rho(|x-y|)  \gamma(dx, dy).
\]
We apply Remark \ref{rem:modulus} to infer  
\begin{equation} \label{eq:aug18.2017.5}
|V| \leq W_2(m, \nu) \Bigl( \|\zeta_2\|_m  \epsilon_2(W_2(m, \nu))+ 
\|\zeta_2\|_{L^\infty}\rho\bigl(W^{1\over 2}_2(m, \nu)\bigr)  \bigl( W^{1\over 2}_2(m, \nu) +1\bigr) \Bigr)
\end{equation}
Checking also that   
\[
\|\tilde A[m](x) (y-x) +\int_{\Mb^2} A_{mm} (x,a) (b-a) \gamma(da, db) \|_\gamma \leq  \bigl(\|\tilde A[m]\|_{L^\infty(m)}+ \|A_{mm}\|_{L^2(m \otimes m)} \bigr)  W_2(m, \nu) \,,
\]
we conclude 
\begin{equation}\label{eq:aug18.2017.6}
|IV| \leq  \bigl(\|\tilde A[m]\|_{L^\infty(m)}+ \|A_{mm}\|_{L^2(m \otimes m)} \bigr)   W_2^2(m, \nu).  
\end{equation}
We combine (\ref{eq:aug18.2017.2}-\ref{eq:aug18.2017.6}) and make the substitution $(a,b) \leftrightarrow (x, y)$ to obtain
\begin{align}
\Lambda(\nu)-\Lambda(m)&=   \int_{\Mb^2}  \langle\tilde A[m](x)+ \nabla \zeta_2^T (x)A_m(x), (y-x)  \rangle \gamma(dx, dy) \nonumber\\ 
&+\int_{\Mb^2} \Bigl\langle  \int_{\Mb^2}  (y-x) , A^T_{mm} (a,x)\zeta_2(a) \Bigr\rangle  \gamma(da, db)  \gamma(dx, dy) \nonumber\\
&+ o \bigl(  W_2(m, \nu) \bigr) \nonumber 
\end{align} 
Thus, 
\[
\nabla_w \Lambda[m](x)= \tilde A[m](x)+ \nabla \zeta_2^T (x)A_m(x)+  \int_{\Mb} A^T_{mm} (a,x)\zeta_2(a)m(da) 
\]
Consequently, 
\begin{align}
\zeta_1\cdot (\Lambda[m] ) &= \int_\Mb \langle \tilde A[m](x) , \zeta_1(x) \rangle m(dx) 
+\int_{\Mb^2} \langle A_{mm} (a,x) \zeta_1(x), \zeta_2(a) \rangle m(da) m(dx) \nonumber\\
&+\int_\Mb \langle A_m(x) , \nabla \zeta_2 (x) \zeta_1(x) \rangle m(dx). \label{eq:aug18.2017.7}
\end{align} 
Making the substitution $a \leftrightarrow x$ \eqref{eq:aug18.2017.7} and using definition \ref{defn:infdef2hessian}, reads off 
\begin{equation}\label{eq:aug18.2017.8}
{\rm \bar Hess} \,U [m](\zeta_1, \zeta_2)=\int_{\Mb} \langle \tilde A[m](x) \zeta_1(x) , \zeta_2(x)\rangle m(dx) +  \int_{\Mb^2} \langle A_{mm}(x,a) \zeta_1(a) , \zeta_2(x)\rangle m(dx) m(da).
\end{equation}
Obviously, there is a constant $C$ such that  $|{\rm \bar Hess} \,U [m](\zeta_1, \zeta_2)| \leq C \, \|\zeta_1\|_m \, \|\zeta_2\|_m$. Thus ${\rm Hess} \, U [m]$ is well--defined and \eqref{eq:aug18.2017.8} remains valid for $\zeta_1, \zeta_2 \in T_m \calPtwod.$ \endproof

\begin{definition}\label{defn:Laplacian} Under the assumptions of Theorem \ref{th:hessian1}, we say that $U$ is twice differentiable at $m\in {\rm dom}(U).$ If $A_{m m}$ satisfies \eqref{eq:aug19.2017.1}, so does $\pi_m \bigl(A_{m m}\bigr)$ which is the matrix whose rows are the orthogonal projections onto $T_m \calPtwod$ of the rows of $A_{m m}$. Note that although  $A_{m m}$ may not be unique, $\pi_m \bigl(A_{m m}\bigr)$ is uniquely determined. In the sequel, we tacitly assume that $A_{m m} \equiv \pi_m \bigl(A_{m m}\bigr).$    
\begin{enumerate}  
\item[(i)] We define the Wasserstein Laplacian operator $\triangle_w$ such that $\triangle_w U:\calPtwod \rightarrow \Reals$, is given by 
\[
(\triangle_w U)[m] =\sum_{i=1}^d {\rm Hess} \,U[m](\zeta_i, \zeta_i), 
\]
where $\zeta_i$ is the Wasserstein gradient of the $i$--th moment $m \mapsto \int_{\mathbb R^d} x_i m(dx).$
\item[(ii)] Given $\epsilon\geq 0$, we call $\triangle_{w, \epsilon}$ such that $\triangle_{w, \epsilon} U:\calPtwod \rightarrow \Reals$, the $\epsilon$--perturbation of $\triangle_w$ such that  
\[
(\triangle_{w, \epsilon} U)[m] =\sum_{i=1}^d \Bigl( {\rm Hess} \,U[m](\zeta_i, \zeta_i) +\epsilon \int_{\Mb} \langle \tilde A[m](x) \zeta_i(x) , \zeta_i(x)\rangle m(dx)\Bigr)
\]
\item[(iii)] We define the second order Wasserstein gradient of $U$ at $m$ to be $\pi_m \bigl(A_{m m}\bigr)$ and we denote it as $\nabla_w^2 U[m].$ 
\item[(iv)] Suppose further that $U$ is twice differentiable in a neighborhood of $m$. If $(x, \nu) \mapsto \tilde A[\nu](x)$ and $(x, y,\nu) \mapsto \pi_\nu \bigl(A_{\nu \nu}\bigr)(x,y)$ are continuous, we say that $U$ is twice continuously differentiable on that neighborhood. 
\end{enumerate}  
\end{definition}

\begin{proposition}\label{pro:hessian1} Let $U$ be as in Theorem \ref{th:hessian1}, which in particular means that we have fixed $m \in \calPtwod$ such that $U$ is differentiable in a neighborhood of $m\in {\rm dom}(U)$ and $U$ is twice differentiable at $m$. Let $T>0$ and suppose $\sigma \in AC_2(0,T,\calPtwod)$ is a path which has a velocity of minimal norm ${\mathbf v} \in C^1((0,T) \times \Mb)$ which is bounded and has bounded first order time and space derivatives.  If $s \in (0,T)$ and $m=\sigma_s$ then 
\[
{d^2 \over dt^2} U(\sigma_t))\Big|_{t=s}= {\rm Hess}\, U[\sigma_s]({\mathbf v}_s, {\mathbf v}_s)+ \langle \partial_t {\mathbf v}_s + \nabla {\mathbf v}_s {\mathbf v}_s; \nabla_w U[\sigma_s]\rangle_{\sigma_s}
\]
\end{proposition}
\proof{} We skip the proof of this proposition since it is similar to that of Theorem \ref{th:hessian1}. The only new ingredient to use here is the following additional remark:  if $\gamma_h \in \Gamma_0(\sigma_t, \sigma_{t+h})$ and $\pi^1, \pi^2: \Mb^2 \rightarrow \Mb$ denote the standard projections then 
\[
\lim_{h\rightarrow 0} \Bigl( \pi^1, {\pi^2-\pi^1 \over h}\Bigr)_\# \gamma_h=(\id , {\mathbf v}_t)_\#  \sigma_t  \qquad \text{in} \quad \mathcal P_2(\Mb^2).
\] \endproof

\begin{remark}\label{re:Laplacian2} Suppose the assumptions in Theorem \ref{th:hessian1} holds, i.e. $U$ is twice differentiable at $m\in {\rm dom}(U).$ Then  
\begin{enumerate} 
\item[(i)] $ A_m:=\nabla_w U[m]$ is differentiable on $\Mb$ and its gradient (w.r.t. the $x$ variable) is $\tilde A[m],$ whose rows belong to $T_m \calPtwod.$   
\item[(ii)] We have 
 \[
 \triangle_w U[m]= \int_{\Mb} {\rm div}_x \bigl(\nabla_w U[m](x)\bigr) m(dx) +  \int_{\Mb^2}  {\rm Tr} \Bigl( \nabla_w^2 U [m](x,a)\Bigr) m(dx) m(da)
 \]
and
 \[
(\triangle_{w, \epsilon} U)[m] = (1+\epsilon) \int_{\Mb} {\rm div}_x \bigl( \nabla_w U[m] (x)\bigr) m(dx) +  \int_{\Mb^2}  {\rm Tr} \Bigl( \nabla_w^2 U [m](x,a)\Bigr) m(dx) m(da)
 \]

\item[(iii)] Note that the expressions in (ii a) continue to make sense if we merely assume that 
 \[{\rm div}_x \bigl( A_{m}(x)\bigr) \in L^1(\Mb, m) \quad \text{and} \quad {\rm Tr} \Bigl( \nabla_w^2 U [m](\cdot,\cdot)\Bigr) \in L^1(\Mb^2, m \otimes m).\]
\item[(iv)] Let $u:\Mb^k \rightarrow \Reals$ be defined for $x:=(x_1, \cdots, x_k) \in \Mb^k$ by 
 \[
u(x)=U\bigl(m_x \bigr), \quad m_x={1 \over k}\sum_{j=1}^k \delta_{x_j}.
\]
If $m=m_x$, then $u$ is differentiable in a neighborhood of $x$, $\nabla u$ is differentiable at $x$ and 
\[
\triangle_w U\bigl[m_x \bigr]=\sum_{j, l=1}^k{\rm div}_j\bigl( \nabla_{x_l} u\bigr)(x). 
\] 

\end{enumerate}  
\end{remark}
\proof{} (i),(ii) and (iii) are obvious and so, we will only present the proof of (iv). 

We fix an open ball in $\Mb^k$ centered at $a$ such that if $y$ is in the ball with twice the radius then $U$ is differentiable  at $m_y.$ We fix $x=(x_1, \cdots, x_k)$ in the ball centered at $a.$  Whenever we choose another point $y=(y_1, \cdots, y_k)$ in the ball, we will always assume that we have reordered $y$ so that 
\[
\gamma^{x, y}:={1 \over k} \sum_{j=1}^k \delta_{ (x_j, y_j)} \in \Gamma_0(\mu^x, \mu^y).
\]

(a) We have 
\begin{equation}\label{eq:sep14.2017.2} 
u(y)-u(x)=U(m_y)-U(m_x)=\int_{\Mb^2} \nabla_w U\bigl(\mu^x \bigr)(q) \cdot (p-q) \gamma^{x, y}(dq, dp)+o\bigl(W_2(m_x, m_y) \bigr), 
\end{equation}
\begin{equation}\label{eq:sep14.2017.3} 
\int_{\mathbb R^{2d}} \nabla_w U\bigl(\mu^x \bigr)(q) \cdot (p-q) \gamma^{x, y}(dq, dp)= {1 \over k} \sum_{j=1}^k  \nabla_w U(m_x)(x_j) \cdot (y_j-x_j)
\end{equation}
and 
\begin{equation}\label{eq:sep14.2017.4} 
 o\bigl(W_2(m_x, m_y) \bigr)=o(|x-y|).
\end{equation}
We combine \eqref{eq:sep14.2017.2}, \eqref{eq:sep14.2017.3} and \eqref{eq:sep14.2017.4} to obtain
\begin{equation}\label{eq:sep14.2017.1} 
\nabla_{x_j} u(x)={1 \over k} \nabla_w U[m_x](x_j)
\end{equation}

(b) Observe that if $x:=a$ for any $j \in \{1, \cdots, k\}$ 
\begin{eqnarray} 
& &\nabla_w U\bigl(m_y \bigr)(y_j)-\nabla_w U\bigl(m_x \bigr)(x_j)\nonumber\\ 
& =&\tilde A[m_x](x_j)(y_j-x_j) 
+  \int_{\Mb^2} \nabla_w^2 U(m_x)(x_j, q) \cdot (p-q) \gamma^{x, y}(dq, dp)\nonumber \\
&+&  o(|x_j-y_j|)+o\bigl(W_2(m_x, m_y) \bigr) \nonumber  
\end{eqnarray}
and so, by \eqref{eq:sep14.2017.4} 
\begin{eqnarray} 
\nabla_w U\bigl(m_y \bigr)(y_j)-\nabla_w U\bigl(m_x \bigr)(x_j)
& =&\tilde A[m_x](x_j)(y_j-x_j) \nonumber\\
&+& {1 \over k} \sum_{l=1}^k \nabla_w^2 U\bigl(m_x \bigr)(x_j, x_l) \cdot (y_l-x_l)+o\bigl(|x-y|) \bigr) \label{eq:sep14.2017.5} 
\end{eqnarray}
Thanks to (i), \eqref{eq:sep14.2017.5} implies 
\[
\nabla^2_{x_l x_j} u (x)= \left\{
\begin{array}{ll}
{1\over k} \tilde A \bigl[m_x \bigr](x_j) +{1\over k^{2}} \nabla^2_w U\bigl[\mu^x \bigr](x_j, x_j)& \mbox{if} \quad j=l,\\ \\
{1\over k^{2}} \nabla^2_w U\bigl[m_x \bigr](x_j, x_l)& \mbox{if} \quad j\not =l
\end{array}
\right.
\]
This, together with (ii) yields (iv).\endproof

\subsection{Particular case: Hessians of  functions belonging to $\cSk$}\label{subsec:particular}
Let $\Phi \in \Skt$. Set $A_m=\nabla_{x_1} \Phi$ when $k=1$ and for  $x_1 \in \Mb$ set 
\[
A_m(x_1)= \int_{\Mb^{k-1}} \nabla_{x_1} \Phi(x_1, \cdots, x_k)m(dx_2) \cdots m(dx_k)  \quad\mbox{if} \quad k\geq 2. 
\]
If $x_1, x_2 \in \Mb$ we set 
\[
A_{mm} (x_1,x_2):= \left\{
\begin{array}{ll}
0  & \mbox{if} \quad k=1,\\ 
\nabla^2_{x_2 x_1} \Phi(x_1, x_2)  & \mbox{if} \quad k=2,\\ \\
A_{mm} (x_1,x_2):=(k-1) \int_{\Mb^{k-2}}  \nabla_{x_2 x_1} \Phi(x) m(dx_3) \cdots m(dx_k)& \mbox{if} \quad k\geq 3
\end{array}
\right.
\]
For $\gamma \in \mathcal P_2(\Mb^2)$ we set 
\[
P_\gamma[m](x_1, y_1):=\nabla_{x_1} A_m(x_1) (y_1-x_1)+ \int_{\Mb^2} A_{mm} (x_1,x_2) (y_2-x_2) \gamma(dx_2, dy_2),
\]

\noindent
Then we have the following lemma.
\begin{lemma}\label{le:grad-hess1} For any  $m \in \calPtwod$, the following hold. 
\begin{enumerate} 
\item[(i)] The function $F_\Phi$ is differentiable in the sense of Wasserstein at any $m \in \calPtwod$ and  $\nabla_w F_\Phi[m]=A_m \in T_m \calPtwod.$  
\item[(ii)] 
Further assume that $\nabla^2 \Phi$ has a modulus of continuity $\rho: [0,\infty) \rightarrow [0,\infty)$ which is concave. If $\nu \in \calPtwod$, then there exists 
$\tilde{\rho}, \epsilon_2:[0,\infty) \rightarrow \Reals$ are nonnegative function (depending on $m$ and $\Phi$)  such that $\lim_{t \rightarrow 0^+} \tilde{\rho}(t)=\lim_{t \rightarrow 0^+} \epsilon_2(t) =0$ and $\tilde{\rho}$ is a concave modulus and
\[
\sup_{\gamma \in \Gamma_0(m, \nu)} \Bigl| A_\nu(y_1)-A_m(x_1)- P_\gamma[m, \nu](x_1, y_1)\Bigr|
\leq |x-y|\tilde{\rho}\bigl(|x_1-y_1|\bigr)+ W_2(m, \nu)\epsilon_2\bigl(W_2(m, \nu)\bigr) \,.
\]
\end{enumerate}  
\end{lemma}
\proof{} Note that since $\nabla^2 \Phi$ is bounded there exists a constant $C>0$ such that 
\begin{equation}\label{eq:aug16.2017.1}
|\Phi(x)| \leq C(1+|x|^2), \quad |\nabla \Phi(x)| \leq C(1+|x|), \quad |\nabla^2 \Phi(x)| \leq C \qquad \forall \quad x \in \Mb.
\end{equation} 
The second inequality in \eqref{eq:aug16.2017.1} ensures that for any $m \in \calPtwod$, $\nabla_{x_1} \Phi \in L^1(m^{\otimes(k-1)})$ and so, $A_m$ is well--defined. Furthermore, $\Phi$ is bounded. Similarly, the second and third inequalities in \eqref{eq:aug16.2017.1} ensure that $\nabla_{x_1} A_m$ and  $A_{mm}$ are  well--defined and bounded.  

(i) The proof of (i) is easier when $k=1$. We assume in the sequel that $k\geq 2.$ Using the fact that $\Phi$ is symmetric, for any $i \in \{2, \cdots, k\}$ if $\sigma$ is the permutation such that $\sigma(1)=i$, $\sigma(i)=1$ and $\sigma(j)=j$ for any $j \not \in \{1, i\}$ we have  
\begin{equation}\label{eq:aug16.2017.2}
\nabla_{x_1} \Phi(x)=\nabla_{x_i} \Phi(x^\sigma). 
\end{equation}
Applying Taylor expansion, thanks to the third inequality in \eqref{eq:aug16.2017.1} there exists a uniformly  continuous function $f: \Mb^k \times \Mb^k \rightarrow \Reals$ bounded by $C$ such that 
\begin{equation}\label{eq:aug16.2017.1.5}
\Phi(y)=\Phi(x)+ \sum_{i=1}^k \nabla_{x_i} \Phi(x) \cdot (y_i-x_i)+{1\over 2} f(x,y)|x-y|^2.
\end{equation} 

Let $m, \nu \in \calPtwod$ and let $\gamma \in \Gamma_0(m, \nu)$.  By changing variables  
\[
\int_{\Mb^2 \times \Mb^2} \nabla_{x_i} \Phi(x) \cdot (y_i-x_i) \gamma(dx_1, dy_1) \gamma(dx_i, dy_i) = \int_{\Mb^2} \nabla_{x_i} \Phi(x^\sigma) \cdot (y_1-x_1) \gamma(dx_i, dy_i) \gamma(dx_1, dy_1). 
\] 
Using \eqref{eq:aug16.2017.2} we conclude that 
\begin{equation}\label{eq:aug16.2017.3}
\int_{\Mb^2 \times \Mb^2} \Bigl(\nabla_{x_i} \Phi(x) \cdot (y_i-x_i)-\nabla_{x_1} \Phi(x) \cdot (y_1-x_1) \Bigr)\gamma(dx_1, dy_1) \gamma(dx_i, dy_i)=0. 
\end{equation}
We have  
\[
F_\Phi[\nu]-F_\Phi[m]={1\over k} \int_{\Mb^k \times \Mb^k} (\Phi(y)-\Phi(x))\gamma^{\otimes k}(dx, dy).
\]
This, together with \eqref{eq:aug16.2017.1.5} and \eqref{eq:aug16.2017.3} implies 
\begin{equation}\label{eq:aug29.2017.1}
F_\Phi[\nu]-F_\Phi[m]= \int_{\Mb^2} A_m(x_1) \cdot (y_1-x_1) \gamma(dx_1, dy_1)+  {1\over 2}\int_{\Mb^k \times \Mb^k} f(x,y)|x-y|^2\gamma^{\otimes k}(dx, dy)
\end{equation}
and so, 
\[
\Bigl|F_\Phi[\nu]-F_\Phi[m]- \int_{\Mb^2} A_m(x_1) \cdot (y_1-x_1) \gamma(dx_1, dx_2)\Bigr| \leq  {Ck \over 2}W_2^2(m, \nu).
\]
This proves that $A_m \in \partial F_\Phi[m].$ Note that $A_m$ is the gradient of 
\[x_1 \mapsto \Phi_1(x_1):=\int_{\Mb^{k-1}} \Phi(x_1, \cdots, x_k)m(dx_2) \cdots m(dx_k)\] 
which is a bounded function with bounded first derivatives. Thus, $A_m \in T_m \calPtwod.$ For any $\zeta \in T_m \calPtwod$ we have $\pi_m(\zeta-A_m)=0$, which means that $\pi_m(\zeta)=\pi_m(A_m)=A_m.$ In particular $\|\zeta\|_m \geq \|\pi_m(\zeta)\|_m= \|A_m\|_m,$ which proves that $A_m$ is the element of minimal norm in $\partial F_\Phi[m].$ 

(ii) Since the proof of (ii) is easier in the case $k=2$ compared to the case when $k\geq 3$, we assume in the sequel that $k\geq 3.$ 

Let $i \in \{3, \cdots, k\}$ and let $\sigma$ be the permutation such that $\sigma(2)=i$, $\sigma(i)=2$ and $\sigma(j)=j$ for any $j \not \in \{2, i\}.$ Given  $x=(x_1, \cdots, x_k) \in \Mb^k$, using the fact that $\Phi$ is symmetric, we have   
\begin{equation}\label{eq:aug16.2017.4}
(\nabla^2_{x_2 x_1} \Phi)(x)= (\nabla^2_{x_1 x_2} \Phi)(x) =   (\nabla^2_{x_1 x_i} \Phi)(x^\sigma)= (\nabla^2_{x_i x_1} \Phi)(x^\sigma). 
\end{equation}
For any  $y=(y_1, \cdots, y_k) \in \Mb^k$, 
\begin{equation}\label{eq:aug16.2017.5}
(\nabla_{x_1} \Phi)(y)-(\nabla_{x_1} \Phi)(x)  =\int_0^1 \sum_{i=1}^k (\nabla^2_{x_i x_1} \Phi)(x+t(y-x))(y_i-x_i) dt.
\end{equation}
Let $m, \nu \in \calPtwod$ and let $\gamma \in \Gamma_0(m, \nu)$. The change of variables which exchanges $x_2$ with $x_i$ is used to obtain  
\[
\int_{\Mb^2 \times \Mb^2} \Bigl((\nabla^2_{x_i x_1} \Phi)(x)(y_i-x_i) -(\nabla^2_{x_i x_1} \Phi)(x^\sigma)(y_2-x_2) \Bigr) \gamma(dx_2, dy_2) \gamma(dx_i, dy_i)= 0.
\]
Combining the latter with \eqref{eq:aug16.2017.4} we infer
\begin{equation}\label{eq:aug16.2017.6}
\int_{\Mb^2 \times \Mb^2} \Bigl((\nabla^2_{x_i x_1} \Phi)(x)(y_i-x_i) -(\nabla^2_{x_2 x_1} \Phi)(x)(y_2-x_2) \Bigr) \gamma(dx_2, dy_2) \gamma(dx_i, dy_i)= 0.
\end{equation}
Set 
\[
e_t(x,y):= \int_0^1 \sum_{i=1}^k\Bigl( (\nabla^2_{x_i x_1} \Phi)(x+t(y-x))- (\nabla^2_{x_i x_1} \Phi)(x) \Bigr)(y_i-x_i)  dt 
\]
By \eqref{eq:aug16.2017.5} 
\[
A_\nu(y_1)-A_m(x_1)= \int_{\Mb^{k-1} \times \Mb^{k-1} } \int_0^1 \sum_{i=1}^k (\nabla^2_{x_i x_1} \Phi)(x+t(y-x))(y_i-x_i) dt  \gamma(dx_2, dy_2) \cdots \gamma(dx_k, dy_k)
\]
This, together with \eqref{eq:aug16.2017.6} yields 
\begin{align}
A_\nu(y_1)-A_m(x_1) 
&=\nabla_{x_1} A_m(x_1) (y_1-x_1)+ \int_{\Mb^2} A_{mm} (x_1,x_2) (y_2-x_2) \gamma(dx_2, dy_2)\nonumber\\
&+  \int_{\Mb^{k-1} \times \Mb^{k-1} } \int_0^1 e_t(x,y)  dt  \gamma(dx_2, dy_2) \cdots \gamma(dx_k, dy_k).\label{eq:aug16.2017.7b}
\end{align}
Let $\rho: [0,\infty) \rightarrow [0,\infty)$ be a concave function, modulus of continuity of $\nabla^2 \Phi$. Then there exists a constant $C_k$ depending only on $k$ and $d$ such that 
\[
|e_t(x,y)| \leq C_k \sum_{i=1}^k \rho(|x_i-y_i|)|x_i-y_i| .
\]
Thus, if we  use the notation $\hat x^1=(x_2, \cdots, x_k),$ we have  
\begin{align}
 \int_{\Mb^{k-1} \times \Mb^{k-1} } \int_0^1 |e_t(x,y)|  dt  \gamma^{\otimes (k-1)}(d \hat x^1,d\hat y^1) 
&\leq  C_k \sum_{i=1}^k  \int_{\Mb^2  } \rho(|x_i-y_i|) |x_i-y_i| \gamma(d x_i,dy_i) \nonumber\\
&=  C_k \rho(|x_1-y_1|) |x_1-y_1|  \nonumber\\
&+ C_k \sum_{i=2}^k \int_{\Mb^2  } \rho(|a-b|) |a-b| \gamma(d a,db).\label{eq:aug16.2017.8}
\end{align}
Thanks to Remark \ref{rem:modulus} we conclude
\begin{align}
 \int_{\Mb^{k-1} \times \Mb^{k-1} } \int_0^1 |e_t(x,y)|  dt  \gamma^{\otimes (k-1)}(d \hat x^1,d\hat y^1) 
&\leq  C_k \rho(|x_1-y_1|) |x_1-y_1|  \nonumber\\
&+(k-1) C_k  \rho\bigl(W^{1\over 2}_2(m, \nu)\bigr) W_2(m, \nu) \Bigl( W^{1\over 2}_2(m, \nu) +1\Bigr)    \Bigr)\label{eq:aug16.2017.9}.
\end{align}
This, together with \eqref{eq:aug16.2017.7b} and \eqref{eq:aug16.2017.9}  proves (ii), after setting setting $\tilde{\rho}(t) = C_k \rho(t)$ and $\epsilon_2(t) = (k-1) C_k  \rho\bigl(t\bigr)  \Bigl( t^{1\over 2} +1\Bigr) $.\endproof

\begin{remark}\label{re:grad-hess2} Using the notation in Lemma \ref{le:grad-hess1}, we obtain the following.
\begin{enumerate} 
\item[(i)]   First 
\[ 
\nabla_{x_1}\bigl( A_{m} (x_1)\bigr) = \int_{\Mb^{k-1}} \nabla^2_{x_1x_1} \Phi(x_1, \cdots, x_k)m(dx_2) \cdots m(dx_k)
\]
is a symmetric matrix such that each row is an element of $T_m \calPtwod.$ 
\item[(ii)] Second, the rows of $A_{mm}$ are in $T_m \calPtwod$. 
\end{enumerate}  
\end{remark}

\begin{corollary}\label{co:grad-hess1} Further assume that $\nabla^2 \Phi$ has a modulus of continuity $\rho: [0,\infty) \rightarrow [0,\infty)$ which is concave. Then, 
\begin{enumerate} 
\item[(i)]  $F_\Phi$ is twice differentiable at any $m \in \calPtwod$, $\nabla_w U[m]=A_m$ and $\nabla_w^2 F_\Phi[m]=A_{mm}.$
\item[(ii)] We have $\triangle_w F_\Phi= F_{\Theta_0}.$
\item[(iii)] Similarly, if $\epsilon>0$ then  $\triangle_{w, \epsilon} F_\Phi=F_{\Theta_\epsilon}.$ Here 
\[
\Theta_\epsilon:=  \biggl((1+\epsilon)  \sum_{j=1}^k \triangle_{x_j} \Phi + \sum_{j\not = n}\sum_{l=1}^d {\partial^2 \Phi \over \partial (x_n)_l \partial (x_j)_l}\biggr).
\]
\end{enumerate}  

\end{corollary} 
\proof{} Thanks to Lemma \ref{le:grad-hess1}, we apply Theorem \ref{th:hessian1} to obtain that $U$ is twice differentiable. Remark \ref{re:Laplacian2} gives an explicit expression of $\triangle_w F_\Phi[m]$ and that of $\triangle_{w, \epsilon} F_\Phi[m].$  We use the symmetric properties of the second derivatives of $\Phi$ to obtain 
\[
 \int_{\Mb^k} \triangle_{x_1} \Phi m(dx_1) \cdots m(dx_k)=  {1\over k} \int_{\Mb^k} \sum_{j=1}^k \triangle_{x_j} \Phi m(dx_1) \cdots m(dx_k)
\]
and 
\[
 \int_{\Mb^k}\sum_{n=1}^d {\partial^2 \Phi \over \partial (x_2)_n \partial (x_1)_n}  m(dx_1) \cdots m(dx_k)=  
 {1 \over k(k-1)} \int_{\Mb^k} \sum_{n=1}^d \sum_{j\not= l} {\partial^2 \Phi \over \partial (x_j)_n \partial (x_l)_n} m(dx_1) \cdots m(dx_k)
\]
to conclude the proof of the Corollary.\endproof

\subsection{Convergence theorem for the Wasserstein Hessian} In this subsection $\mu \in \calPtwod$, $C_\mu>0$  and  $\mathcal O$, an open ball centered at $\mu \in \calPtwod$. Suppose $G_N: \calPtwod \rightarrow \Reals$ is a sequence of continuous functions converging uniformly to  $G: \calPtwod \rightarrow \Reals$ and $G_N$ is twice differentiable (cf. Definition \ref{defn:Laplacian}) on $\mathcal O$. Suppose 
$$ 
\nabla_w G_N[m] \in C(\Mb; \mathbb R^{d}), \quad \nabla\bigl(\nabla_w G_N\bigr)[m]) \in C(\Mb; \mathbb R^{d \times d}), \quad \nabla^2_w G_N[m]) \in C(\Mb^2; \mathbb R^{d \times d})
$$ 
for any $m \in \mathcal O$ and  $|\nabla_w G_N[m](x)| \leq C_\mu(1+|x|)$ for any $(N, x, m) \in \mathbb N \times \Mb \times \mathcal O.$ Suppose $(x, m) \mapsto \nabla\bigl(\nabla_w G_N\bigr)[m](x)$ is continuous bounded on $\Mb \times \mathcal O$ and  $(x, y, m) \mapsto \nabla^2_w G_N[m](x,y)$ is continuous and bounded on $\Mb^2 \times \mathcal O$ independently on $N.$ Suppose  
\begin{equation}\label{eq:infdef1bis} 
\sup_{\gamma\in\Gamma_0(m, \nu)} \Bigl| G_N[\nu]-G_N[m] -\int_{\Mb^2} \nabla_w G_N[m](x) \cdot (y-x) \gamma(dx, dx)\Bigr|  \leq C_m W_2^2(m,\nu),
\end{equation}
for any $m, \nu \in \mathcal O$, and 
\begin{equation}\label{eq:infdef1ter} 
\sup_{\gamma \in \Gamma_0(m, \nu)}\Bigl|\nabla_w G_N[\nu](y)-\nabla_w G_N[m](x) - P_\gamma^N[m](x, y)\Bigr| \leq C_m \bigl(|x-y|^2+ W_2^2(m, \nu) \bigr)
\end{equation} 
Here, for $\gamma \in \mathcal P(\Mb^2)$ and $x, y \in \Mb,$ we have set 
\[
P_\gamma^N[m](x, y):= \nabla\bigl(\nabla_w G_N\bigr)[m](x) (y-x)+ \int_{\Mb^2} \nabla^2_w G_N[m] (x,a) (b-a) \gamma(da, db),
\]

\begin{theorem}\label{th:grad-hess-series} Suppose \eqref{eq:infdef1bis} and \eqref{eq:infdef1ter} hold and 
\begin{enumerate} 
\item[(a)] $\bigl(\nabla_w G_N \bigr)_N$ converges uniformly on $\Mb \times \mathcal O$ to $(x, m) \mapsto A_m(x), $  
\item[(b)] $\bigl(\nabla\bigl(\nabla_w G_N\bigr)\bigr)_N$ converges uniformly on $\Mb \times \mathcal O$  to $ (x, m) \mapsto \tilde A[m](x)$   
\item[(c)] $\bigl( \nabla^2_w G_N \bigr)_N$ converges uniformly on $\Mb^2 \times \mathcal O$  to  $(x, y, m) \mapsto A_{m m}(x,y).$
\end{enumerate}  
Then  
\begin{enumerate} 
\item[(i)] $G$ is differentiable on $\mathcal O$ and $A= \nabla_w G.$
\item[(ii)] $G$ is twice continuously differentiable  on $\mathcal O$, $\tilde A= \nabla\bigl(\nabla_w G\bigr)$ and $A_{m m}\equiv \nabla^2_w G[m].$ 
\end{enumerate}  

\end{theorem} 
\proof{}  Note that $(x, m) \mapsto A_m(x),$ $\tilde A[m](x)$ and  $(x, y, m) \mapsto A_{m m}(x,y)$ are continuous and the latter two functions are bounded as limits of bounded functions. Let $m, \nu \in \mathcal O$. 

(i) For any $\gamma\in\Gamma_0(m, \nu)$ we have  
\[
\Bigl| \int_{\Mb^2} \Bigl(\nabla_w G_N[m](x)- A_m(x) \bigr) \cdot (y-x) \gamma(dx, dx)\Bigr| \leq \Big\| \nabla_w G_N[m]- A_m\Big\|_{L^\infty} W_2(m, \nu)
\]
and so, letting $N$ tend to $\infty$ in \eqref{eq:infdef1bis}, we obtain 
\[ 
 \Bigl| G [\nu]-G [m] -\int_{\Mb^2} A_{m}(x) \cdot (y-x) \gamma(dx, dx)\Bigr|  \leq C_m W_2^2(m,\nu).
\] 
Since $\gamma\in\Gamma_0(m, \nu)$ is arbitrary, we conclude $A_{m} \in \partial G[m]$.  Observe that since $|\nabla_w G[m](x)| \leq C_m(1+|x|)$ for any $(x, m) \in \Mb \times \mathcal O,$  and as $\nabla_w G_N[m] \in T_m \calPtwod$, we conclude that $A_{m} \in T_m \calPtwod$ and so, $A_m= \nabla_w G[m].$  

(ii) Since $\bigl(\nabla\bigl(\nabla_w G_N\bigr)\bigr)_N$ converges uniformly to  $\tilde A$, we have $\tilde A= \nabla\bigl(\nabla_w G\bigr).$ Observe if $\gamma\in\Gamma_0(m, \nu)$ and $x_1, y_1 \in \Mb$ then 
\[
\Bigl| \int_{\Mb^2}\Bigl(\bigl( \nabla_w^2 G_N[m] - A_{m m}\bigr)(x_1, x_2)(y_2-x_2)\Bigr) \gamma(dx_2, dy_2) \Bigr| \; \leq \; \Big\| \nabla_w G_N[m]- A_{m m}\Big\|_{L^\infty} W_2(m, \nu).
\] 
As above, we conclude that 

\[
\Bigl|\nabla_w G[\nu](y)-\nabla_w G[m](x) - P_\gamma[m](x, y)\Bigr| \leq C_m \bigl(|x-y|^2+ W_2^2(m, \nu) \bigr)
\]
where 
\[
P_\gamma[m](x, y):= \tilde A[m](x) (y-x)+ \int_{\Mb^2} A_{m m} (x,a) (b-a) \gamma(da, db).
\]
Since $\gamma\in\Gamma_0(m, \nu)$ is arbitrary, we apply Theorem \ref{th:hessian1} to obtain that $G$ is twice differentiable at $m$, $\tilde A[m]=\nabla\bigl(\nabla_w G\bigr)[m]$ and $A_{m m}= \nabla^2_w G[m].$ The identity $\pi_m\bigl(  \nabla^2_w G_N[m]\bigr)= \nabla^2_w G_N[m]$ implies $A_{m m}=\pi_m\bigl( A_{m m}\bigr)$.  We conclude the proof of (ii) by setting $\rho(t)= C_m t$ and $\epsilon_2(t) = C_m t$. \endproof

\section{Fourier transform and expansions} 

\subsection{Polynomial eigenfunctions of the Laplacian operator; Harmonic functions} Let $\Phi_\xi^k \in \Sktcomplex$ be defined by 
\[
\Phi^k_{\xi}(x):= {1\over k!} \sum_{\sigma \in P_k} \exp\Bigl(-2\pi i \sum_{j=1}^k \langle \xi_{\sigma (j)}, x_j\rangle \Bigr), \quad \forall \, \xi=(\xi_1, \cdots, \xi_k) \in \Mb^k.
\]
The function $\Phi^k_{\xi}$ is obtain as the symmetrization of $x \mapsto \exp\Bigl(-2\pi i \sum_{j=1}^k \langle \xi_{j}, x_j\rangle \Bigr).$ 

Let  $F^k_\xi: \calPtwod \rightarrow \mathbb C$  be  $F_{\Phi^k_{\xi}}.$ In other words, 
\[
F^k_\xi[m]= {1 \over k} \int_{\Mb^k} \Phi^k_{\xi}(x)m(dx_1) \cdots m(dx_k)
\]
Set 
\begin{equation}\label{eq:lambdak}
\lambda^2_k(\xi):= 4\pi^2 \Big|  \sum_{j=1}^k \xi_j\Big|^2 \quad \text{and}\quad \widehat m(\xi_j):= \int_{\Mb} \exp\bigl( -2\pi \langle \xi_j, z\rangle \bigr) m(dz).
\end{equation} 
Note $\widehat m$ is the Fourier transform of $m.$

\begin{lemma}\label{le:computationFk} The following hold for any $m \in \calPtwod:$   
\begin{enumerate} 
\item[(i)]  We have $F^k_\xi[m]= k^{-1}\;\; \prod_{j=1}^k  \widehat m(\xi_j)$. 
\item[(ii)] We have $\triangle_w F^k_\xi[m]= -\lambda^2_k(\xi) F^k_\xi[m]$. 
\item[(iii)] We have $\triangle_{w,\epsilon} F^k_\xi[m]= -\Bigl(\lambda^2_k(\xi) +4\pi^2 \epsilon \sum_{j=1}^n |\xi_j|^2\Bigr) F^k_\xi[m]$.
\end{enumerate}  
\end{lemma} 
\proof{} (i) It suffices to observe that for any $\sigma \in P_k$, we have 
\[
\int_{\Mb^k } \exp\Bigl(-2\pi i \sum_{j=1}^k \langle \xi_{\sigma (j)}, x_j\rangle \Bigr) m^{\otimes k}(dx) =  \int_{\Mb^k } \exp\Bigl(-2\pi i \sum_{j=1}^k \langle \xi_{j}, x_j\rangle \Bigr)  m^{\otimes k}(dx) = \prod_{j=1}^k  \widehat m(\xi_j).
\]

(ii) Direct computations reveal that  
\begin{align}
\nabla_{x_1} \Phi^k_{\xi}(x) 
&={-2\pi i \over k!}  \sum_{\sigma \in P_k} \xi_{\sigma(1) }\exp\Bigl(-2\pi i \sum_{j=1}^k \langle \xi_{\sigma (j)}, x_j\rangle \Bigr)     \label{eq:sep01.2017.1}\\
\nabla^2_{x_1 x_1} \Phi^k_{\xi}(x)  &=  
{-4\pi^2 \over k!}  \sum_{\sigma \in P_k} \xi_{\sigma(1) } \otimes \xi_{\sigma(1) } \exp\Bigl(-2\pi i \sum_{j=1}^k \langle \xi_{\sigma (j)}, x_j\rangle \Bigr) \label{eq:sep01.2017.2}\\
 \nabla^2_{x_2 x_1} \Phi^k_{\xi}(x)  &=  
{-4\pi^2 \over k!}  \sum_{\sigma \in P_k} \xi_{\sigma(1) } \otimes \xi_{\sigma(2) } \exp\Bigl(-2\pi i \sum_{j=1}^k \langle \xi_{\sigma (j)}, x_j\rangle \Bigr) \label{eq:sep01.2017.3}
\end{align}
Thus, since $m(dx_1) \cdots m(dx_k)$ is permutation invariant, the expression  
\[
\int_{\Mb^k} {\rm Tr}\Bigl(\nabla^2_{x_2 x_1} \Phi^k_{\xi}(x) \Bigr)  m(dx_1) \cdots m(dx_k) 
\] 
equals  
\begin{align}
 &= 
-{4\pi^2 \over k!} \sum_{\sigma \in P_k} \langle \xi_{\sigma(1)},\xi_{\sigma(2)} \rangle \int_{\Mb^k}  \exp\Bigl(-2\pi i \sum_{j=1}^k \langle \xi_{j}, x_j\rangle \Bigr) m(dx_1) \cdots m(dx_k) \nonumber\\ 
&= -{4\pi^2(k-2)! \over k! } \sum_{j\not =n}^k \langle \xi_{j },  \xi_{n } \rangle  \int_{\Mb^k}  \exp\Bigl(-2\pi i \sum_{j=1}^k \langle \xi_{j}, x_j\rangle \Bigr) m(dx_1) \cdots m(dx_k).\nonumber
\end{align}
This reads off 
\begin{equation} \label{eq:aug24.2017.4}
\int_{\Mb^k} {\rm Tr}\Bigl(\nabla^2_{x_2 x_1} \Phi^k_{\xi}(x) \Bigr)  m(dx_1) \cdots m(dx_k)= {-4\pi^2  \over k-1}  \sum_{j\not =n}^k \langle \xi_{j },  \xi_{n } \rangle F^k_\xi(m).
\end{equation}
Similarly, 
\begin{equation} \label{eq:aug24.2017.5}
\int_{\Mb^k} \triangle_{x_1} \Phi^k_{\xi}(x) m(dx_1) \cdots m(dx_k)  =  -4\pi^2 \sum_{j=1}^k |\xi_j|^2 F^k_\xi(m).
\end{equation}
Combining this with \eqref{eq:aug24.2017.4}, thanks to Corollary \ref{co:grad-hess1}, we complete the proof of (ii).  Similarly, we obtain the proof of (iii). \endproof

\begin{remark}\label{re:eigen} The followings hold.
\begin{enumerate} 
\item[(i)]  If $\lambda_k(\xi)=0$ then $F^k_\xi$ belongs to the kernel of $\triangle_w.$ We say that  $F^k_\xi$  is a harmonic function.
\item[(ii)] ({\bf non smooth harmonic functions}) Let $f \in {\rm Sym_2[\Mb]}$ be an even function and set $\Phi(x,y)= f(x-y).$ Note that ${\rm div}_x (\nabla_x \Phi)= {\rm Tr}(\nabla_{yx} \Phi)$. Using Corollary \ref{co:grad-hess1}, we conclude that  $\triangle_w F_\Phi[m] \equiv 0$ and so, $U$ is a harmonic function. Starting with $f \not \in C^3(\Mb)$, we obtain that $\Phi \not\in C^3(\Mb)$ and so, $F_\Phi$ is a harmonic function which  is not  regular up to the third order.  However, if $\epsilon>0$, $\triangle_{w, \epsilon} F_\Phi \not \equiv 0$ and $\triangle_{w, \epsilon}$ has a smoothing effect.
\item[(iii)] A direct consequence of (ii) is that $\triangle_w$ is not a smoothing operator (except on $H^s\bigl(\calPtwod\bigr)$: cf. Theorem \ref{th:superposition}). 
\end{enumerate}  
\end{remark}

\subsection{$H^s$--spaces and spaces of Fourier transforms}  Throughout this subsection, $s \geq 0$ is a real number.

\begin{definition}\label{defnHsbis} Let $\lambda_k$ be the function defined in \eqref{eq:lambdak}. 
\begin{enumerate}
\item[(i)] 
We call $\mathcal A$ the set of sequences of functions $(a_k)_{k=1}^\infty$ such that $a_k: \Mb^k \rightarrow \mathbb C$ is a Borel function that are symmetric in the sense that $a_k(\xi)=a_k(\xi^\sigma)$ for any $\xi \in \Mb^k$ and any $\sigma \in P_k.$ In other words, $a_k$ is defined on $\Mb^k/P_k$, the $k$--symmetric product of $\Mb.$ 
\item[(ii)] We call $(a_k)_{k=1}^\infty \in \mathcal A$ the Fourier transform of $F:\calPtwod \rightarrow \Reals,$ if there exist $\Phi_k \in \cSk \cap L^2(\Mb^k)$ is such that the series 
\[
 \sum_{k=1}^\infty {1 \over k! \; k} \int_{\Mb^k}  \Phi_k(x) m(dx_1) \cdots m(dx_k)
\]
converges to $F$ and $\widehat a_k=  \Phi_k$. 
\end{enumerate}  
\end{definition}

\begin{definition}\label{defnHs} We have the following definition.  
\begin{enumerate}
\item[(i)] We call $\mathcal A^s$ the set of sequences $A:=\bigl( a_k\bigr)_{k=1}^\infty\subset \mathcal A$ such that 
\begin{equation}\label{eq:aug23.2017.1}
||A||^2_{H^s}:=\sum_{k=1}^\infty {1\over k! }\int_{\Mb^k } |a_k(\xi)|^2 \bigl( 1+ \lambda^2_k(\xi)\bigr)^s d\xi<\infty.
\end{equation} 
\item[(ii)] If $B=\bigl( b_k\bigr)_{k=1}^\infty \in \mathcal A^s,$ then the following is a well--defined sesquilinear form (cf. Lemma \ref{de:sesquilinear}): 
\[
\langle A; B \rangle_{{H^s}}:= \sum_{k=1}^\infty {1\over k! }\int_{\Mb^k } a_k(\xi) b^*_k(\xi) \bigl( 1+ \lambda^2_k(\xi)\bigr)^s d\xi
\]
\end{enumerate}  
\end{definition}

\begin{lemma}\label{de:sesquilinear} The sesquilinear form $\langle \cdot; \cdot\rangle_{H^s}: \mathcal A^s  \times \mathcal A^s \rightarrow \mathbb C$ is well defined.
\end{lemma}
\proof{} Let $A, B \in \mathcal A^s$ be as in Definition \ref{defnHs}.  Then for any $\lambda>0$ we have 
\begin{equation}\label{eq:aug22.2017.3}
{1\over k! } \Bigl|\int_{\Mb^k } a_k(\xi)  b_k^*(\xi) \bigl( 1+ \lambda^2_k(\xi)\bigr)^s d\xi \Bigr| \leq {1 \over 2 k!} \int_{\Mb^k } \biggl({|a_k(\xi)|^2 \over \lambda^2} +  \lambda^2 | b_k(\xi)|^2\biggr)  \bigl( 1+ \lambda^2_k(\xi)\bigr)^s d\xi. 
\end{equation} 
Therefore, the series produced by the left hand side of \eqref{eq:aug22.2017.3} converges absolutely, which concludes the proof. \endproof

\begin{lemma}\label{le:poisson} Suppose $A:=\bigl( a_k\bigr)_{k=1}^\infty$ and $B=\bigl( b_k\bigr)_{k=1}^\infty$ belong to $\mathcal A^s.$ Then 
\begin{enumerate} 
\item[(i)] (H\"older's inequality) 
\[
|\langle A; B\rangle_{H^s}| \leq ||A||_{H^s} \cdot ||B||_{H^s}
\]
\item[(ii)] (triangle inequality) 
\[
||A+B||_{H^s} \leq ||A||_{H^s} + ||B||_{H^s}
\]
\end{enumerate}  
\end{lemma} 
\proof{} Assume without loss of generality that $ ||B||_{H^0} \not =0.$ 

(i) By \eqref{eq:aug22.2017.3} 
\[
2 |\langle A; B\rangle_{H^s}| \leq { ||A||^2_{H^s} \over \lambda^2} + \lambda^2  ||B||^2_{H^s}.
\]
We use $\lambda:=||A||_{H^s}^{1\over 2}  ||B||_{H^s}^{-{1\over 2}}$ to conclude the proof of (i).

(ii) We use (i) and the  identity 
\[
||A+B||^2_{H^s}=||A||^2_{H^s} + ||B||^2_{H^s}+ \langle A; B\rangle_{H^s} +\langle G; F\rangle_{H^s}
\]
to conclude the proof of (ii). \endproof

\begin{lemma}\label{le:Hs} Let $\Phi \in \Sk \cap L^2(\Mb^k)$ and let $\bigl( a_k\bigr)_{k=1}^\infty, \bigl( b_k\bigr)_{k=1}^\infty \in \mathcal A$ be such that 
\[
\sum_{k=1}^\infty {1\over k! }\int_{\Mb^k } \Bigl( |a_k(\xi)|^2 + |b_k(\xi)|^2 \Bigr) d\xi<\infty.
\]
\begin{enumerate}
\item[(i)] Since the Fourier transform is an isometry of $L^2(\Mb^k; \mathbb C)$ onto $L^2(\Mb^k;\mathbb C)$, we obtain  $a:=\check \Phi \in L^2(\Mb^k; \mathbb C)$. One check that $a$ is symmetric and so, if we further assume that $\Phi \in   L^1(\Mb^k)$ then $a \in \cSkcomplex \cap L^\infty(\Mb^k; \mathbb C).$  
\item[(ii)] Observe   for any $k\geq 1$, $a_k, b_k \in L^2(\Mb^k; \mathbb C)$. Let $\Phi_k, \Psi_k \in \cSkcomplex \cap L^2(\Mb^k; \mathbb C)$ be such that $a_k= \check \Phi_k$ and $b_k= \check\Psi_k.$  For any $N \geq 1$,   
\[
\sum_{k=1}^N {1 \over k!} \int_{\Mb^k } a_k(\xi) b_k^*(\xi) d \xi= \sum_{k=1}^N {1 \over k!} \int_{\Mb^k}  \Phi_k(x) \Psi_k^*(x) dx.
\]
\item[(iv)] Further assume for any integer $k \geq 1$,  $a_k \in L^1(\Mb^k; \mathbb C).$ Then for any $N\geq 1$ 
\[
\sum_{k=1}^N {1 \over k!} \int_{\Mb^k } a_k(\xi) F^k_\xi[m] d \xi= \sum_{k=1}^N {1 \over k! \;\; k} \int_{\Mb^k}  \Phi_k(x) m(dx_1) \cdots m(dx_k)
\]
The series converge uniformly on $\calPtwod$ if there exist  constant $C, \delta>$ independent of $m$ and $k$ such that 
\begin{equation}\label{eq:uniform-bound-ak}
\int_{\Mb^k}|a_k(\xi)| d\xi \leq {C k! \over k^{\delta}}
\end{equation}
\end{enumerate}  
\end{lemma}
\proof{}   (i) is straightforward to check. (ii) is a consequence of Plancherel's theorem and the fact that  the Fourier transform is an isometry of $L^2(\Mb^k; \mathbb C)$.

(iii) Since $\Phi_k=\widehat{a}_k$, when $a_k \in L^1(\Mb^k; \mathbb C),$ we may use Fubini's theorem to obtain 
\begin{align}
\int_{\Mb^k} \Phi_k(x) m(dx_1) \cdots m(dx_k)
&= \int_{\Mb^{k}} \ m(dx_1) \cdots m(dx_k)\int_{\Mb^k} a_k(\xi_1, \cdots, \xi_k) \exp\Bigl({-2\pi i \sum_{j=1}^k \langle \xi_j,x_j \rangle}\Bigr) d\xi \nonumber\\
&=  \int_{\Mb^{k}} a_k(\xi_1, \cdots, \xi_k) d\xi \int_{\Mb^k}  \exp\Bigl({-2\pi i \sum_{j=1}^k \langle \xi_j,x_j \rangle}\Bigr) m(dx_1) \cdots m(dx_k) \nonumber\\
&=k \int_{\Mb^{k}} a_k(\xi_1, \cdots, \xi_k) F_{\xi}[m]d\xi.\label{eq:aug23.2017.2}
\end{align}
By the fact that $|F^k_\xi[m]|\leq k^{-1},$ we have  
\[
\Big| \sum_{k=1}^N {1 \over k!} \int_{\Mb^k } a_k(\xi) F^k_\xi[m] d \xi \Big| \leq   \sum_{k=1}^N {1 \over k! \; k} \int_{\Mb^k } |a_k(\xi) | d \xi   \leq \sum_{k=1}^N  {C \over k^{1+\delta}}.  
\]
This concludes the proof.  \endproof

\begin{definition}\label{de:defncurlyHs} We have the following definition.  
\item[(i)] We call $H^s(\calPtwod)$ the set of $F: \calPtwod \rightarrow [-\infty, \infty]$ for which there exist $(a_k)_{k=1}^\infty\subset \mathcal A^s$ such that     
\[
\sum_{k=1}^\infty {1\over k! }\int_{\Mb^k } a_k(\xi) F^k_\xi[m] d \xi
\]
converges to $F[m]$ for any $m \in \calPtwod$. 
\item[(ii)] We define $\mathcal{H}^s(\calPtwod) $ the set of $F: \calPtwod \rightarrow [-\infty, \infty]$ for which there exist $(a_k)_{k=1}^\infty\subset \mathcal A^s$, $\delta, C>0$ such that \eqref{eq:uniform-bound-ak} holds for all $k$ natural number and  
\begin{equation}\label{eq:sep24.2017.6}
\sum_{k=1}^\infty {1\over k! }\int_{\Mb^k } a_k(\xi) F^k_\xi[m] d \xi 
\end{equation}  
converges to $F[m]$ for any $m \in \calPtwod$.  Thanks to Remark \ref{re:F_lambda_m}, the following definition is meaningful. 
\end{definition}

From definition, we have
\begin{eqnarray*}
\mathcal{H}^s(\calPtwod)  \subset H^s(\calPtwod) \cap C(\calPtwod)
\end{eqnarray*}
and the second inclusion results from the fact that the convergence of the series converges uniformly on $m \in \calPtwod$ (cf. Lemma \ref{le:Hs} (iv)).

\subsection{Integrations by parts; Hessians in terms of Fourier transforms}\label{subsection:super-grad-hess} Throughout this subsection  $s, \delta>0$ , $\epsilon\geq 0$ and $\bigl( a_k
\bigr)_{k=1}^\infty \in \mathcal A^s$. When needed we shall make various  assumptions such as 
\begin{equation}\label{eq:sep01.2017.4} 
\int_{\Mb^k} |a_k(\xi)| \cdot |\xi_1| d\xi, \;   \int_{\Mb^k} |a_k(\xi)| \cdot | \xi_1|^2 d\xi  \; \leq {C k! \over k^{1+ \delta} }, 
\end{equation}
\begin{equation}\label{eq:sep01.2017.5} 
   \int_{\Mb^k} |a_k(\xi)| \cdot | \xi_1| \cdot  | \xi_2| d\xi  \quad \leq {C k! \over k^{2+\delta} }.
\end{equation}
or 
\begin{equation}\label{eq:sep01.2017.5bis} 
 \int_{\Mb^k} |a_k(\xi)| \cdot |\xi_1|^3 d\xi \leq {Ck! \over k^{1+\delta} }, \quad   \int_{\Mb^k} |a_k(\xi)| \cdot | \xi_1|^2 \cdot  | \xi_2| d\xi  \quad \leq {Ck! \over k^{3+\delta} }.
\end{equation}
When \eqref{eq:uniform-bound-ak}  is in force then the series  
\begin{equation}\label{eq:U-0}
U_0[m]:=\sum_{k=1}^\infty {1 \over k!} \int_{\Mb^k } a_k(\xi) F^k_\xi[m] d \xi, \qquad m \in \calPtwod
\end{equation}
converges uniformly  (cf. Lemma \ref{le:Hs} (iv)).
\begin{corollary}\label{co:superposition2} Assume \eqref{eq:uniform-bound-ak}, \eqref{eq:sep01.2017.4} and \eqref{eq:sep01.2017.5}  hold. 
\begin{enumerate} 
\item[(i)] Then $U_0$ is continuously differentiable on $\calPtwod$, and using the notation $\langle \xi, x \rangle$ in place of  $\sum_{j=1}^k \langle \xi_{j}, x_j \rangle$, we have
\begin{equation}\label{eq:sep30.2017.7.5}
\nabla_w U_0[m](x_1) \equiv  \sum_{k=1}^\infty  {-2 \pi  i \over k!} \int_{\Mb^{k-1} \times \Mb^k } a_k(\xi) \xi_1  e^{-2\pi i \langle \xi, x \rangle}  m(dx_2) \cdots m(dx_k) d \xi
\end{equation}
\item[(ii)] If we further assume \eqref{eq:sep01.2017.5bis} holds,  then $U_0$ is twice continuously differentiable on $\calPtwod.$  We have 
\begin{equation}\label{eq:sep30.2017.7}
 \nabla \bigl( \nabla_w U_0[m](x_1) \bigr)  \equiv  \sum_{k=1}^\infty  {- 4\pi^2  \over k!}  \int_{\Mb^{k-1} \times \Mb^k } a_k(\xi)\xi_1 \otimes \xi_1  e^{-2\pi i \langle \xi, x \rangle} m(dx_2) \cdots m(dx_k) d \xi
\end{equation}
and 
\begin{equation}\label{eq:sep30.2017.8}
 \nabla^2_w U_0[m](x_1, x_2) \bigr)\equiv  \sum_{k=1}^\infty  {- 4\pi^2(k-1)  \over k!} \int_{\Mb^{k-2} \times \Mb^k } a_k(\xi) \, \xi_1 \otimes \xi_2  e^{-2\pi i \langle \xi, x \rangle} m(dx_3) \cdots m(dx_k) d \xi.
\end{equation}
\item[(iii)] Under the same assumptions as in (ii), 
\[
\triangle_{w,\epsilon} U_0[m]= \sum_{k=1}^\infty {1 \over k!} \int_{\Mb^k } a_k(\xi) \triangle_{w,\epsilon}(F^k_\xi[m]) d \xi.
\]
\end{enumerate}  
\end{corollary} 
\proof{} (i) Let $ m \in \calPtwod$ and set 
\[
G_N[m]= \sum_{k=1}^N {1 \over k!} \int_{\Mb^k } a_k(\xi) F^k_\xi[m] d \xi= \sum_{k=1}^N {1 \over k!} {1 \over k} \int_{\Mb^k } \Phi_k(x) m(dx_1) \cdots m(dx_k), 
\] 
where $\Phi_k=\widehat a_k.$
Observe 
\[
\nabla_{x_1} \Phi_k(x)= -2 \pi i  \int_{\Mb^k } a_k(\xi) \xi_1 \exp\Bigl(-2\pi i \sum_{j=1}^k \langle \xi_{j}, x_j \rangle \Bigr) d \xi,
\]
Thus, by Lemma \eqref{le:grad-hess1} and the linearity of the Wasserstein gradient, 
\[
\nabla_w G_N[m](x_1)= -2 \pi i \sum_{k=1}^N g_k^0[m](x_1)
\]
where 
\[
g_k^0[m](x_1):=  {1 \over k!} \int_{\Mb^{k-1} \times \Mb^k } a_k(\xi) \xi_1 \exp\bigl(-2\pi i \sum_{j=1}^k \langle \xi_{j}, x_j \rangle \bigr)  m(dx_2) \cdots m(dx_k) d \xi.
\]
If $\nu \in \calPtwod$ and $\gamma \in \Gamma_0(m, \nu)$, the first order expansion of $t \rightarrow e^{-2\pi t}$ yields the first order Taylor expansion of $F_{\Phi_k}$ around $m$,   given by 
\[
\Big|F_{\Phi_k}[\nu] -F_{\Phi_k}[m]- \int_{\Mb^2} \langle \nabla_w F_{\Phi_k}[m](x_1) ; y_1-x_1 \rangle \gamma(dx_1, dy_1)\Big| \leq C[k] W_2^2(m, \nu).
\]
Here 
\[
C[k]:= {2\pi^2 \over  k} \sum_{j, l=1}^k \int_{\Mb^k} |\xi_j| |\xi_l| |a_k(\xi)|d\xi
\]
Since $a_k$ is symmetric,  
\begin{equation}\label{eq:sep29.2017.6}
\sum_{j, l=1}^k \int_{\Mb^k} |\xi_j| |\xi_l| |a_k(\xi)|d\xi=k \int_{\Mb^k} |\xi_1|^2 |a_k(\xi)|d\xi +k(k-1) \int_{\Mb^k} |\xi_1| |\xi_2| |a_k(\xi)|d\xi.
\end{equation}
We combine \eqref{eq:uniform-bound-ak}, \eqref{eq:sep01.2017.4} and \eqref{eq:sep01.2017.5} and use that $k-1 \leq k$ to conclude that 
\[
{C[k] \over k!} \leq {4 \pi^2 C \over k^{1+\delta}}. 
\] 
Thus, by the above first order Taylor expansion of $F_{\Phi_k}$ around $m$ we obtain 
\begin{equation}\label{eq:sep29.2017.1}
\Big|G_N[\nu] -G_N[m]- \int_{\Mb^2} \langle \nabla_w G_N[m](x_1) ; y_1-x_1 \rangle \gamma(dx_1, dy_1)\Big| \leq C_m W_2^2(m, \nu) 
\end{equation}
where 
\[
C_m:= \sum_{k=1}^N {1 \over k! }  C[k]  <\infty.     
\]
But $(x_1, m) \mapsto g_k^0[m](x_1)$ are continuous functions such that 
\[
\Bigl| g_k^0[m](x_1) \Bigr| \leq {1 \over k!}  \int_{\Mb^k } |a_k(\xi) \xi_1| d\xi.
\]
Thanks to \eqref{eq:uniform-bound-ak} and \eqref{eq:sep01.2017.4}  we obtain that the series $\bigl(-2 \pi i \sum_{k=1}^N g_k^0[m](x_1)\bigr)_N$ is a Cauchy sequence for the uniform convergence and so,  it converges uniformly to a continuous function given by the function at the right handside of \eqref{eq:sep30.2017.7.5}, which we denote as $A.$ We let $N$ tend to $\infty$ in \eqref{eq:sep29.2017.1} to conclude $A \equiv \nabla_w U_0$ and conclude the proof of (i).

For any $n \in \{1, \cdots, k\}$ we have    
\begin{equation}\label{eq:aug31.2017.1}
\nabla_{x_n x_1} \Phi_k(x)= -4 \pi^2  \int_{\Mb^k } a_k(\xi)\xi_1 \otimes \xi_n \exp\Bigl(-2\pi i \sum_{j=1}^k \langle \xi_{j}, x_j \rangle \Bigr) d \xi,
\end{equation} 
Hence by Remark \ref{re:grad-hess2} 
\begin{equation}\label{eq:sep30.2017.2}
\nabla \bigl( \nabla_w G_N[m](x_1) \bigr)= -4 \pi^2 \sum_{k=1}^N g_k^1[m](x_1)
\end{equation} 
where 
\[
g_k^1[m](x_1) :=  {1 \over k!} \int_{\Mb^{k-1} \times \Mb^k } a_k(\xi)\xi_1 \otimes \xi_1 \exp\Bigl(-2\pi i \sum_{j=1}^k \langle \xi_{j}, x_j \rangle \Bigr) m(dx_2) \cdots m(dx_k) d \xi
\]
Note, 
\[
\Bigl| g_k^1[m](x_1) \Bigr| \leq {1 \over k!}  \int_{\Mb^k } |a_k(\xi)| \;  |\xi_1|^2 d\xi.
\]
Thanks to \eqref{eq:uniform-bound-ak} and \eqref{eq:sep01.2017.4} again, we obtain that the series $\bigl(-4\pi^2\sum_{k=1}^N g_k^1[m](x_1)\bigr)_N$ is a Cauchy sequence for the uniform convergence and so,  it converges uniformly to the continuous function $\tilde A$ at the right handside of \eqref{eq:sep30.2017.7}. We will soon see that it is legitimate to denote this limit as $\bigl( \nabla_w U_0[m](x_1) \bigr).$

By Lemma \ref{le:grad-hess1} and the linearity of $\nabla^2_w$,  
\begin{equation}\label{eq:sep30.2017.3}
 \nabla^2_w G_N[m](x_1, x_2) \bigr)= -4 \pi^2 \sum_{k=2}^N g_k^2[m](x_1, x_2)
\end{equation} 
where 
\[
g_k^2[m](x_1, x_2) :=  {k-1 \over k!} \int_{\Mb^{k-2} \times \Mb^k } a_k(\xi) \, \xi_1 \otimes \xi_2 \exp\Bigl(-2\pi i \sum_{j=1}^k \langle \xi_{j}, x_j \rangle \Bigr) m(dx_3) \cdots m(dx_k) d \xi.
\]
The first order expansion of $t \rightarrow e^{-2\pi t}$ yields the first order Taylor expansion of $A_k :=\nabla_w F_{\Phi_k}$ around $(x_1, m)$ given by 
\begin{align}
A_k[\nu](y_1) &= A_k[m](x_1)- 4 \pi^2 k!  g_k^1[m](x_1)(y_1-x_1) \nonumber\\
& -4 \pi^2 k!  \int_{\Mb^2} g^2_k[m]^1(x_1, x_2)(y_2-x_2) \gamma(dx_2, dy_2) +B_k \nonumber
\end{align} 
Here, the remainder $B_k$ is such that 
\[
|B_k| \leq 4 \pi^2 \sum_{j, l=1}^k \int_{\Mb^{2(k-1)} \times \Mb^k } |a_k(\xi)| |\xi_l| |\xi_j| |y_l-x_l| |y_j-x_j|  \gamma(dx_2, dy_2) \cdots \gamma(dx_k, dy_k) d \xi.
\]
Using the fact that $a_k$ is symmetric, we argue as in \eqref{eq:sep29.2017.6} to express the upper bound on $B_k$ in terms of integrals involving  just the variables $(\xi_1, \xi_2,\xi_3).$ We obtain  
\begin{align}
|B_k| & \leq 4 \pi^2  |x_1-y_1|^2 \int_{\Mb^k} |a_k(\xi)| |\xi_1|^3 d\xi  \nonumber\\
&+   4 \pi^2 (k-1) W_2^2(m, \nu) \Bigl(  \int_{\Mb^k} |a_k(\xi)| |\xi_1|  |\xi_2|^2 d\xi  +(k-2) \int_{\Mb^k} |a_k(\xi)| |\xi_1|  |\xi_2||\xi_3| d\xi    \Bigr) \nonumber\\
&  8 \pi^2(k-1) |x_1-y_1|W_2(m, \nu) \int_{\Mb^k} |a_k(\xi)| |\xi_1|^2  |\xi_2| d\xi  \nonumber
\end{align} 
We use that $2|\xi_2| |\xi_3| \leq |\xi_2|^2+|\xi_3|^2$ and use again the fact that $a_k$ is symmetric  and argue as in \eqref{eq:sep29.2017.6} to eliminate the variables $\xi_3$ from the previous integral. We obtain 
\begin{align}
|B_k| & \leq 4 \pi^3  |x_1-y_1|^2 \Bigl( \int_{\Mb^k} |a_k(\xi)| |\xi_1|^3 d\xi +(k-1)\int_{\Mb^k} |a_k(\xi)| |\xi_1|^2  |\xi_2| d\xi \Bigr) \nonumber\\
&+   4 \pi^2 (k-1) kW_2^2(m, \nu)   \int_{\Mb^k} |a_k(\xi)| |\xi_1|^2  |\xi_2| d\xi   \nonumber
\end{align} 
We exploit the first order Taylor expansion of $\nabla_w F_{\Phi_k}$ around $(x_1, m)$, use \eqref{eq:sep30.2017.2} and \eqref{eq:sep30.2017.3} to conclude that by linearity that 
\begin{align}
\nabla_w G_N[\nu](y_1)-  \nabla_w G_N[m](x_1) &= \nabla \bigl( \nabla_w G_N[m](x_1) \bigr)(y_1-x_1) \nonumber\\
&+   \int_{\Mb^2} \nabla^2_w G_N[m](x_1, x_2) \bigr)(y_2-x_2) \gamma(dx_2, dy_2) +R_N\label{eq:sep30.2017.5}
\end{align} 
where the remainder $R_N$ satisfies
\begin{align}
|R_N| & \leq 4\pi^2 |x_1-y_1|^2 \sum_{k=1}^N \Bigl({1\over k!}  \int_{\Mb^k} |a_k(\xi)| |\xi_1|^3 d\xi +{k-1\over k!}\int_{\Mb^k} |a_k(\xi)| |\xi_1|^2  |\xi_2| d\xi \Bigr)  \nonumber\\
&+   4 \pi^2 \sum_{k=1}^N {(k-1) k \over k!} W_2^2(m, \nu)   \int_{\Mb^k} |a_k(\xi)| |\xi_1|^2  |\xi_2| d\xi   \nonumber
\end{align} 
We combine \eqref{eq:uniform-bound-ak}, \eqref{eq:sep01.2017.4}, \eqref{eq:sep01.2017.5}  and \eqref{eq:sep01.2017.5bis} to obtain that a universal constant $\bar C$ such that 
\[
|R_N| \leq \bar C_m:=  \sum_{k=1}^\infty {\bar C \over k^{1+\delta}}.
\]
If necessary, we replace $C_m$ by $\max\{C_m, \bar C_m\}$. We use \eqref{eq:sep30.2017.5} to obtain 
\begin{equation}\label{eq:sep30.2017.5new2}
\Bigl| \nabla_w G_N[\nu](y_1)-  \nabla_w G_N[m](x_1) -P^N_\gamma[m](x_1, y_1)\Bigr| \leq C_m\Bigl( |x_1-y_1|^2+W_2^2(m,\nu)\Bigr)
\end{equation}
Here, 
\[
P^N_\gamma[m](x_1, y_1):= \nabla_{x_1}  \bigl( \nabla_w G_N[m](x_1) \bigr) (y_1-x_1) +  \int_{\Mb^2} \nabla_{w}^2 G_N[m](x_1,x_2) (y_2-x_2) \gamma(dx_2, dy_2),
\]

We have  
\[
\Bigl| g_k^2[m](x_1, x_2) \Bigr| \leq {k-1 \over k!}  \int_{\Mb^k } |a_k(\xi)| \;  |\xi_1|^2 d\xi.
\]
Once again, thanks to  \eqref{eq:uniform-bound-ak} and \eqref{eq:sep01.2017.4}  we obtain that the series $\bigl(-4 \pi^2 \sum_{k=1}^N g_k^2[m](x_1, x_2)\bigr)_N$ is a Cauchy sequence for the uniform convergence and so,  it converges uniformly to the continuous function $\bar A$ at the right handside of \eqref{eq:sep30.2017.8}. We let $N$ tend to $\infty$ in \eqref{eq:sep30.2017.5new2} to conclude that 
\begin{equation}\label{eq:sep30.2017.5ter}
\Bigl| \nabla_w U_0[\nu](y_1)-  \nabla_w U_0[m](x_1) -P_\gamma[m](x_1, y_1)\Bigr| \leq C_m\Bigl( |x_1-y_1|^2+W_2^2(m,\nu)\Bigr)
\end{equation}
Here, 
\[
P_\gamma[m](x_1, y_1):=\tilde A[m](x_1) (y_1-x_1) +  \int_{\Mb^2} \bar A[m](x_1,x_2) (y_2-x_2) \gamma(dx_2, dy_2),
\]
We use Theorem \ref{th:grad-hess-series} to obtain that 
\[
\tilde A[m](x_1)=\nabla \bigl( \nabla_w U_0[m](x_1) \bigr)  \quad \bar A[m](x_1, x_2)=\nabla^2_w U_0[m](x_1, x_2) \bigr)
\]
and conclude the proof of (ii). 

(iii) By Corollary \ref{co:grad-hess1} and the above uniform convergences, we have 
\[
\triangle_{w, \epsilon} U_0(m) \equiv \sum_{k=1}^\infty \triangle_{w, \epsilon} \biggl( {1 \over k!} \int_{\Mb^k } a_k(\xi) F^k_\xi[m] d \xi\biggr)= 
\sum_{k=1}^\infty  {1 \over k!} \int_{\Mb^k } a_k(\xi) \triangle_{w, \epsilon} \bigl( F^k_\xi[m] \bigr) d \xi.
\]
This concludes the proof. \endproof

In the next proposition, we assume that we are given  $(b_k)_{k=1}^\infty \subset \mathcal A^0$ be such that \eqref{eq:uniform-bound-ak}, \eqref{eq:sep01.2017.4}, \eqref{eq:sep01.2017.5} , \eqref{eq:sep01.2017.5bis}. We assume that 
\begin{equation}\label{eq:oct06.2017.0} 
 \int_{\Mb^k} \bigl(|a_k(\xi)| +|b_k(\xi)| \bigr) \lambda_k^2(\xi) d\xi \leq {C k!\over k^{\delta}}
\end{equation}
and 
\begin{equation}\label{eq:oct06.2017.0new} 
\sum_{k=1}^\infty {1 \over k!} \int_{\Mb^k} |a_k(\xi)|^2 \lambda_k^4(\xi) d\xi <\infty.
\end{equation}
Assume $(b_k)_{k=1}^\infty \subset \mathcal A^0$ is such that the analogous of \eqref{eq:uniform-bound-ak}, \eqref{eq:sep01.2017.4}, \eqref{eq:sep01.2017.5} , \eqref{eq:sep01.2017.5bis}. Define 
\[
F:= \sum_{k=1}^\infty {1 \over k!}  \int_{\Mb^k} a_k(\xi) F_\xi^k[m] d\xi, \qquad G:= \sum_{k=1}^\infty {1 \over k!}  \int_{\Mb^k} b_k(\xi) F_\xi^k[m] d\xi.
\]

\begin{proposition}\label{pr:integration-by-parts}  Assume $s=0$, and $(a_k)_{k=1}^\infty$ and $(b_k)_{k=1}^\infty$ \eqref{eq:uniform-bound-ak}, \eqref{eq:sep01.2017.4}, \eqref{eq:sep01.2017.5} , \eqref{eq:sep01.2017.5bis},  and $(a_k)_{k=1}^\infty$ further satisfies  \eqref{eq:sep30.2017.8} hold. If \eqref{eq:oct06.2017.0} and \eqref{eq:oct06.2017.0new} hold then 
\begin{enumerate}
\item[(i)] The following functions belong to the $d$--cartesian product of $\mathcal H^0(\calPtwod)$ by itself: 
\[
m \rightarrow \int_{\Mb} \nabla_w F[m](x_1)m(dx_1),\quad \int_{\Mb} \nabla_w G[m](x_1)m(dx_1) 
\]
\item[(ii)] We have $\triangle_w F \in \mathcal H^0(\calPtwod).$
\item[(iii)]  We have the integration by parts formula:
\[
-\langle \triangle_w F; G \rangle_{H^0}= \left \langle \int_{\Mb} \nabla_w F[m](x)m(dx) ; \int_{\Mb} \nabla_w G[m](x)m(dx)  \right \rangle_{H^0} 
\]
\item[(iv)] In particular, 
\[
-\langle \triangle_w F; F \rangle_{H^0}= \left \|  \int_{\Mb} \nabla_w F[m](x)m(dx) \right \|^2_{H^0} 
\]
\end{enumerate}

\end{proposition} 
\proof{} Corollary \ref{co:superposition2} ensures that $F$ is twice continuously differentiable, $G$ is continuously differentiable. Setting 
\[
f_k(\xi):= \sum_{j=1}^k \xi_j  a_k, \quad g_k(\xi):= \sum_{j=1}^k \xi_j b_k,
\]
we use Lemma \ref{le:computationFk} and Corollary \ref{co:superposition2} to obtain the explicit expressions 
\begin{equation}\label{eq:oct06.2017.01} 
\int_{\Mb} \nabla_w F[m](x_1)m(dx_1)= \sum_{k=1}^\infty {-2\pi i \over k!}  \int_{\Mb^k} f_k(\xi) F_\xi^k[m] d\xi, 
\end{equation}
\begin{equation}\label{eq:oct06.2017.02}  
\int_{\Mb} \nabla_w G[m](x_1)m(dx_1)= \sum_{k=1}^\infty {-2\pi i \over k!}  \int_{\Mb^k} g_k(\xi) F_\xi^k[m] d\xi
\end{equation}
and 
\begin{equation}\label{eq:oct06.2017.03} 
\triangle_w F[m](x_1)=-\sum_{k=1}^\infty \int_{\Mb^k} a_k(\xi)\lambda^2_k(\xi) F_\xi^k[m] d\xi
\end{equation} 
Combining \eqref{eq:uniform-bound-ak} and \eqref{eq:oct06.2017.0} we have 
\begin{equation}\label{eq:oct06.2017.1} 
\int_{\Mb^k} |f_k(\xi)| d\xi\leq \int_{\Mb^k} |a_k(\xi)| d\xi + \int_{\Mb^k} |a_k(\xi)| \lambda_k^2(\xi) d\xi\leq {2 C k!\over k^{\delta}}
\end{equation}
Since  $(a_k)_{k=1}^\infty \subset \mathcal A^0$
\begin{equation}\label{eq:oct06.2017.1b} 
\sum_{k=1}^\infty \int_{\Mb^k} {|f_k(\xi)|^2\over k!} d\xi=  \int_{\Mb^k} \sum_{k=1}^\infty {|a_k(\xi)|^2\over k!} \lambda_k^2(\xi) d\xi<\infty.
\end{equation}
We combine \eqref{eq:oct06.2017.1} and \eqref{eq:oct06.2017.1b} to conclude the the proof of (i) for $F$. Similarly, we conclude the proof of (i) for $G.$ 

(ii) is obtained as a consequence of \eqref{eq:oct06.2017.0}  and \eqref{eq:oct06.2017.0new}.

(iii) Since $\triangle_w F \in \mathcal H^0(\calPtwod), G$ we use their expressions to obtain that have 
\begin{equation}\label{eq:oct06.2017.2} 
\langle \triangle_w F; G \rangle_{H^0}= -\sum_{k=1}^\infty {1 \over k!} \int_{\Mb^k} a_k(\xi) (b_k(\xi))^* \lambda_k^2(\xi)d\xi 
\end{equation}
We use the expressions in \eqref{eq:oct06.2017.01}  and \eqref{eq:oct06.2017.02} to conclude thet 
\[
\langle \nabla_w F; \nabla_w G \rangle_{H^0}= \sum_{k=1}^\infty {4\pi^2 \over k!} \int_{\Mb^k} a_k(\xi) (b_k(\xi))^* \biggl|  \sum_{j=1}^k \xi_j  \biggr|^2 d\xi 
\]
This, together with \eqref{eq:oct06.2017.2} concludes the proof of (iii). \endproof

\section{Recovery of $k$--polynomial of $\mathcal H^s(\calPtwod)$ functions} 
In this section, we study two type of problems. The first one consists to know if we can write any symmetric function $\Phi_k \in C(\Mb^k)$ in terms of $F_{\Phi_k}$. The second question consists in knowing when we can write $F_{\Phi_k}$ in terms of $\sum_{k=1}^\infty F_{\Phi_k}.$ The second problem can be formulated in terms of Fourier transforms. Given $(a_k)_{k=1}^\infty \in \mathcal A^s$ such that   $a_k  \in L^1(\Mb^k)$,  define 
\[
F_N: m \mapsto  \sum_{k=1}^N  {1\over k! } \int_{\Mb^k } a_k(\xi) F^k_\xi[m] d \xi.
\]
for  any natural number $N$. For $k\in \{1, \cdots, N\}$ we are able to recover $\int_{\Mb^k } a_k(\xi) F^k_\xi[m] d \xi$  from $F_N.$ For instance, the recovering allows to conclude that  if $F_N \equiv 0$ then $a_k \equiv 0$ for any $k \in \{1, \cdots, N\}.$ In this section, we endeavour  to prove a more general statement by allowing $N=\infty,$ at the expense of imposing additional growth conditions on the $\|a_k\|_{L^1}.$ Further assume  there exist $C, \delta>$ such that  \eqref{eq:uniform-bound-ak} holds and 
\[
\sum_{k=1}^\infty {1\over k! }\int_{\Mb^k } a_k(\xi) F^k_\xi[m] d \xi
\]
converges uniformly on $\calPtwod$ on $\calPtwod$ (cf. Lemma \ref{le:Hs}) to a function we denote as $F$ . Set   
\begin{equation}\label{eq:sep23.2017.3} 
\Phi_k:= \widehat a_k \in \Sk \cap L^2(\Mb^k)
\end{equation} 
Note, $\Phi_k $ is continuous and by Riemann--Lebesgue lemma (cf. e.g. \cite{SteinS}  Exercise 22 pp. 94)   
\begin{equation}\label{eq:riemann-lebesgue2}
\lim_{|x| \rightarrow \infty} \Phi_k(x)=0.
\end{equation}

\subsection{The inverse of the restriction $\Phi \rightarrow F_\Phi$ to polynomial} A natural question we address in the subsection is the reconstruction of $\Phi$ from $F=F_\Phi.$ For example, assume $k=2$ and $F=F_{2\Phi}$. We have 
\begin{equation}\label{eq:sep03.2017.2}
F(\delta_{a_1})=   \int_{\Mb^2} \Phi(x_1, x_2) \delta_{a_1}(dx_1) \delta_{a_1}(dx_2)= \Phi(a_1, a_1).
\end{equation}
Similarly,  
\begin{equation}\label{eq:sep03.2017.3}
F\Bigl({\delta_{a_1}+ \delta_{a_2} \over 2} \Bigr)= {1 \over 4} \Bigl(\Phi(a_1, a_1)+ \Phi(a_2, a_2) + 2 \Phi(a_1, a_2)\Bigr)
\end{equation} 
We combine \eqref{eq:sep03.2017.2} and \eqref{eq:sep03.2017.3} to obtain the polarization identity 
\begin{equation}\label{eq:sep03.2017.4}
\Phi(a_1, a_2)= 2 F \Bigl({\delta_{a_1}+ \delta_{a_2} \over 2} \Bigr) - {1 \over 2}F \bigl(\delta_{a_1} \bigr)-{1\over 2} F \bigl(\delta_{a_2} \bigr).
\end{equation}
Observe that if $F_{2\Phi} \equiv 0$, \eqref{eq:sep03.2017.4} implies $\Phi \equiv 0$ and so, $\Phi \rightarrow F_\Phi$ is an injective map of $C(\Mb^2/P_2)$.  In general, we could determine $\Phi$ applying the idea of coming from the construction of polar forms, either by a construction using differentiation or the inclusion--exclusion principle. To avoid differentiating, we chose here to use the inclusion--exclusion principle. 

Given two positive integers $1\leq r \leq k$ we defined the index set of multi--indexes 
\[
C_r^k := \{ (i_1,\cdots, i_r) \; | \;  i_1, \cdots, i_r \in \{1, \cdots ,k \} , i_1 < i_2 < \cdots < i_r \}
\]
Now, given $x = (x_1,x_2,\cdots ,x_k) \in \Mb^k$, and for a given multi--index $(i_1,\cdots .i_r)=I \in C_r^k$, we define $m_{x_I}$ as follows:
\[
m_{x_I} := \frac{1}{r} \sum_{j =1}^r \delta_{x_{i_j}}.
\]
Given a continuous function $F: \calPtwod \rightarrow \Reals$, we define  
\[
O_k(F) (x_1,\cdots ,x_k)  = \frac{1}{k!} \sum_{r=1}^k   \biggl( (-1)^{k-r} r^k  \sum_{I \in C_r^k} F (m_{x_I})  \biggr).
\]

\begin{theorem}\label{th:inverseOk}  The map $kO_k$ is the inverse map of \ \ $\Phi\rightarrow F_{\Phi}.$ In other words, we have 
\[
 \Phi(x_1,\cdots ,x_k)  = \frac{1}{k!}  \sum_{r=1}^k   \biggl( (-1)^{k-r} r^k  \sum_{I \in C^k_r} F_{k \Phi }(m_{x_I})  \biggr)=O_k(F_{k \Phi }) (x_1,\cdots ,x_k)
\] 
for any $x_1, \cdots, x_k \in \Mb.$ In other words, $O_k(F_\Phi)=\Phi/k.$ 
\end{theorem} 
\proof{} Note first that for any $F=F_{k\Phi}$ such that $\Phi \in C(\Mb^k/P_k)\bigr\}$, $O_k(F)$ is continuous and symmetric in the sense that it is defined on the quotient space $\Mb^k/P_k.$ Let $\mathcal{M}_c(\Mb)$ denote the set of signed Radon measures of compact support on $\Mb.$ This is a vector space which contains the set of Radon probability measures on $\Mb.$ We define $\alpha :  \mathcal{M}^k_c(\Mb) \rightarrow \mathbb{R}$ by 
\[ 
 \alpha(m_1,...,m_k) :=  \int_{\Mb^k } \Phi(x_1, x_2,... , x_k)m_1(dx_1)m_2(dx_2)... m_k(dx_k), 
\]
for $m_1, \cdots, m_k \in \mathcal{M}_c(\Mb)$. This is a $k$--multilinear form  and so 
\[ 
m \rightarrow \tilde{\alpha} (m) := \alpha(m,...,m)\,.
\]
is a $k$--homogeneous functional on $ \mathcal{M}_c(\Mb).$ We apply the polarization identity coming from the inclusion--exclusion principle, which goes back to \cite{MazurO} \cite{Nelson} \cite{Schetzen} (cf. \cite{Thomas} for a recent and simple proof, and \cite{Greenberg} \cite{Ward} for a formulation in terms of $n$-th defects of $F$). We obtain  
\begin{eqnarray}\label{polarpolar}
\alpha(m_1,...,m_k) ={1\over k! } \sum_{r=1}^k (-1)^{k-r} \sum_{I \in C^k_r}   \tilde{\alpha} \left( \sum_{i \in I } m_i \right) \,.
\end{eqnarray}
Setting $m_i = \delta_{x_i} $, using the definition of $\alpha$ and the fact that $\alpha$ is $k$--multilinear, we have 
\[
\Phi(x_1, x_2,... , x_k)=  \alpha(\delta_{x_1},\cdots, \delta_{x_k})= r^k \alpha\Bigl({\delta_{x_1} \over r},\cdots ,{\delta_{x_k}\over r}\Bigr).
\]
This, together with \eqref{polarpolar} yields 
\[
\Phi(x_1, x_2,... , x_k)={r^k \over k!} \sum_{r=1}^k (-1)^{k-r} \sum_{I \in C^k_r}   \tilde{\alpha} \left(  m_{x_I} \right).
\]
Since $F_{k \Phi}$ and $ \tilde{\alpha}$ coincide on the set of Radon probability measure, we conclude the proof of the theorem. \endproof

\subsection{One dimensional analytical extension} For any $ \lambda \in (0,1)$, $m \in  \calPtwod$ and $y \in \Mb^k$, we have 
\begin{eqnarray}
& &  F_{ \Phi_k} [\lambda m + (1-\lambda) \delta_y] \nonumber\\
&=& \frac{1}{k} \sum_{l=0}^k  \frac{k!}{ l! (k-l)! } \lambda^l (1-\lambda)^{k-l}  \int_{\Mb^l} \Phi_k(x_1,..,x_l, \underbrace{   y,...,y}_{k-l \text{ times}})
 m (dx_1) \cdots m (dx_l) \label{eq:sep23.2017.1}
\end{eqnarray}
Now since 
\begin{equation}\label{eq:sep23.2017.4}
\| \Phi_k\|_{L^\infty} \leq \|\widehat a_k\|_{L^1}<\infty
\end{equation} 
we may apply the dominated convergence theorem and use \eqref{eq:riemann-lebesgue2} to obtain that all the terms in \eqref{eq:sep23.2017.1}, except corresponding to $k=l,$ tend to $0$ as $|y|\rightarrow \infty.$ Thus, 
\begin{equation}\label{eq:sep23.2017.2}
 \lim_{|y| \rightarrow \infty} F_{\Phi_k} \left(\lambda m + (1-\lambda) \delta_y \right)=  \lambda^k F_{ \Phi_k} [m].
\end{equation}

\begin{remark}\label{re:F_lambda_m} Suppose $(a_k)_{k=1}^\infty \in \mathcal A^s$ satisfies \eqref{eq:uniform-bound-ak} so that \eqref{eq:sep23.2017.3} holds and $\Phi_k  \in C_0 (\Mb^k)$.   
\begin{enumerate}
\item[(i)] For any $(\lambda, m) \in (0,1) \times \calPtwod$ 
\begin{equation}
\label{def_F_lambda_m}
\lim_{|y| \rightarrow \infty} F \left(\lambda m + (1-\lambda) \delta_y \right) =  \sum_{k=1}^\infty {1 \over k! } \lambda^k F_{ \Phi_k} [m]=: \mathcal F[\lambda, m].
\end{equation}  
\item[(ii)] Hence, $\lambda \mapsto \mathcal F[\lambda, m]$ admits a extension denoted the same way, which is continuously differentiable at $0$. For any $l\geq 1$,  
\begin{equation}\label{diff_F_lambda_m}
{\partial^l  \mathcal F\over \partial \lambda^l} [\lambda, m]  \bigg|_{\lambda = 0}=   F_{ \Phi_l} [m] \,.
\end{equation}
\end{enumerate}
\end{remark}
\proof{} (i) By \eqref{eq:sep23.2017.1} and \eqref{eq:sep23.2017.4} we obtain for any integers $1\leq M<N$, 
\begin{eqnarray}
\sum_{k=M}^N {1\over k!} \Bigl|F_{ \Phi_k} [\lambda m + (1-\lambda) \delta_y]- \lambda^k F_{ \Phi_k} [m] \Bigr|
&\leq & \sum_{k=M}^N {1\over k! \; k} 
 \sum_{l=0}^{k-1}  \frac{k!}{ l! (k-l)! } \lambda^l (1-\lambda)^{k-l}  \bigl\|\widehat a_k\bigr\|_{L^1}    
 \nonumber\\
&\leq & \sum_{k=M}^N {\bigl\|\widehat a_k\bigr\|_{L^1} \over k! \; k}    \nonumber.
\end{eqnarray}
This, together with \eqref{eq:uniform-bound-ak} implies 
\[
\sum_{k=M}^\infty {1\over k!} \Bigl|F_{ \Phi_k} [\lambda m + (1-\lambda) \delta_y]- \lambda^k F_{ \Phi_k} [m] \Bigr| \leq \sum_{k=M}^\infty {C \over k^{\delta+1}}
\]
Thus by \eqref{eq:sep23.2017.2},
\[
\limsup_{|y|\rightarrow \infty}  \Bigl| F \left(\lambda m + (1-\lambda) \delta_y \right) -\mathcal F[\lambda, m] \Bigr| \leq \sum_{k=M}^\infty {C \over k^{\delta+1}}.
\]
We let $M$ tend to $\infty$ to obtain (i).

(ii) Observe the domain of convergence of the analytic function $(\sum_{k=1}^\infty z/k^{1+\delta})$ in the complex plane $\mathbb C$ is the unit disk. Since by \eqref{eq:uniform-bound-ak}   
\[
{1 \over k! } \sup_{\calPtwod} \big|F_{ \Phi_k}\big| \leq {C \over k^{1+\delta}  }
\] 
we conclude $\mathcal F[\cdot, m]$ extends to an analytic function on the unit disk. Therefore, it is differentiable at $0$ and one  checks that (ii) holds.  \endproof

\subsection{Projections of a subset of $ \bigoplus_{k}  \cSk$ onto $\cSkcomplex$}

\begin{definition}\label{proj_curlyHs}  For any natural number $k$, thanks to Remark \ref{re:F_lambda_m}, we may define the 
following operator $\pi_k: \mathcal{H}^s(\calPtwod) \rightarrow \mathcal{S}ym[k](\Reals)$ as follow: 
\[ 
 \pi_k (F) [m]  : = {\partial^l  \mathcal F\over \partial \lambda^l} [\lambda, m]  \bigg|_{\lambda = 0}  \quad \forall \; m \in \calPtwod.
\]
\end{definition}

Then the following Corollary is a direct consequence of Remark \ref{re:F_lambda_m}.

\begin{corollary}\label{F_lambda_m} Suppose $(a_k)_{k=1}^\infty \in \mathcal A^s$ satisfies \eqref{eq:uniform-bound-ak} so that \eqref{eq:sep23.2017.3} holds and $\Phi_k := \widehat a_k \in C_0 (\Mb^k)$.   Let $F \in \mathcal{H}^s(\calPtwod)$ be the continuous function obtained as the uniform limit of the series 
\[
\sum_{k=1}^\infty  {1\over k! } \int_{\Mb^k } a_k(\xi) F^k_\xi[m] d \xi= \sum_{k=1}^\infty  {1\over k! } F_{\Phi_k}[m] \quad \forall m \in \calPtwod.
\] 
\begin{enumerate}
\item[(i)] We have \begin{eqnarray*}
\pi_k (F) =  \frac{1}{k!}  F_{ \Phi_k}  \quad \text{and} \quad  F =  \sum_{k=1}^\infty  \pi_k (F).
\end{eqnarray*} 
\item[(ii)] In particular $\pi_k \circ \pi_k (F) =  \pi_k (F) .$
\item[(iii)] By (i) and Theorem \ref{th:inverseOk} we obtain   
\[
\Phi_k = k\; k! \;\; O_k(\pi_k (F)).
\]
\item[(iv)] Using (iii) and the Fourier transform inverse formula we have 
\[
a_k(\xi) =  k\; k! \;\; \widehat{O_k ( \pi_k (F))}(-\xi).
\]
\end{enumerate}
\end{corollary}
Because of (i) and (ii) in Corollary \ref{F_lambda_m}  we refer $\pi_k$ in Definition \ref{proj_curlyHs} as projection operator.

\begin{remark} (Sufficient conditions for uniqueness of Fourier coefficients)\label{re:unique-proj} Suppose $(a_k)_{k=1}^\infty \in \mathcal A^s$ satisfies \eqref{eq:uniform-bound-ak} and  let $F:\calPtwod \rightarrow \mathbb C$ be defined by 
\[
F[m]=\sum_{k=1}^\infty  {1\over k! } \int_{\Mb^k } a_k(\xi) F^k_\xi[m] d \xi \quad \forall m \in \calPtwod.
\] 
If $F\equiv 0$ then for any natural number $k,$ we have $a_k \equiv 0$ 

\end{remark}
\proof{} The remark is obtained as a direct consequence of Corollary \ref{F_lambda_m}. \endproof

\begin{definition}\label{de:norm-innner-product}  Let $s \geq 0$ and let $\lambda_k$ be the function defined in \eqref{eq:lambdak}.
\begin{enumerate} 
\item[(i)] Let $(a_k)_{k=1}^\infty, (b_k)_{k=1}^\infty \in \mathcal A^s$ be such that \eqref{eq:uniform-bound-ak}. Let $F, G:\calPtwod \rightarrow \mathbb C$ be defined by 
\[
F[m]=\sum_{k=1}^\infty  {1\over k! } \int_{\Mb^k } a_k(\xi) F^k_\xi[m] d \xi
\]
and 
\[ 
G[m]=\sum_{k=1}^\infty  {1\over k! } \int_{\Mb^k } b_k(\xi) F^k_\xi[m] d \xi \quad \forall m \in \calPtwod.
\] 
By Remark \ref{re:unique-proj}, $(a_k)_{k=1}^\infty, (b_k)_{k=1}^\infty$ are uniquely determined. We define 
\[
\langle F; G\rangle_{H^s}= \sum_{k=1}^\infty {1 \over k!} \int_{\Mb^k} a_k(\xi) b_k^*(\xi)  \Bigl( 1+ \lambda^2_k(\xi)\Bigr)^s d\xi
\]
and 
\[ 
  \| F \|^2_{H^s}=\langle F; F\rangle_{H^s}.
\]
\item[(ii)] Let $(V_k^n)_{k=1}^\infty, (W^n_k)_{k=1}^\infty \in \mathcal A^s$ for $ n=1, \cdots, d$ and let $V_k\in L^2(\Mb^k; \mathbb C^d)$ (resp. $W_k \in L^2(\Mb^k; \mathbb C^d)$) be the vector field whose components  $(V_k^n)_{n=1}^d,$ (resp. $(W^n_k)_{n=1}^d$ satisfy \eqref{eq:uniform-bound-ak}. Let 
\[
V[m]=\sum_{k=1}^\infty  {1\over k! } \int_{\Mb^k } V_k(\xi) F^k_\xi[m] d \xi
\] 
\[
 W[m]=\sum_{k=1}^\infty  {1\over k! } \int_{\Mb^k } W_k(\xi) F^k_\xi[m] d \xi.
\]
We define 
\[
\langle V; W\rangle_{H^s}= \sum_{n=1}^d \langle V^n; W^n\rangle_{H^s}, \quad \| V \|^2_{H^s}=\langle V; V\rangle_{H^s}.\] 

\end{enumerate}   
\end{definition}

\begin{definition}(Conditional continuity of  $O_k \circ \pi_k$)\label{de:continuity1}  Let $s \geq 0$ and let $N \geq 1$ be an integer. Suppose there exists $(a_k)_k \in \mathcal A^s$ such that  $a_k \equiv 0$ for any $k>N$,    
\begin{equation}\label{eq:oct01.2017.0}
a_k \in  L^1(\Mb^k/P_k; \mathbb C) \quad \text{and} \quad \Phi_k:=\widehat a_k \in \Sk \quad \forall k \in \{1, \cdots, N\}.
\end{equation}
We say that $F \in \mathcal{H}^s_N(\calPtwod)$ provided that 
\begin{equation}\label{eq:oct01.2017.1}
F=\sum_{k=1}^N {1\over k!} F_{\Phi_k}.
\end{equation}
\end{definition}

\begin{remark}\label{re:continuity0}  Let $N \geq 1$ be an integer. Suppose $(a_k)_k \in \mathcal A^s$ is such that \ref{eq:oct01.2017.0} holds. If there exists an integer $N$ such that $a_k \equiv 0$ for any $k>N$  then there are constants $C, \delta>0$ such that \eqref{eq:uniform-bound-ak}. In other words, $\mathcal{H}^s_N(\calPtwod) \subset \mathcal{H}^s(\calPtwod).$
\end{remark} 

\begin{lemma}\label{le:continuity}  Let $s\geq 0$ be a real number and $N$ be a natural number. There exists a constant $C_N$ such that for any natural number $k\leq N$ and any $F \in \mathcal{H}^s_N(\calPtwod)$  we have 
\[
  \|O_k \circ \pi_k(F)\|_{C(\Mb^k)} \leq C_N \sup_{m} |F[m]|. 
\]
In other words, if we endow $\mathcal{H}^s_N(\calPtwod)$ with the supremum norm then $O_k \circ \pi_k: \mathcal{H}^s_N(\calPtwod) \rightarrow C(\Mb^k/P_k)$ is continuous.
\end{lemma} 
\proof{} Let $(a_k)_k \in \mathcal A^s$ be such that \eqref{eq:oct01.2017.0} holds and $a_k \equiv 0$ for any $k>N$. Suppose $F$ satisfies \eqref{eq:oct01.2017.1} where $\Phi_k:=\widehat a_k$. Note 
\[
k O_k \circ \pi_k(F)= \Phi_k
\]
and so, we need to estimate the norms of the $\Phi_k$ in terms of the norm of $F.$ Set 
\[
\Psi_N(x_1, \cdots, x_N):= {N\over k\;  k!} {k! (N-k)! \over N!} \sum_{I \in C_k^N} \Phi_k(x_I).
\]
Observe that $\Psi_N \in C(\Mb^N/P_N)$ and $F=F_{\Psi_N}.$ By Theorem \ref{th:inverseOk} $\Psi_N=N O_N(F)$ and so, 
\begin{equation}\label{eq:oct01.2017.2}
|\Psi_N| \leq {N^{N+2} \over N!} \sup_{m} |F[m]|.
\end{equation}
Recall that as $a_k \in L^1$, \eqref{eq:riemann-lebesgue2} holds: 
\[
\lim_{ |(x_1, \cdots, x_k)| \rightarrow \infty} \Phi(x_1, \cdots, x_k)=0
\]
and so, $\lim_{|x_I| \rightarrow \infty}   \Psi_N(x_1, \cdots, x_N)$ exists. We conclude that 
\[
\lim_{ x_1  \rightarrow \infty} \Psi_N(x_1, \cdots, x_N)={\Phi_1(x_1) \over N}.
\] 
This, together with \eqref{eq:oct01.2017.2} implies 
\begin{equation}\label{eq:oct01.2017.3}
{1 \over N!} \|\Phi_1\|_{C(\Mb)} \leq  {N^{N+2} \over N!} \sup_{m} |F[m]|.
\end{equation} 
Observe 
\[
\lim_{ x_1, x_2  \rightarrow \infty}\Psi_N(x_1, \cdots, x_N)-{\sum_{j=1}^2\Phi_1(x_j) \over N!}= {N \over 2 2!} {2! (N-2)! \over N!}{\Phi_2(x_1, x_2)}.
\]
Hence, 
\[
{N \over 2\; 2!} {2! (N-2)! \over N!}  \|\Phi_2\|_{C(\Mb^2)} \leq \| \Psi_N  \|_{C(\Mb^N)}+ {2 \over N!} \| \Phi_1  \|_{C(\Mb)}
\] 
We combine \eqref{eq:oct01.2017.2} and \eqref{eq:oct01.2017.3} to conclude that 
\[
{N \over 2\; 2!} {2! (N-2)! \over N!}  \|\Phi_2\|_{C(\Mb^2)} \leq  {N^{N+2} +2N^{N+2}\over N!} \sup_{m} |F[m]|.
\] 
We repeat the same procedure $(N-2)$ times to conclude the proof of the Lemma.  \endproof

\section{Brownian motion on the Wasserstein space.} 
\noindent Throughout this section, $(W_t)_t$ is the standard $d$--dimensional brownian motion starting at the origin. Given $m \in \calPtwod$ and $\beta \in \Reals$ we set 
\[
\mathbb B_t^{m, \beta}:= (\id + \sqrt{2 \beta} W_t)_\# m, \qquad \sigma^\epsilon_t:= \sigma^{\epsilon, \beta}_t[m]:=(\id + \sqrt{2 \beta} W_t)_\# (G^\epsilon_t \ast m).
\] 
Here, $G^\epsilon_t$ is the heat kernel for the heat equation (cf. \eqref{eq:aug.19.2017.ter}) so that 
\[
\partial_t \bigl(G^\epsilon_t \ast u_0\bigr)=\epsilon \triangle \bigl(G^\epsilon_t \ast u_0\bigr), \qquad \text{on} \quad (0,\infty) \times \Mb.
\] 
Note that $\sigma_t:=\mathbb B_t^{m,  \beta}$ solves the system of stochastic differential equations  
\begin{equation}\label{eq:aug.19.2017.2}
\left\{
\begin{array}{ll}
d\sigma   &= {\rm div}(\sigma dS) \qquad \qquad \qquad \qquad \;\; \text{on} \quad [0,T] \times \calPtwod \qquad \sigma_0=m\\ 
dS  & =\beta  \nabla \ln(\sigma)  dt-\sqrt{2  \beta}dW\qquad \quad \text{on} \quad [0,T] \times \calPtwod,
\end{array}
\right.
\end{equation}  
The paths $\sigma^\epsilon$ satisfy the system of stochastic differential equations 
\begin{equation}\label{eq:general-heat2} 
\left\{
\begin{array}{ll}
d\sigma^\epsilon  &= {\rm div}(\sigma^\epsilon dS) \qquad \qquad \qquad \qquad \;\; \text{on} \quad [0,T] \times \calPtwod \qquad \sigma_0^\epsilon=m\\ 
dS  & =(\epsilon+ \beta)  \nabla \ln(\sigma^\epsilon)  dt-\sqrt{2  \beta}dW \;\; \text{on} \quad [0,T] \times \calPtwod.
\end{array}
\right.
\end{equation} 
These reduce to a single equation 
\begin{equation}\label{eq:one-eq}
d\sigma^\epsilon= {\rm div} \Bigl( (\epsilon+\beta)  \nabla \sigma^\epsilon  dt-\sqrt{2  \beta}\sigma^\epsilon dW \Bigr) \;\; \text{on} \quad (0,T) \times \calPtwod, \quad \sigma_0^\epsilon=m
\end{equation}

By definition, we say that $\sigma^\epsilon$ satisfies \eqref{eq:one-eq} if for every $\phi \in C^1\bigl((0,T); C^2_c(\Mb)\bigr)$ and every $0<s<r<T$ we have 
\begin{align}
&\int_{\Mb} \phi(r,y) \sigma_r^\epsilon(dy )- \int_{\Mb} \phi(s,x) \sigma_s^\epsilon(dx)  \nonumber\\ 
=& \int_s^r \Bigl(\int_{\Mb} \partial_t \phi(t, x)\sigma_t^\epsilon(dx) dt +\sqrt{2\beta} \langle \nabla \phi(t,x), \sigma_t^\epsilon(dx) d W \rangle  + (\epsilon+\beta) \triangle \phi(t, x)\sigma_t^\epsilon(dx)dt \Bigr)\label{eq:aug.19.2017.2def}\\
0=& \lim_{t\rightarrow 0^+} \mathbb E\bigl(W_2(\sigma_t^\epsilon, m) \bigr).\nonumber
\end{align} 

This section  lays down arguments favoring the fact $\mathbb B_t^{m}:=\mathbb B_t^{m,1}$   can be called the  Brownian  starting at $m.$ This is a random path in $\calPtwod$, starting at $m$ and corresponds to a common noise for finitely or infinitely many particles in $\Mb,$ dependent on whether or not the support of $m$ is finite.

\subsection{It\^o formula} Let $T>0$ and let $V:[0,T] \times \calPtwod \rightarrow \Reals$ be a continuous map such that $V$ is differentiable on $(0,T) \times \calPtwod.$ We further assume that for each $t>0,$ $V(t, \cdot)$ is twice continuously differentiable and $\partial_t V$, $\nabla_w V$, $\nabla (\nabla_w V)$ and $\nabla_w^2 V$ are continuous on $[0,T] \times \calPtwod$ in the sense that they have a continuous extension. Suppose that for any $m \in \calPtwod$ there exists $C_m$ and a neighborhood $\mathcal O_m$ of $m$ such that 
\begin{equation}\label{eq:aug20.2017.5}
 |\nabla_w V(t, m)(x)| \leq C_m(1+|x|) \quad \forall t \in [0,T], \quad \forall x \in \Mb, \quad \forall \nu \in \mathcal O_m
\end{equation}
and 
\begin{equation}\label{eq:aug20.2017.5b}
\bigl| \nabla (\nabla_w V[\nu])\bigr|, \; \bigl| \nabla_w^2 V \bigr| \leq C_m \quad \forall  \nu \in \mathcal O_m.
\end{equation}
Given $m \in \calPtwod$, let $t \mapsto m_t \in \calPtwod$ be the unique solution to the equation 
\begin{equation}\label{eq:heateq}
\partial_t m_t =\epsilon \triangle m_t \quad \text{in} \quad (0,T) \times \Mb, \qquad m_0=m.
\end{equation}
Setting 
$$M(a, m):=  (\id +\sqrt{2\beta} a)_\# m,$$ 
we define the maps 
\begin{align}
(t, a) \mapsto A(t,a) &:=V\bigl(t, M(a, m_t) \bigr)  \nonumber\\
(t, a) \mapsto \Lambda_1(t,a) &:=V\bigl(t, M(a, m) \bigr)  \nonumber\\ 
(t, a) \mapsto \Lambda_2(t,a) &:=\partial_t V\bigl(t, M(a, m)\bigr)  \nonumber\\ 
(t, a) \mapsto \Lambda_3(t,a) & := \int_\Mb \nabla_w V\bigl(t, \bigl[M(a, m) \bigr](x)\bigr) M(a, m)(dx) \nonumber\\ 
(t, a) \mapsto \Lambda_4(t,a) & :=  \triangle_w V\bigl(t, \bigl[M(a, m)\bigr]\bigr) \nonumber \\
(t, a) \mapsto \Lambda_5(t,a) & :=  \partial_t V\bigl(t, M(a, m_t) \bigr)+  \int_\Mb \nabla_w V\bigl(t, \bigl[M(a, m_t) \bigr](x)\bigr) M(a, m_t)(dx)\nonumber
\end{align} 
 %
 %
\begin{lemma}\label{le:measurable} The maps $\Lambda_1, \cdots, \Lambda_5$ are continuous. 
\end{lemma}
\proof{} The map $M(\cdot, m): \Mb \rightarrow  \calPtwod$  is $\sqrt{2\beta}$--Lipschitz and the map $m \mapsto M(a, \cdot)$ is $1$--Lipschitz. Let $G_t^\epsilon$ be the Green function for equation \eqref{eq:heateq}, i.e. the heat kernal \eqref{eq:aug.19.2017.ter}. Since $G_{t+h}^\epsilon = G_{h}^\epsilon \ast G_{t}^\epsilon$, we have $m_{t+h}=G^\epsilon_h \ast m_t$ and so, by Lemma 5.17 \cite{GangboKP}, $t \mapsto m_t$ is a Lipschitz map. As a consequence, $t\mapsto M(a, m_t)\in  \calPtwod$ is a Lipschitz map. In particular, since all the previous maps are continuous, we obtain that  $A$, $\Lambda_1$ and $\Lambda_2$ are continuous as a composition of continuous maps.  Similarly, the maps  
\[
(t, x,a) \mapsto  \nabla_w V\bigl(t, M(a, m) \bigr)(x), \quad {\rm div} (\nabla_w V) \bigl(t, M(a, m) \bigr)(x)\bigr) 
\]
\[
(t, x,a) \mapsto \nabla_w V\bigl(t, M(a, m_t) \bigr)(x), \quad {\rm div} (\nabla_w V) \bigl(t, M(a, m_t) \bigr)(x)\bigr) 
\]
as well as
\[ 
 (t, x, y, a, m) \mapsto   \nabla^2_w V\bigl(t, M(a, m)\bigr)(x,y), \; \nabla^2_w V\bigl(t, M(a, m_t)\bigr)(x,y).
\]
are continuous as compositions of continuous maps .

In order to prove that $\Lambda_3, \Lambda_4, \Lambda_5$ are continuous, we consider an arbitrary sequence. We are to show that for any $\epsilon>0$ 
\begin{equation}\label{eq:aug20.2017.5.5}
\limsup_{n\rightarrow \infty} |\Lambda_i(t_n,a_n)-\Lambda_i(t,a)|\leq \epsilon \quad i \in \{3, 4, 5\}.
\end{equation}
We shall only present the proof when $i=3$  since the cases $i=4, 5$ are easier. Indeed, since $\nabla\bigl(\nabla_w V \bigr)$ and  $\nabla^2_w V$ are uniformly bounded, following the same lines of arguments we would obtain a proof in the cases $i=4, 5$.

Since $(M(a_n, m))_n$ converges to $M(a, m)$ in $\calPtwod$, we assume without loss of generality that $M(a_n,m ) \in \mathcal O_{M(a,m)}.$ Note the second moment of  $M(a_n, m)$ converges to that of $M(a, m).$ Thus we may choose $R$ large enough such that for any positive integer $n$,
\begin{equation}\label{eq:aug20.2017.6}
C_{M(a,m)}\Bigl( \sqrt{\int_{\{|x|\geq R \}} (1+|x|^2) M(a_n, m)(dx)} +\sqrt{\int_{\{|x|\geq R \}} (1+|x|^2) M(a,m)(dx)}  \Bigr)<\epsilon. 
\end{equation}
Let $g_R \in C_c(\Mb)$ be such that $0 \leq g_R \leq 1$ on $\Mb$, $g_R \equiv 1$ on $\{|x|\leq R \}$ and $g_R \equiv 0$ on $\{|x|\geq R+1 \}.$ We have 
\begin{align}
& \Lambda_3(t_n,a_n)-\Lambda_3(t,a) = \Lambda(R, n)  \nonumber\\
 & + \int_\Mb \nabla_w V\bigl(t, \bigl[M(a_n, m) \bigr](x)\bigr) g_R(x)M(a_n, m)(dx)- \int_\Mb \nabla_w V\bigl(t, \bigl[M(a, m) \bigr](x)\bigr) g_R(x)M(a, m)(dx),  \nonumber
\end{align} 
where we have set $ \Lambda(R, n)$ to be the  expression  
\[
\int_\Mb \nabla_w V\bigl(t, \bigl[M(a_n, m) \bigr](x)\bigr) (1-g_R)(x)M(a_n,m )(dx)- \int_\Mb \nabla_w V\bigl(t, \bigl[M(a, m) \bigr](x)\bigr) (1-g_R)(x)M(a, m)(dx).
\]
We combine \eqref{eq:aug20.2017.5} and \eqref{eq:aug20.2017.6} and use H\"older's inequality to obtain 
\begin{equation}\label{eq:aug20.2017.7}
|\Lambda(R, n)| \leq C_{M(a, m)}\Bigl( \sqrt{\int_{\{|x|\geq R \}} (1+|x|^2) M(a_n, m)(dx)} +\sqrt{\int_{\{|x|\geq R \}} (1+|x|^2) M(a, m)(dx)}  \Bigr)< \epsilon.
\end{equation}
By \eqref{eq:aug20.2017.5} and \eqref{eq:aug20.2017.5b} $\Bigl( \nabla_w V\bigl(t, \bigl[M(a_n, m) \bigr]\bigr) \Bigr)_n$ is equicontinuous on $\{|x|\leq R \}$ and so, it converges uniformly to $\nabla_w V\bigl(t, \bigl[M(a, m) \bigr]\bigr)$ there. Since convergence in $\calPtwod$ implies narrow convergence, we conclude that 
\[
\lim_{n\rightarrow \infty} \Bigl| \int_\Mb \nabla_w V\bigl(t, \bigl[M(a_n, m) \bigr](x)\bigr) g_R(x)M(a_n, m)(dx)- \int_\Mb \nabla_w V\bigl(t, \bigl[M(a, m) \bigr](x)\bigr) g_R(x)M(a, m)(dx)\Bigr|=0.
\] 
This, together with \eqref{eq:aug20.2017.7} yields \eqref{eq:aug20.2017.5.5} for $i=3.$ \endproof

%
 %
\begin{lemma}\label{le:measurable2} We have  
\[
\nabla \Lambda_1(t, a)= \sqrt{2\beta} \Lambda_3(t, M(a, m)), \qquad \triangle \Lambda_1(t, a)=2\beta \triangle_w V\bigl(t, \bigl[M(a, m)\bigr]\bigr)
\]  
and  
\[
\partial_t A(t,a)= \partial_t V(t, M(a, m_t)) +\epsilon \int_{\Mb} {\rm div} \Bigl( \nabla_w V(t, M(a, m_t))  \Bigr)  M(a, m_t)(dx).
\]
Hence, using the continuity properties obtained on these functions in Lemma \ref{le:measurable}, we deduce that $\Lambda_1(t, \cdot)$ and $A(t, \cdot)$  are twice continuously differentiable.
\end{lemma}
\proof{} (i) Let $a_0, a_1 \in \Mb$ and set 
\[a_\tau:=(1-\tau) a_0+\tau a_1, \quad {\bf v}(\tau,y)\equiv \sqrt{2\beta}(a_1-a_0), \quad T_\tau(x)= x+\sqrt{2\beta}(a_1-a_0) t.\] 
Since $T_\tau$ is the gradient of a convex function and pushes $M(a_0)$ forward to $M(a_\tau)$ we obtain that it is a geodesic in $\calPtwod.$ The vector field ${\bf v}$ is the velocity of minimal norm for the path $\tau \mapsto M(a_\tau).$ We have 
\[
{d \over d \tau} V(t, M(a_\tau))= \langle \nabla_w V(t,M(a_\tau));  {\bf v}(\tau,\cdot )\rangle_{M(a_\tau)}=\sqrt{2\beta} \langle \Lambda_3(t, M(a_\tau)), a_1-a_0 \rangle.
\]
This proves the first identity. 

Since $\partial_\tau {\bf v}+ \nabla {\bf v} {\bf v}=0$, we use Proposition \ref{pro:hessian1} to conclude that 
\[
{d^2 \over d \tau^2} V(t, M(a_\tau)) = {\rm Hess} V(t, M(a_\tau) )({\mathbf v}, {\mathbf v})=2\beta {\rm Hess} V(t, M(a_\tau) )(a_1-a_0, a_1-a_0).
\] 
Thus, 
\begin{align}
\Bigl(\nabla_a^2 V\bigl(t, M \bigr)\Bigr)(a_\tau) & =2\beta \int_\Mb \nabla\bigl(\nabla_w V(t, M(a_\tau) \bigr) M(a_\tau)(dx)  \nonumber\\
 & + 2\beta \int_{\Mb^2} \nabla_{m m} V\bigl(t, M(a_\tau)\bigr) M(a_\tau)(dx)M(a_\tau)(dy).  \nonumber
\end{align} 
Hence, 
\begin{equation}\label{eq:aug20.2017.9}
\Bigl(\triangle_a V\bigl(t, M \bigr)\Bigr)(a_\tau)= 2\beta \triangle_w V\bigl(t, \bigl[M(a_\tau)\bigr]\bigr),
\end{equation}
which proves the second identity.

For $t>0$, $m_t$ is absolutely continuous with respect to Lebesgue measure. In the reminder of the proof, we don't distinguish between $m_t$ and its Radon--Nikodym derivative with respect to Lebesgue measure. For the sake of not making the notation more cumbersome we abuse notation by writing of $m_t(dx)=m_t(x)dx.$ The velocity of $t \mapsto M(a, m_t)$ is 
\[
{\bf v}_t(y)=-\epsilon {\nabla m_t(y -\sqrt{2\beta} a) \over m_t(y -\sqrt{2\beta} a)}
\]
and so, it is uniformly bounded on $[\delta, T] \times \Mb$ for any $\delta>0.$ It is convenient to temporarily use the notation $Tx:=x -\sqrt{2\beta} a.$ The chain rule yields 
\begin{align}
\partial_t A(t,a)&= \partial_t V(t, M(a, m_t)) + \langle \nabla_w V(t, M(a, m_t)); {\bf v}_t\rangle_{M(a, m_t)} \nonumber\\
&=  \partial_t V(t, M(a, m_t)) -\epsilon \Bigl\langle \nabla_w V(t, M(a, m_t)); {\nabla m_t \circ T^{-1} \over m_t\circ T^{-1}}  \Bigr\rangle_{M(a, m_t)}    \nonumber\\
&=\partial_t V(t, M(a, m_t)) -\epsilon \Bigl\langle \nabla_w V(t, M(a, m_t)) \circ T; \nabla m_t    \Bigr\rangle_{L^2}    \nonumber\\ 
&= \partial_t V(t, M(a, m_t)) +\epsilon \int_{\Mb} {\rm div} \Bigl( \nabla_w V(t, M(a, m_t)) \circ T(y)   \Bigr) m_t(dy)\nonumber\\
&= \partial_t V(t, M(a, m_t)) +\epsilon \int_{\Mb} {\rm div} \Bigl( \nabla_w V(t, M(a, m_t)) \Bigr)  M(a, m_t)(dx) \label{eq:Itoa}
\end{align}
This is the third desired identity. \endproof

\begin{theorem}\label{th:measurable2} Setting $\sigma_t:= (\id +\sqrt{2\beta}W_t)_\# m$, and $\sigma_t^\epsilon:= (\id +\sqrt{2\beta}W_t)_\# (G^\epsilon_t \ast m).$ Then for any $0<s <r<T$, we have  
\begin{enumerate} 
\item[(i)] 
\[
V(r, \sigma_r)- V(s, \sigma_s)= \int_s^r \Bigl(\partial_t V(t, \sigma_t) dt +\sqrt{2\beta} \langle \nabla_w V(t, \sigma_t); dW_t\rangle_{\sigma_t}+ \beta \triangle_w V(t, \sigma_t)dt \Bigr).
\]  
\item[(ii)] 
\[
V(r, \sigma_r^\epsilon)- V(s, \sigma_s^\epsilon)= \int_s^r \Bigl(\partial_t V(t, \sigma_t^\epsilon) dt +\sqrt{2\beta} \langle \nabla_w V(t, \sigma_t^\epsilon); dW_t\rangle_{\sigma_t^\epsilon}+ \beta \triangle_{w, {\epsilon \over \beta}} V(t, \sigma_t^\epsilon)dt \Bigr).
\] \end{enumerate}  
\end{theorem}
\proof{} Applying It\^o's formula to $\Lambda_1(t, W_t)$ and using the differentiability properties of $\Lambda_1$ obtained in Lemma \ref{le:measurable2}  we have  (i). Similarly, we apply It\^o's formula to $A(t, W_t)$ and use the differentiability properties of $A$ obtained in Lemma \ref{le:measurable2} to obtain  (ii).  \endproof

Theorem \ref{th:measurable2} is in fact an extension of It\^o's formula to the set pf probability measures. An extension based on probabilistic arguments (hence different from ours) was proposed in books \cite{CarmornaDeI}  \cite{CarmornaDeII}. For very general processes, \cite{CardaliaguetDLL} gave a proof which uses a heavy machinery. We have opted to offer the above proof since our setting is different from that of the above cited prior works.

\subsection{Heat equation}\label{subsec:heat}  Assume  $U_0:\calPtwod \rightarrow \Reals$ is twice continuously differentiable (cf. Definition \ref{defn:Laplacian}). We assume for any $m \in \calPtwod,$ $\nabla \bigl( \nabla_w U_0[m]\bigr): \Mb  \rightarrow \mathbb R^{d \times d}$ and $\nabla^2_w U_0[m]: \Mb \times \Mb \rightarrow \mathbb R^{d \times d}$ are uniformly bounded and have a concave modulus of continuity independent of $m.$ We fix $\epsilon, \beta>0$ and define 
\[
U(t, m):= \mathbb E \Bigl(U_0\bigl(\sigma^{\epsilon, \beta}_t[m]\bigr) \Bigr) \qquad \forall (t, m) \in (0,\infty) \times \calPtwod.
\]
The arguments starting three lines after (81) in \cite{CardaliaguetDLL} lead to the following statements:  
\begin{enumerate} 
\item[(i)] $U$ is continuously differentiable on $(0,\infty) \times \calPtwod$ and for any $t>0$, $U(t, \cdot)$ is twice continuously differentiable on $\calPtwod.$ 
\item[(ii)] $U$ satisfies the heat equation 
\[
\partial_t U = \beta \triangle_{w, \frac{\epsilon}{\beta}} U \quad \text{on} \quad (0,\infty) \times \calPtwod, \quad U(0, \cdot)=U_0.
\] 
\end{enumerate}  
The reader should not be misled by the fact that in the statement of Theorem 4.3 \cite{CardaliaguetDLL}, the authors refer to a smooth Hamiltonian $H: \Mb^2 \rightarrow \Reals$ such that $\nabla_{pp} H(x,p)>0$, whereas $H\equiv 0$ in the current manuscript (see statements starting three lines after (81) \cite{CardaliaguetDLL}). Statements similar to the above ones have also been made in \cite{PLLions2007-2008} \cite{PLLions2008-2009}. 

When $U_0 \in \mathcal{H}^s(\calPtwod)$  (see Definition \ref{de:defncurlyHs}), in Theorem \ref{th:superposition} (see also Remark \ref{re:superposition}), we reach  conclusions much stronger than those in \cite{CardaliaguetDLL} \cite{PLLions2007-2008} \cite{PLLions2008-2009}. 

\subsection{Actions of Brownian motions on eigenfunctions}\label{subsection:effet}  Recall that $G^\epsilon_t$ denote the heat kernel for the heat equation, given in \eqref{eq:aug.19.2017.ter}. If $\beta>0$ and $\epsilon \geq 0$ then 
\[
\sigma^{\epsilon, \beta}_t[m]:=(\id+ \sqrt{2\beta} W_t)_\#(G^\epsilon_t \ast m) \qquad \forall m \in \calPtwod.
\]
In Lemma \ref{le:computationFk} we obtained a family $\{F^k_\xi, \; | \; \xi \in \Mb^k\}_{k=1}^\infty$ of eigenfunctions of  $\triangle_{w, \epsilon}$ with eigenvalues 
\begin{equation}\label{eq:aug26.2017.6}
-\lambda^2_{k, \epsilon}(\xi):= -4\pi^2 \Bigl(  \Big|  \sum_{j=1}^k \xi_j\Big|^2 + \epsilon \sum_{j=1}^n |\xi_j|^2\Bigr).
\end{equation} 
What is obvious is that 
\begin{equation}\label{eq:aug26.2017.5}
V_\xi^k(t, m):= \exp\Bigl( - \beta \lambda^2_{k, \frac{\epsilon}{\beta}}(\xi) t\Bigr) \; F^k_\xi.
\end{equation}
satisfies the heat equation 
\[
\partial_t V_\xi^k= \beta \triangle_{w, \frac{\epsilon}{\beta}} V_\xi^k, \qquad V_\xi^k(0,\cdot)= F^k_\xi.
\]
Under appropriate convergence conditions, its superposition 
\begin{equation}\label{eq:aug26.2017.8}
V(t, m):= \sum_{k=1}^\infty {1 \over k!} \int_{\Mb^k } a_k(\xi) \exp\bigl( - \beta \lambda^2_{k, \frac{\epsilon}{\beta}}(\xi) t  \bigr)  F^k_\xi[m] d \xi
\end{equation}
is also expected to satisfy 
\begin{equation}\label{eq:aug26.2017.9}
\partial_t V = \beta \triangle_{w, \frac{\epsilon}{\beta}} V , \qquad V (0,\cdot)= \sum_{k=1}^\infty {1 \over k!} \int_{\Mb^k } a_k(\xi) F^k_\xi[m] d \xi:=U_0.
\end{equation} 
According to  Subsection \ref{subsec:heat}, the function  
\[
U(t, m):=  \mathbb E\Bigl(U_0\bigl[\sigma^{\epsilon, \beta}_t[m]\bigr] \Bigr)
\]
also satisfies the same heat equation. This manuscript does not state a general uniqueness result for the heat equation and so, it is legitimate to ask if $U=V.$

\begin{lemma}\label{le:heat-special-case} We have  $V_\xi^k=U_\xi^k$ where, $U_\xi^k(t, m):=  \mathbb E\Bigl(F^k_\xi\bigl[\sigma^{\epsilon, \beta}_t[m]\bigr] \Bigr).$ 
\end{lemma}
 \proof{} To alleviate the notation, we write $m^{\epsilon}_t=G^\epsilon_t \ast m$ and don't display $m$ and $\beta$ in $\sigma^{\epsilon, \beta}_t[m]$. We have 
\begin{align}
k\; F^k_\xi[\sigma^\epsilon_t]  
&=\int_{\Mb^k} \exp\Bigl(-2\pi i \sum_{j=1}^k \langle \xi_j, x_j \rangle \Bigr) \sigma^\epsilon_t(dx_1) \cdots \sigma^\epsilon_t(dx_k) 
\nonumber\\
&= \exp\Bigl(-2\pi i \sum_{j=1}^k \langle \xi_j, \sqrt{2\beta} W_t \rangle \Bigr) 
\int_{\Mb^k} \exp\Bigl(-2\pi i \sum_{j=1}^k \langle \xi_j, x_j  \rangle \Bigr) m^\epsilon_t (dx_1) \cdots m^\epsilon_t (dx_k) \nonumber\\
&=\exp\Bigl(-2\pi i \sum_{j=1}^k \langle \xi_j, \sqrt{2\beta} W_t \rangle \Bigr) 
\prod_{j=1}^k \int_{\Mb} G^\epsilon_t \ast \exp\Bigl(-2\pi i \langle \xi_j, x_j  \rangle \Bigr) m(dx_j).\label{eq:aug26.2017.1}
\end{align}
One check that 
\[
u(t,x_j):=\exp\bigl(-2\pi i\langle \xi_j, x_j \rangle\bigr)  \exp\bigl(-4\pi^2 \epsilon t |\xi_j|^2 \bigr) 
\]
satisfies the equation 
\[
\partial_t u=\epsilon \triangle u, \qquad u(0, x_j)= \exp\bigl(-2\pi i\langle \xi_j, x_j \rangle\bigr).
\]
and so, 
\[
G_t^\epsilon \ast \exp\Bigl(-2\pi i\langle \xi_j, x_j \rangle \Bigr)=u(t,x_j).
\] 
This, together with \eqref{eq:aug26.2017.1} yields 
\begin{align}
k \; F^k_\xi[\sigma^\epsilon_t]   
&=  \exp\Bigl(-2\pi i \sum_{j=1}^k \langle \xi_j, \sqrt{2\beta} W_t \rangle \Bigr)  \prod_{j=1}^k \int_{\Mb} \exp\bigl(-4\pi^2 \epsilon t |\xi_j|^2 -2\pi i\langle \xi_j, x_j \rangle\bigr) m(dx_j) \nonumber\\
&=  \exp\Bigl(-2\pi i \sum_{j=1}^k \langle \xi_j, \sqrt{2\beta} W_t \rangle \Bigr)   \int_{\Mb^k} \exp\Bigl(-4\pi^2 \epsilon t \sum_{j=1}^k |\xi_j|^2 -2\pi i\Bigl\langle \sum_{j=1}^k \xi_j, x_j \Bigr\rangle\Bigr) m^{\otimes k}(dx) \nonumber\\
&=  \exp\bigl(-2\pi i \sum_{j=1}^k \langle \xi_j, \sqrt{2\beta} W_t \rangle \bigr)  \exp\bigl(-4\pi^2 \epsilon t \sum_{j=1}^k |\xi_j|^2 \bigr)
 \int_{\Mb^k} \exp\bigl( -2\pi i\bigl\langle \sum_{j=1}^k \xi_j, x_j \bigr\rangle\bigr) m^{\otimes k}(dx).\label{eq:aug26.2017.2}
\end{align}
We have 
\begin{align}
\mathbb E \Bigl( \exp\Bigl(-2\pi i \sum_{j=1}^k \langle \xi_j, \sqrt{2\beta} W_t \rangle \Bigr)\Bigr)   
&=\int_{\Mb} \exp\Bigl(-2\pi i \sum_{j=1}^k \langle \xi_j, \sqrt{2\beta} z \rangle \Bigr) G_t^{0.5}(z) dz
\nonumber\\
&=  \int_{\Mb} {1\over \sqrt{2\beta}^d } \exp\Bigl(-2\pi i  \Bigl\langle \sum_{j=1}^k \xi_j, w \Bigr\rangle \Bigr) G_t^{0.5}\Bigl({w \over \sqrt{2\beta}}\Bigr) dw\nonumber \\
&=  \int_{\Mb} \exp\Bigl(-2\pi i  \Bigl\langle \sum_{j=1}^k \xi_j, w \Bigr\rangle \Bigr) G_t^{\beta}( w) dw  .\label{eq:aug26.2017.3}
\end{align} 
We then face the computation of the Fourier transform of  $G_t^{\beta}$ which by the Fourier transform table (cf. e.g. \cite{LiebL} page 125)  is
\[
\widehat G_t^{\beta}(\sum_{j=1}^k \xi_j)= \exp\Bigl(-4\pi^2 \Big|\sum_{j=1}^k \xi_j\Big|^2 \beta t  \Bigr)
\]
This, together with  \eqref{eq:aug26.2017.3} yields 
\begin{equation}\label{eq:aug26.2017.4}
\mathbb E \Bigl( \exp\Bigl(-2\pi i \sum_{j=1}^k \langle \xi_j, \sqrt{2\beta} W_t \rangle \Bigr)\Bigr) = \exp\Bigl(-4\pi^2 \Big|\sum_{j=1}^k \xi_j\Big|^2 \beta t  \Bigr).
\end{equation}
Combining \eqref{eq:aug26.2017.2}, \eqref{eq:aug26.2017.3} and \eqref{eq:aug26.2017.4} we infer 
\[
 \mathbb E\Bigl( F^k_\xi[\sigma^\epsilon_t]  \Bigr)=  
\exp\Bigl(-\beta \lambda^2_{k, \frac{\epsilon}{\beta}}(\xi) t \Bigr)  F_\xi^k[m].
\]
This reads $U_\xi^k(t, m)= V_\xi^k(t, m),$ which concludes the proof.  \endproof

\subsection{Actions of Brownian motions on superpositions of eigenfunctions}\label{subsection:effet2} Throughout this subsection  $s, \delta>0$ , $\epsilon\geq 0$ and $\bigl( a_k\bigr)_{k=1}^\infty \in \mathcal A^s$. Assume \eqref{eq:uniform-bound-ak} holds and let $U_0$ be as in \eqref{eq:U-0}.    For $t>0$, set 
\[
V(t,m):=  \sum_{k=1}^\infty {1 \over k!} \int_{\Mb^k } b_k(t, \xi) F^k_\xi[m] d \xi
\] 
where for $t >0$ and $\epsilon \geq 0$, 
\begin{equation} \label{eq:agu26.2017.9bisb}
b_k(t, \xi):= a_k(\xi) \exp\Bigl(- \beta \lambda^2_{k, \frac{\epsilon}{\beta}}(\xi) t \Bigr) .
\end{equation} 
and $\lambda^2_{k,\epsilon}(\xi)$ is given by \eqref{eq:aug26.2017.6}. We will at some point need a stronger assumption than \eqref{eq:uniform-bound-ak}: 
\begin{equation}\label{eq:sep02.2017.1bis}
\int_{\Mb^k}|a_k(\xi)| d\xi \leq {C k! \over k^{3+\delta}}, 
\end{equation}

\begin{theorem}[superposition]\label{th:superposition} Assume $s>0$  and \eqref{eq:uniform-bound-ak}  holds.  Then the followings hold:
\begin{enumerate} 
\item[(i)] For any $l, t>0,$  the series $V(t, \cdot) \in H^{l}(\calPtwod)$ and converges uniformly.
\item[(ii)] We have $V(t, m)= \mathbb E\Bigl( U_0(\sigma^\epsilon_t[m])\Bigr)$ 
\item[(iii)]  If $\epsilon=0$, \eqref{eq:uniform-bound-ak}, \eqref{eq:sep01.2017.4}, \eqref{eq:sep01.2017.5} and \eqref{eq:sep01.2017.5bis}  hold, then  for $t>0,$ $V(t, \cdot)$ is twice continuously differentiable and 
\begin{equation}\label{eq:sep01.2017.8b}
 \beta \triangle_{w,\frac{\epsilon}{\beta}} V(t,m) =  \sum_{k=1}^\infty {1 \over k!}  \int_{\Mb^k } b_k(t, \xi) \beta \triangle_{w,\frac{\epsilon}{\beta}} F^k_\xi[m] d \xi. 
\end{equation} 
Furthermore, $V$ satisfies the heat equation \eqref{eq:aug26.2017.9}.
\item[(iv)]  If $\epsilon>0$ and \eqref{eq:sep02.2017.1bis} holds then $V(t, \cdot)$ is twice continuously differentiable, \eqref{eq:sep01.2017.8b} holds,  and $V$ satisfies the heat equation \eqref{eq:aug26.2017.9} 
\end{enumerate}  
\end{theorem} 
\proof{} (i)  
We  need to prove (i) only for $l>s.$ Fix $l>s$ and let $n \geq l$ be  integer and set 
\[
\bar C := \inf_{a \geq 0} {  (n)!+{(a t)^n   } \over  (n)!(1+a)^n}>0.
\]  
We have 
\[
|b_k(t, \xi)|^2\Bigl( 1+\lambda^2_{k,0}(\xi)\Bigr)^l 
\leq | a_k(\xi)|^2{\Bigl( 1+\lambda^2_{k,0}(\xi)\Bigr)^l \over \Bigl( 1+{\left( \beta \lambda^2_{k,\frac{\epsilon}{\beta}}(\xi) t\right)^n \over (n)!   }\Bigr)}  \leq | a_k(\xi)|^2{\Bigl( 1+\lambda^2_{k,0}(\xi)\Bigr)^l \over \Bigl( 1+{(\lambda^2_{k,0}(\xi) t)^n \over (n)!   }\Bigr)} 
\]
and so, 
\begin{equation} \label{eq:agu26.2017.9}
|b_k(t, \xi)|^2\Bigl( 1+\lambda^2_{k,0}(\xi)\Bigr)^l \leq   {| a_k(\xi)|^2\over \bar C}.
\end{equation} 
Thus, 
\[
\sum_{k=1}^\infty {1 \over k!} \int_{\Mb^k } |b_k(t, \xi)|^2\Bigl( 1+\lambda^2_{k,0}(\xi)\Bigr)^l d \xi   \leq  
\sum_{k=1}^\infty {1 \over k!} \int_{\Mb^k } {| a_k(\xi)|^2\over \bar C} d \xi <\infty.
\] 
Similarly, 
\[
\sum_{k=N+1}^\infty {1 \over k!}\biggl| \int_{\Mb^k } a_k(\xi) \exp\bigl(- \beta \lambda^2_{k,\frac{\epsilon}{\beta}}(\xi) t \bigr)  F^k_\xi[m] d \xi \biggr| \leq  
\sum_{k=N+1}^\infty {1 \over k!}  \int_{\Mb^k } |a_k(\xi)|   d \xi \leq \sum_{k=N+1}^\infty {C \over k^{1+\delta}}.
\] 
These prove (i). 

(ii) By Lemma \ref{le:heat-special-case}, given an integer $N>1$ we have 
\[
\mathbb E \biggl( \sum_{k=1}^N {1 \over k!} a_k(\xi)   F_\xi^k\bigl[\sigma^\epsilon_t[m] \bigr] \biggr)= 
\sum_{k=1}^N {1 \over k!} b_k(t, \xi) F_\xi^k[m].
\] 
and so, 
\begin{equation}\label{eq:aug26.2017.7}
\mathbb E \biggl( \int_{\Mb^k }\sum_{k=1}^N {1 \over k!} a_k(\xi)   F_\xi^k\bigl[\sigma^\epsilon_t[m] \bigr]  d\xi \biggr) = 
\sum_{k=1}^N {1 \over k!} \int_{\Mb^k } b_k(t, \xi) F_\xi^k[m] d\xi.
\end{equation} 
By (i), the expression on the right hand side of \eqref{eq:aug26.2017.7} converges uniformy to $V(t, m).$ Since $|F_\xi^k| \leq 1$ 
\[
\sum_{k=1}^N {1 \over k!} \bigg|\int_{\Mb^k} a_k(\xi)   F_\xi^k\bigl[\sigma^\epsilon_t[m] \bigr] d\xi \bigg| \leq  \sum_{k=1}^N {1 \over k!} \int_{\Mb^k} |a_k(\xi)| d\xi.
\] 
We apply the Lebesgue dominated convergence theorem to obtain 
\[
\mathbb E \biggl( \int_{\Mb^k }\sum_{k=1}^\infty {1 \over k!} a_k(\xi)   F_\xi^k\bigl[\sigma^\epsilon_t[m] \bigr]  d\xi \biggr)=
\sum_{k=1}^\infty \mathbb E \biggl( \int_{\Mb^k } {1 \over k!} a_k(\xi)   F_\xi^k\bigl[\sigma^\epsilon_t[m] \bigr] d\xi \biggr).
\] 
In conclusion, we have proven that letting $N$ tend to $\infty$ in \eqref{eq:aug26.2017.7} yields 
\[
\mathbb E\Bigl( U_0(\sigma^\epsilon_t[m])\Bigr)=V(t,m).
\]

Under the sole assumptions in \eqref{eq:uniform-bound-ak}, \eqref{eq:sep01.2017.4}, 
\[
 \Big| \beta \lambda^2_{k,\frac{\epsilon}{\beta}}(\xi)|a_k(\xi)| \exp\Bigl(-\beta \lambda^2_{k,\frac{\epsilon}{\beta}}(\xi) t \Bigr)  F^k_\xi[m] \Big| \leq 
  {\beta \lambda^2_{k,\frac{\epsilon}{\beta}}(\xi) \over \exp\Bigl(\beta \lambda^2_{k,\frac{\epsilon}{\beta}}(\xi) t \Bigr)} |a_k(\xi)| \leq { |a_k(\xi)| \over k t}. 
\]
Thus, 
\[
 \sum_{k=1}^\infty {1 \over k!} \Big| \int_{\Mb^k } \beta \lambda^2_{k,\frac{\epsilon}{\beta}}(\xi) a_k(\xi) \exp\Bigl(- \beta \lambda^2_{k,\frac{\epsilon}{\beta}}(\xi) t \Bigr)  F^k_\xi[m] d \xi\Big| \leq \sum_{k=1}^\infty {C \over t\; k^{1+\delta}}.
\] 
This proves the uniform convergence of the following series:   
\begin{equation}\label{eq:sep01.2017.6late}
 \partial_t V(t,m)=- \sum_{k=1}^\infty {1 \over k!}  \int_{\Mb^k } \beta \lambda^2_{k,\frac{\epsilon}{\beta}}(\xi) a_k(\xi) \exp\Bigl(- \beta \lambda^2_{k,\frac{\epsilon}{\beta} }(\xi) t \Bigr)  F^k_\xi[m] d \xi.
\end{equation} 
Since $F^k_\xi[m] $ is an eigenfunction of $\beta \triangle_{w,\frac{\epsilon}{\beta}}$ and $- \beta \lambda^2_{k,\frac{\epsilon}{\beta}}(\xi)$ is the associate eigenvalue, \eqref{eq:sep01.2017.6late} reads off 
\begin{equation}\label{eq:sep01.2017.6}
 \partial_t V(t,m)= \sum_{k=1}^\infty {1 \over k!}  \int_{\Mb^k } b_k(t, \xi)  \beta \triangle_{w,\frac{\epsilon}{\beta}} F^k_\xi[m] d \xi.
\end{equation} 

(iii) Suppose $\epsilon=0$, \eqref{eq:uniform-bound-ak}, \eqref{eq:sep01.2017.4}, \eqref{eq:sep01.2017.5} and \eqref{eq:sep01.2017.5bis}  hold.  The identity $|b_k| \leq |a_k|$ yields   
\begin{equation}\label{eq:sep01.2017.7} 
\int_{\Mb^k} |b_k(t, \xi)| \leq {C k! \over k^{\delta} }, \quad 
\int_{\Mb^k} |b_k(t, \xi)| \bigl(|\xi_1|+|\xi_1|^2  \bigr)d\xi \leq {2C k! \over k^{1+\delta} },
\end{equation}
\begin{equation}\label{eq:sep01.2017.8}
\int_{\Mb^k} |b_k(t, \xi)| \cdot | \xi_1|\; | \xi_2| d\xi  \quad \leq {C k! \over k^{2+\delta} } 
\end{equation}
and 
\begin{equation}\label{eq:sep01.2017.8ter}
 \int_{\Mb^k} |b_k(\xi)| \cdot |\xi_1|^3 d\xi \leq {Ck! \over k^{1+\delta} }, \quad   \int_{\Mb^k} |b_k(\xi)| \cdot | \xi_1|^2 \cdot  | \xi_2| d\xi  \quad \leq {Ck! \over k^{3+\delta} }.
\end{equation}
Thanks to Corollary \ref{co:superposition2}, \eqref{eq:sep01.2017.7} , \eqref{eq:sep01.2017.8} and \eqref{eq:sep01.2017.8ter} yield 
\[
 \sum_{k=1}^\infty {1 \over k!}  \int_{\Mb^k } b_k(t, \xi)  \beta \triangle_{w,\frac{\epsilon}{\beta}} F^k_\xi[m] d \xi= 
  \beta \triangle_{w,\frac{\epsilon}{\beta}} \biggl(\sum_{k=1}^\infty {1 \over k!}  \int_{\Mb^k } b_k(t, \xi)  F^k_\xi[m] d \xi\biggr).
\]
This, together with \eqref{eq:sep01.2017.6} completes the proof of (iii).

(iv) Suppose $\epsilon>0$ and \eqref{eq:sep02.2017.1bis} holds. We have 
\begin{equation}\label{eq:sep30.2017.9}
\Biggl(  1+\Bigl(4 \pi^2 \epsilon t \sum_{j=1}^k |\xi_j|^2\Bigr) +{1\over 2} \Bigl(4 \pi^2 \epsilon t \sum_{j=1}^k |\xi_j|^2\Bigr)^2 \Biggr) |b_k(t, \xi)| \leq   |a_k(\xi)|. 
\end{equation}
This, together with \eqref{eq:sep02.2017.1bis}, yields 
\begin{equation}\label{eq:sep01.2017.9} 
 \min\{ 1,  4 \pi^2 \epsilon t \}\int_{\Mb^k} |b_k(t, \xi)| \bigl( 1+ |\xi_1|+|\xi_1|^2  \bigr)d\xi  \leq 3  \int_{\Mb^k} |a_k(\xi)| d\xi \leq {3 C k! \over k^{3+\delta} }
\end{equation}
Similarly, using the fact that $2|\xi_1| |\xi_2| \leq |\xi_1|^2+|\xi_2|^2$ we obtain 
\[
8 \pi^2 \epsilon t \int_{\Mb^k} |b_k(t, \xi)|\cdot  | \xi_1| \cdot  | \xi_2| d\xi \leq 4 \pi^2 \epsilon t \int_{\Mb^k} |b_k(t, \xi)|( | \xi_1|^2+ | \xi_2|^2) d\xi 
\]
and so, 
\begin{equation}\label{eq:sep01.2017.10}
8 \pi^2 \epsilon t \int_{\Mb^k} |b_k(t, \xi)| \cdot | \xi_1|\cdot  | \xi_2| d\xi \leq   \int_{\Mb^k} |a_k(\xi)| d\xi \leq  
 {C k! \over k^{3+\delta} }.
\end{equation} 
We exploit the crude estimate 
\begin{equation}\label{eq:sep30.2017.10}
 |\xi_1|^2 |\xi_2| +|\xi_1|^3  \leq  {3 \over 2}\Bigl( 1+  |\xi_1|^4 + |\xi_2|^2 \Bigr) 
\end{equation} 
But 
\begin{equation}\label{eq:sep30.2017.11}
\min\{1, 4 \pi^2 \epsilon t , 8 \pi^4 \epsilon^2 t^2 \}\Bigl( 1+  |\xi_1|^4 + |\xi_2|^2 \Bigr) \leq \Biggl(  1+\Bigl(4 \pi^2 \epsilon t \sum_{j=1}^k |\xi_j|^2\Bigr) +{1\over 2} \Bigl(4 \pi^2 \epsilon t \sum_{j=1}^k |\xi_j|^2\Bigr)^2 \Biggr)
\end{equation} 
Thanks to \eqref{eq:sep02.2017.1bis} and \eqref{eq:sep30.2017.9}, we combine \eqref{eq:sep30.2017.10} and \eqref{eq:sep30.2017.11} to obtain a constant $\bar \epsilon(t) \in (0,\infty)$ depending only on $\epsilon$ and $t>0$ such that 
\begin{equation}\label{eq:sep01.2017.12}
\int_{\Mb^k} |b_k(t, \xi)| \cdot | \xi_1|^3 \cdot   | \xi_1|^2 | \xi_2| d\xi \leq      {\bar \epsilon(t) \; k! \over k^{3+\delta} }.
\end{equation}  
Thanks to \eqref{eq:sep01.2017.9}, \eqref{eq:sep01.2017.10} and \eqref{eq:sep01.2017.12}, we may apply Corollary \ref{co:superposition2} to obtain that $V(t, \cdot)$ is twice continuously differentiable and the identity in \eqref{eq:sep01.2017.8b} holds.  This, together with  \eqref{eq:sep01.2017.6}, proves (iv). \endproof

\begin{remark}[smoothing effects] \label{re:superposition}  Suppose $\epsilon, s>0$,  and let $\bigl( a_k\bigr)_{k=1}^\infty \in \mathcal A^s$ be such that for any $k\geq 1$, $a_k \in L^2(\Mb^k; \mathbb C)$ and \eqref{eq:sep02.2017.1bis} holds. Let $U_0$ and $V(t, \cdot)$ be the uniformly convergent series in Theorem \ref{th:superposition}. Then for any $t>0$, $\left(  \beta \triangle_{w,\frac{\epsilon}{\beta}} \right)^r V(t, \cdot)$ is twice continuously  differentiable for all $r \in \mathbb{N}$. 
\end{remark}
\proof{} The statement is clear for $r = 0$ from Theorem \ref{th:superposition}. We start with $r = 1$. Recall that $b_k$ is given by \eqref{eq:agu26.2017.9bisb}. Set
\[
c_k(t, \xi):= - \beta \lambda^2_{k,  \frac{\epsilon}{\beta} }(\xi) b_k(t, \xi)= - {a_k(\xi) \beta \lambda^2_{k, \frac{\epsilon}{\beta}}(\xi)  \over \exp\bigl( \beta \lambda^2_{k, \frac{\epsilon}{\beta}}(\xi) t \bigr) }
\]
By Theorem \ref{th:superposition}  
\[
\beta \triangle_{w,\frac{\epsilon}{\beta}} V(t, m):= \sum_{k=1}^\infty {1 \over k!} \int_{\Mb^k } c_k(t, \xi)  F^k_\xi[m] d \xi.
\] 
We have 
\begin{equation}\label{eq:sep01.2017.11}
|c_k(t, \xi)| \leq {|a_k(\xi)| \over t}, 
\end{equation}
and setting $|\xi|^2= \sum_{j=1}^k|\xi_j|^2$ we obtain 
\[
 |c_k(t, \xi)|\cdot |\xi|^4 =
 {|a_k(\xi)| \cdot  \beta \lambda^2_{k, \frac{\epsilon}{\beta}} (\xi)\;  |\xi|^4 \over \exp\Bigl( { \beta \lambda^2_{k, \frac{\epsilon}{\beta}} (\xi) t \over 2} \Bigr) \cdot \exp\Bigl( { \beta \lambda^2_{k, \frac{\epsilon}{\beta}} (\xi) t \over 2} \Bigr) } \leq 
{|a_k(\xi)| \cdot \beta \lambda^2_{k, \frac{\epsilon}{\beta}} (\xi) |\xi|^4  \over  \exp\Bigl( {  \beta \lambda^2_{k, \frac{\epsilon}{\beta}} (\xi)t \over 2} \Bigr)\cdot \exp\Bigl(      2 \pi^2 \epsilon t |\xi|^2  \Bigr) } 
\]
Thus, 
\begin{equation}\label{eq:sep01.2017.12new}
 |c_k(t, \xi)| \cdot |\xi|^4 \leq {|a_k(\xi)| \over \pi^4 t^3 \epsilon^2}
\end{equation} 
We use the crude estimate 
\[
1+ |\xi_1|+  |\xi_1|^2 +  |\xi_1|  |\xi_2| + |\xi_1|^2  |\xi_2| +  |\xi_1|^3 \leq 6(1+ |\xi|^4).
\]
to conclude that 
\[
\int_{\Mb^k} |c_k(t, \xi)| \Bigl( 1+ |\xi_1|+  |\xi_1|^2 +  |\xi_1|  |\xi_2| + |\xi_1|^2  |\xi_2| +  |\xi_1|^3 \Bigr) d\xi \leq 6\int_{\Mb^k} |c_k(t, \xi)| (1+ |\xi|^4) d\xi
\]
This, together with \eqref{eq:sep01.2017.12new} and \eqref{eq:sep01.2017.11} implies there is a constant $\bar \epsilon(t) \in (0,\infty)$ depending only on $\epsilon$ and $t>0$ such that 
\begin{equation}\label{eq:sep30.2017.12}
\int_{\Mb^k} |c_k(t, \xi)| \Bigl( 1+ |\xi_1|+  |\xi_1|^2 +  |\xi_1|  |\xi_2| + |\xi_1|^2  |\xi_2| +  |\xi_1|^3 \Bigr) d\xi \leq {\epsilon(t) k! \over k^{3+\delta}}. 
\end{equation}

Thanks to Corollary \ref{co:superposition2},  \eqref{eq:sep30.2017.12} implies $\beta \triangle_{w, \frac{\epsilon}{\beta}} V(t, \cdot)$ is twice differentiable. 

To arrive at the statement with $r >1$, one only need to inductively follow the same argument as above: apply Corollary \ref{co:superposition2} together with noticing $2|\xi_1| |\xi_2| \leq |\xi_1|^2+|\xi_2|^2$ and 
\[
 {|a_k(\xi)| \cdot  \left( \beta \lambda^2_{k, \frac{\epsilon}{\beta}} (\xi) \right)^r \;  |\xi|^4  \over \exp\Bigl( { \beta \lambda^2_{k, \frac{\epsilon}{\beta}} (\xi) t } \Bigr)  } \leq 
{|a_k(\xi)| \cdot \left( \beta \lambda^2_{k, \frac{\epsilon}{\beta}} (\xi) \right)^r   |\xi|^4  \over  \left( 1+  \frac{1}{r!} \left( {  \beta \lambda^2_{k, \frac{\epsilon}{\beta}} (\xi)t \over 2} \right)^r \right)\cdot \exp\Bigl(      2 \pi^2 \epsilon t |\xi|^2  \Bigr) } 
\]\endproof

\section{Measures concentrated on locally compact subset of $\calPtwod \times \calPtwod.$} 
The set $\mathcal M^2:= \calPtwod \times \calPtwod$ is a metric space when endowed with the metric $\mathbb W$ defined by 
\[
\mathbb W^2(m, \tilde m):= W_2^2(m_1, \tilde m_1)+ W_2^2(m_2, \tilde m_2) 
\] 
for $m=(m_1, m_2),$ $\tilde m=(\tilde m_1, \tilde m_2) \in \calPtwod \times \calPtwod.$ If $R>0$ we denote as $B_R$ the ball in $\Mb$ centered at the origin and of radius $R$. Note the infinite dimensional ball in $\mathcal M^2$, of radius $R>0$, centered at $\delta_{(0,0)}$ is  
\[
\bigl\{(m_1, m_2) \in \mathcal M^2 \; | \; \int_\Mb |x|^2 (m_1+m_2)(dx)\leq R \bigr\}.
\]
Let $\mathbb S_R$ denote the set of $m \in \calPtwod$ such that the support of $m$ is contained in $B_R$ and let $\chi_R: \calPtwod \rightarrow \{0,1\}$ be the function which assumes the value $1$ on  $\mathbb S_R$ and the value $0$ on the complement of $\mathbb S_R$. 

Fix a natural number $k,$ If $r, p \leq k$ are natural numbers and let $I \in C_r^k$ and $J \in C_p^k.$ The map $x \rightarrow (m_{x_I}, m_{x_J})$ is a $(r^{-1}+p^{-1})$--Lipschitz map of $\Mb$ onto $\calPtwod \times \calPtwod$. Thus, it  maps compact sets into compact sets and bounded sets into bounded set. 

Let $\mathbb P^{I, J, R}$ denote the push forward of the Lebesgue measure on $B_R^k$ by the map $x \rightarrow (m_{x_I}, m_{x_j})$. If $H: \mathcal M^2 \rightarrow [0,\infty)$ is continuous then 
\[
\int_{\mathcal M^2} H(m_1, m_2) \mathbb P^{I ,J, R}(dm_1, dm_2)= \int_{B_R^k} H(m_{x_I}, m_{x_J}) dx_1 \cdots dx_k.
\] 
We define the signed Borel regular measure 
\[
\mathbb P^{k, R}:={k^2 \over (k!)^3} \sum_{r=1}^k \sum_{p=1}^k (-1)^{r+p} r^k p^k \sum_{I \in C^k_r} \sum_{J \in C^k_p} \mathbb P^{I, J, R}.
\]

\begin{proposition}\label{pr:measureR0}  Let $r, p \leq k$ be integers and let  $I \in C_r^k$, $J \in C_p^k$ be multi--indexes.  
\begin{enumerate}
\item[(i)] The measure  $\mathbb P^{I ,J, R}$ are of finite total mass and are supported by $\mathbb S_R$. 
\item[(ii)] The signed measure is  $\mathbb P^{k, R}$ of finite total variations and is supported by $\mathbb S_R$. 
\end{enumerate}  
\end{proposition} 
\proof{} Since (i) implies (ii), it suffices to show (i).   Setting $H\equiv 1$ on $\mathcal M^2$ we have 
\begin{equation}\label{eq:oct06.2017.3}
 \mathbb P^{I ,J, R}\bigl( \mathcal M^2\bigr)=\mathcal L^{kd}(B_R^k)<\infty.
\end{equation} 
Observe that if $x \in B_R^k$ then $m_{x_I}, m_{x_J} \in \mathbb S_R$ and so, $\chi_R(m_{x_I}) \chi_R(m_{x_J})=1.$ Thus, 
\[
\int_{\mathcal M^2} \chi_R(m_1) \chi_R(m_2)  \mathbb P^{I ,J, R}(dm_1, dm_2)= \int_{B_R^k}  \chi_R(m_{x_I}) \chi_R(m_{x_J})dx_1 \cdots dx_k=\mathcal L^{kd}(B_R^k).
\]  
This, together with \eqref{eq:oct06.2017.3} proves (i). \endproof

\begin{proposition}\label{pr:measureR}  If  $\Phi, \Psi \in L^2(\Mb^k/P_k)$ are supported by $B_R$  then  
$$
\Bigl\langle F_{\Phi}; F_{\Psi } \Bigr\rangle_{H^0 }=\int_{\mathcal M^2} F_\Phi[m_1] F_\Psi[m_2] \mathbb P^{k ,R}(dm_1, dm_2)
$$
\end{proposition} 
\proof{} By Lemma \ref{le:Hs}, since $\Phi$ and $\Psi$ are supported by $B_R$ we have 
\[
\langle F_\Phi; F_\Psi \rangle_{H^0 }={1 \over k!} \int_{B_R^k} \Phi(x) \Psi(x) dx.
\] 
By Theorem \ref{th:inverseOk} we obtain 
\begin{align}
\langle F_\Phi; F_\Psi \rangle_{H^0}  
&={k^2 \over (k!)^3} \int_{B_R^k} 
\sum_{r=1}^k   \Bigl( (-1)^{k-r} r^k  \sum_{I \in C^k_r} F_{\Phi }(m_{x_I})  \Bigr) 
   \sum_{s=1}^k   \Bigl( (-1)^{k-s} s^k  \sum_{J \in C^k_s} F_{\Psi }(m_{x_J})  \Bigr)  dx \nonumber\\
&= \sum_{r, s=1}^k {k^2 \over (k!)^3}  (-1)^{r+s} r^k s^k  \sum_{I \in C^k_r}  \sum_{J \in C^k_s}  \int_{\mathcal M^2} F_{\Phi }(m_1) F_{\Psi }(m_2)  \mathbb P^{I ,J}(dm_1, dm_2).  \nonumber
\end{align}
This concludes the proof. \endproof

\begin{theorem}\label{th:measureR}  If $R>0$ and Let $\Phi, \Psi \in C_c^3(\Mb^k/P_k)$ then 
\[
-\int_{\mathcal M^2} \triangle_w F_\Phi[m_1] F_\Psi[m_2]  d\mathbb P^{k ,R} =\int_{\mathcal M^2} D_2(\nabla_w F_\Phi, \nabla_w G_\Psi)  d\mathbb P^{k ,R} 
\]
where 
\[
D_2(\nabla_w F, \nabla_w G)(m_1, m_2):=  \int_{\Mb^2}  \langle \nabla_w F[m_1](q_1)  ;  \nabla_w G[m_2](q_2) \rangle m_1(dq_1) m_2(dq_2) . 
\]
\end{theorem} 
\proof{} Let $m \in \mathbb S_R.$ Set 
\[
 f:= \sum_{j=1}^k \nabla_{x_j} \Phi =\bigl[f^1, \cdots, f^d \bigr]^T, \qquad  g:= \sum_{j=1}^k \nabla_{x_j} \Psi =\bigl[g^1, \cdots, g^d \bigr]^T
\]
The material presented in  Subsection \ref{subsec:particular} allows to obtain  
\begin{align}\label{eq:oct07.2017.1}
    \int_{\Mb} \nabla_w F_\Phi[m](x_1)m(dx_1)= {1\over k}\int_{\Mb^k} f(x_1, \cdots, x_k) m(dx_1) \cdots m(dx_k) 
    &= \begin{bmatrix}
           F_{f^1}[m] \\
           \vdots \\
          F_{f^d}[m]
         \end{bmatrix}
  \end{align}
Similarly, 
\begin{align}\label{eq:oct07.2017.2}
    \int_{\Mb} \nabla_w F_\Psi[m](x_1)m(dx_1)= {1\over k}\int_{\Mb^k} g(x_1, \cdots, x_k) m(dx_1) \cdots m(dx_k) 
    &= \begin{bmatrix}
           F_{g^1}[m] \\
           \vdots \\
          F_{g^d}[m]
         \end{bmatrix}
  \end{align}
By Proposition \ref{pr:measureR} 
\begin{eqnarray*}
\left \langle \int_{\Mb} \nabla_w F_\Phi (m,x) m(dx); \int_{\Mb} \nabla_w F_\Psi (m,x) m(dx) \right \rangle_{H^0}
&=& \sum_{n=1}^d \langle F_{f^n};  F_{g^n} \rangle_{H^0} \\
&=& \sum_{n=1}^d \int_{\mathcal M^2} F_{f^n} [m_1] F_{g^n}[m_2] \mathbb P^{k ,R}(dm_1, dm_2). 
\end{eqnarray*}
In other words,  
\begin{equation}\label{eq:soct07.2017.4}
\left \langle \int_{\Mb} \nabla_w F_\Phi (m,x) m(dx); \int_{\Mb} \nabla_w F_\Psi (m,x) m(dx) \right \rangle_{H^0} =  \int_{\mathcal M^2} D_2(\nabla_w F, \nabla_w G)(m_1, m_2)  \mathbb P^{k ,R}(dm_1, dm_2)
\end{equation}
By  Corollary \ref{co:grad-hess1} (ii), $\triangle_w F_\Phi= F_\Theta$ where 
\[
\Theta(x_1, \cdots, x_k)= \sum_{j, l=1}^k \sum_{n=1}^d {\partial^2 \Phi \over \partial (x_j)_n \partial (x_l)_n}.
\] 
This, together with Proposition \ref{pr:measureR}  implies 
\[
\langle \triangle_w F_\Phi; F_\Psi \rangle_{H^0}=\langle F_\Theta; F_\Psi \rangle_{H^0}=  \int_{\mathcal M^2}  F_\Theta[m_1]  F_\Psi[m_2] \mathbb P^{k ,R}(dm_1, dm_2).
\] 
Hence, 
\begin{equation}\label{eq:soct07.2017.5}
\langle \triangle_w F_\Phi; F_\Psi \rangle_{H^0}= \int_{\mathcal M^2}  \triangle_w F_\Phi[m_1]  F_\Psi[m_2] \mathbb P^{k ,R}(dm_1, dm_2)
\end{equation} 
Thanks to Proposition \ref{pr:integration-by-parts}, \eqref{eq:soct07.2017.4} and \eqref{eq:soct07.2017.5} yield the desired result.\endproof

\section*{Acknowledgements}
\noindent The research of WG was supported by NSF grant DMS--17 00 202 and the research of YTC was supported by NSF grant ECCS-1462398.
The authors wish to thank  R. Hynd, J. Lott and T. Pacini for fruitful comments or criticisms on an earlier draft of this manuscript.  The author would like to express their gratitude to L. Ambrosio and P. Cardaliaguet  for stimulating discussions that helped improve the manuscript.

%


\begin{thebibliography}{99}

\bibitem{AmbGb}{\sc L. Ambrosio, W. Gangbo}, \emph{Hamiltonian ODE's in the Wasserstein space of probability measures}, Comm. Pure Appl. Math. \textbf{61} (2008), No. 2, pp. 18--53.

\bibitem{AGS}
{\sc L.~Ambrosio, N. Gigli, G. Savar\'e},
\newblock Gradient flows in metric spaces and the Wasserstein spaces of
probability measures.
\newblock {\em Lectures in Mathematics}, ETH Zurich, Birkh\"auser, 2005.

\bibitem{Brenier}{\sc Y. Brenier}, \emph{Polar factorization and monotone rearrangement of vector-valued functions}, Comm. Pure Appl. Math. \textbf{XLIV} (1991), 375--417.

\bibitem{BrenierG2003} {\sc Y. Brenier, W. Gangbo},  \emph{$L\sp p$ approximation of maps by diffeomorphisms},  Calc. Var. Partial Differential Equations  {\bf 16},  no. 2 (2003), 
147--164.

\bibitem{BorweinZhu}
{\sc J.M. Borwein, Q.J. Zhu}, \newblock \emph{A survey of subdifferential calculus with applications}, Nonlin. Anal. {\bf 38} (1999), 687--773.


\bibitem{Caffarelli}
{\sc L. Caffarelli}, \emph{Boundary regularity of maps with convex potentials},  Ann. Math. (Second Series), \textbf{144}, No. 3 (1996), pp. 453--496.

\bibitem{Cardaliaguet} {\sc P. Cardaliaguet}, \emph{Notes on Mean-Field Games}, lectures by P.L. Lions, Coll\`{e}ge de France, 2010.

\bibitem{CardaliaguetDLL} {\sc P. Cardaliaguet, F. Delarue, J-M. Lasry, P-L. Lions}, \emph{The master equation and the convergence problem in mean field games}, (Preprint).

\bibitem{CarmornaDe} {\sc R. Carmona, F. Delarue}, \emph{The master equation for large population in equilibriums}, (Preprint).

\bibitem{CarmornaDeI} Probabilistic Theory of Mean Field Games with Applications I, {\em Probability Theory and Stochastic Modelling} , Springer.

\bibitem{CarmornaDeII} Probabilistic Theory of Mean Field Games with Applications I, {\em Probability Theory and Stochastic Modelling} , Springer.

\bibitem{CrandallLions} {\sc M.G. Crandall, P.-L. Lions}, \emph{Hamilton-Jacobi Equations in Infinite Dimensions I. Uniqueness of Viscosity Solutions}, J. Funct. Anal. \textbf{62} (1985), 379--396.

\bibitem{CL2}{\sc M.G. Crandall, P.-L. Lions}, \emph{Hamilton-Jacobi equations in infinite dimensions II. Existence of viscosity solutions}, J. Funct. Anal. \textbf{65} (1986), 368--405.

\bibitem{CL3}{\sc M.G. Crandall, P.-L. Lions}, \emph{Hamilton-Jacobi equations in infinite dimensions III}, J. Funct. Anal. \textbf{68} (1986), 214--247.

\bibitem{Dorfman}{\sc I. Ja. Dorfman}, \emph{On means and the Laplacian of functions on Hilbert spaces}, Math.  USSR  Sbornik, Tom  81  (123)  (1970),  No. 2.                                                                                                              


\bibitem{GangboKP} {\sc  W. Gangbo, H.K. Kim, T. Pacini},  \emph{Differential forms on Wasserstein space and infinite-dimensional Hamiltonian systems}, Memoirs of the AMS, vol 211 (2011), no 993 (3 of 5)  1--77.

\bibitem{GNT2}{\sc W. Gangbo, T. Nguyen, A. Tudorascu}, \emph{Hamilton-Jacobi equations in the Wasserstein space},
Meth. Appl. Anal. \textbf{15}, No. 2 (2008), pp. 155--184.


\bibitem{GangboS2015} {\sc  W. Gangbo, A. Swiech},  \emph{Metric viscosity solutions of Hamilton--Jacobi equations depending on local slopes}, Calc. Var. (2015) 54:1183--1218.

\bibitem{GangboT2017}{\sc W. Gangbo, A. Tudorascu}, \emph{On differentiability in the Wasserstein space and well--posedness for Hamilton--Jacobi equations},
(in progress).


\bibitem{Gigli2008} {\sc  N. Gigli},  \emph{Second order analysis on $(\mathcal P_2(M), W_2)$}, Memoirs of the AMS, vol 216 (2012), no 1018 (End of volume).

\bibitem{Greenberg} {\sc   M. J. Greenberg}, \emph{Lectures on forms in many variables}, Mathematics Lecture Note Series, W. A. Benjamin, Inc., New York, 1969. 

\bibitem{PLLions2007-2008}  {\sc   P--L. Lions}, \emph{Lectures Coll\`ege de France}, 2007--2008.

\bibitem{PLLions2008-2009}  {\sc   P--L. Lions}, \emph{Lectures Coll\`ege de France}, 2008--2009.

\bibitem{Lott2008}  {\sc  J. Lott},  \emph{Some Geometric Calculations on Wasserstein Space}, Commun. Math. Phys. (2008) 277, 423--437.

\bibitem{HyndK2} {\sc  R. Hynd, H.K. Kim},  \emph{Value functions in the Wasserstein spaces: finite time horizons}, J. Funct. Anal. \textbf{269} (2015), 968--997. 

\bibitem{Levy}  {\sc P. L\'evy}, \emph{Probl\`emes concrets d'analyse fonctionnelle}, avec compl\'ements sur les fonctionnelles analytiques par F. Pellegrino,  Collection de monographies sur la th�orie des fonctions, Paris, Gauthier-Villars, 1951. 

\bibitem{LiebL}  {\sc E. H. Lieb, M. Loss}, \emph{Analysis}, Graduate Studies in Mathematics, 2nd American Mathematical Soc., 2001.

\bibitem{MazurO}  {\sc S. Mazur and W. Orlicz}, \emph{Grundlegende Eigenschaften der polynomischen Operationen}, Studia Math.\textbf{5} (1934) 50--��68. 

\bibitem{Nelson}  {\sc E. Nelson}, \emph{Probability theory and Euclidean field theory. In: G. Velo and A. Wightman (Eds.), Constructive Quantum Field Theory}. Lecture Notes
in Physics 25, pp. 94--��124. Springer-Verlag, Berlin-New York, 1973. 

\bibitem{Schetzen} {\sc M. Schetzen}, \emph{The Volterra and Wiener Theories of Nonlinear Systems}, Wiley, New York, 1980. 

\bibitem{Thomas} {\sc E. G. F. Thomas}, \emph{A polarization identity for multilinear maps}, Indagationes Mathematicae, Volume \textbf{25}, Issue 3, 18  (2014), 468--474. 


\bibitem{Preiss}{\sc D. Preiss}, \emph{Differentiability of Lipschitz functions on Banach spaces}, J. Funct. Anal. \textbf{91} (1990), 312--345.

\bibitem{SteinS}  {\sc E. M. Stein, R. Shakarchi}, \emph{Real Analysis, Measure Theory, Integration, and Hilbert Spaces}, Princeton Lectures in Analysis III, Princeton University Press2001.

\bibitem{VonRenesseS} {\sc  M-K. von Renesse, K-T. Sturm},  \emph{Entropic measure and Wasserstein diffusion}, Ann. Probab. Volume 37, Number 3 (2009), 1114--1191.

\bibitem{Sturm} {\sc K-T. Sturm}, Entropic Measure on Multidimensional Spaces, Seminar on Stochastic Analysis, Random Fields and Applications VI, 261--277, 2011.

\bibitem{Umemura}{\sc Y. Umemura}, \emph{On the infinite dimensional Laplacian operator}, J. Math., Kyoto University 4.3  (1965), 477--492. 

\bibitem{Ward} {\sc H. N. Ward}, \emph{Combinatorial polarization}, Discrete Math. 26 (1979), Volume 26, Issue 2, Pages 185--197.
\end{thebibliography}
\end{document}